\documentclass[10pt]{amsart}
\usepackage{amsmath}
\usepackage{times}
\usepackage{mathabx} 
\usepackage{amsfonts}
\usepackage{hyperref}
\usepackage{amsthm}
 \usepackage{amssymb} 
 \usepackage[all]{xy}
 \usepackage{textcomp}
\usepackage{mathrsfs}

\makeindex
\usepackage{latexsym}
\usepackage[OT2,T1]{fontenc}
\DeclareSymbolFont{cyrletters}{OT2}{wncyr}{m}{n}
\DeclareFontFamily{OT1}{rsfs}{}
\usepackage{color}

\newif\iflongproof
\longprooftrue

\newcommand{\changecolor}{black}

\usepackage{verbatim}
\usepackage{graphicx}
\usepackage{nomencl}
\DeclareFontEncoding{OT2}{}{} %
  
     \def\acts{\mathrel{\reflectbox{$\righttoleftarrow$}}}
\newcommand{\Spec}{\mathrm{Spec}}
\newcommand{\Z}{\mathbf{Z}}
\renewcommand{\indexname}{Index of symbols and important notation}
\newcommand{\Afinite}{\mathbf{A}_{\mathrm{f}}}

\newcommand{\Jinf}{\mathsf{J}}
\newcommand{\usbR}{\overline{\mathrm{R}}}

\newcommand{\mot}{\mathrm{mot}}
\newcommand{\Gal}{\mathrm{Gal}}
\newcommand{\coad}{\mathrm{coad}}
\newcommand{\LG}{^{L}\hat{G}}

\newcommand{\defring}{\mathrm{R}_{\rhobar}}
\newcommand{\Ext}{\mathrm{Ext}}
\newcommand{\univ}{\mathrm{univ}}
\newcommand{\coproduct}{\mathrm{coprod}}
\newcommand{\Frobl}{\mathrm{Frob}^T}

\newcommand{\dualLiedual}{\widetilde{\mathfrak{g}}}
\newcommand{\Winf}{\mathsf{W}}
\newcommand{\locally}{}

\newcommand{\Cores}{\mathrm{Cores}}
\newcommand{\Res}{\mathrm{Res}}

\newcommand{\Ad}{\mathrm{Ad}}
 \newcommand{\Hom}{\mathrm{Hom}}
 
 \newcommand{\Hop}{\mathrm{H}^{\mathrm{op}}}
 \newcommand{\Lie}{\mathrm{Lie}}
 \newcommand{\triv}{\mathrm{triv}}
 \newcommand{\ab}{\mathrm{ab}}

\newcommand{\ur}{\mathrm{ur}}
\newcommand{\rhobar}{\overline{\rho}}
\newcommand{\gHecke}{\tilde{\mathbb{T}}}
\newcommand{\End}{\mathrm{End}}
\newcommand{\Chains}{\mathrm{Chains}}

\newcommand{\ind}{\mathrm{ind}}
\newcommand{\Frob}{\mathrm{Frob}}
\newcommand{\Sym}{\mathrm{Sym}}
\newcommand{\SL}{\mathrm{SL}}

\newcommand{\Vinf}{\mathsf{V}}
\newcommand{\ZZ}{\mathbf{Z}}
\newcommand{\Cinf}{\mathsf{C}}
\newcommand{\Sinf}{\mathsf{S}}
\newcommand{\Iinf}{\mathsf{I}}
\newcommand{\Rinf}{\mathsf{R}}

\newcommand{\Aut}{\mathrm{Aut}}
\newcommand{\aff}{\mathrm{aff}}

\renewcommand{\dim}{\mathrm{dim}}

\newcommand{\Q}{\mathbf{Q}}

\newcommand{\C}{\mathbf{C}}
\newcommand{\adele}{\mathbf{A}}
\newcommand{\mubar}{\overline{\mu}}
\newcommand{\lambdabar}{\overline{\lambda}}
\newcommand{\OO}{\mathcal{O}}
\newcommand{\et}{\mathrm{et}}
\newcommand{\R}{\mathbf{R}}
\newcommand{\G}{\mathbf{G}}

\newcommand{\Hecke}{\mathscr{H}}
\newcommand{\dHecke}{\mathscr{H}}

\newcommand{\PGL}{\mathrm{PGL}}
\newcommand{\ramprimes}{S}

 \DeclareFontShape{OT1}{rsfs}{n}{it}{<-> rsfs10}{}
\DeclareMathAlphabet{\mathscr}{OT1}{rsfs}{n}{it}
\newcommand{\GL}{\mathrm{GL}}

\swapnumbers

\newtheorem{theorem}[subsection]{Theorem}
\newtheorem{Definition}[subsection]{Definition}
  \newtheorem{corollary}[subsection]{Corollary}
  \newtheorem*{corollary*}{Corollary}
\newtheorem*{theorem*}{Theorem}
\newtheorem*{conjecture*}{Conjecture}
\newtheorem{conjecture}[subsection]{Conjecture} 

\newtheorem{lemma}[subsection]{Lemma}

\newtheorem*{lemma*}{Lemma}
\newtheorem{prop}[subsection]{Proposition}

\theoremstyle{remark} 
\newtheorem*{remark}{Remark}

\newcommand{\Tor}{\mathrm{Tor}}

\newcommand{\BWq}{\mathsf{q}}

\newcommand{\F}{\mathbf{F}}
\newcommand{\TT}{\mathbf{A}}
\newcommand{\GG}{\mathbf{G}}

\begin{document} 

\setcounter{tocdepth}{1}
\title{Derived Hecke algebra and   cohomology of arithmetic groups}
\author{Akshay Venkatesh}

 \begin{abstract}
 
We describe a graded extension of the usual Hecke algebra: it acts
in a graded fashion on the cohomology of an arithmetic group $\Gamma$. Under favorable conditions,
the cohomology is freely generated in a single degree over this graded Hecke algebra. 

From this construction we extract an action of certain $p$-adic Galois cohomology 
groups on $H^*(\Gamma, \Q_p)$, and formulate the central conjecture: 
the motivic $\Q$-lattice inside these Galois cohomology groups preserves $H^*(\Gamma,\Q)$.
    \end{abstract}
\maketitle

    \tableofcontents

\section{Introduction}

Let  $\GG$ be a   semisimple $\Q$-group, and let $Y(K)$ 
be the associated arithmetic manifold (see   \eqref{YKdef}). 
 Particularly when $Y(K)$ is not an algebraic variety, it often happens
that the same Hecke eigensystem can occur in several different cohomological degrees (see \S \ref{sec:MV}). 
Our goal is to construct extra endomorphisms of cohomology that partly explain this,
 and give evidence that these extra endomorphisms are related to certain motivic cohomology groups. 
  
\subsection{Derived Hecke algebra}\label{introdha}
 Let $v$ be a prime, $G_v = \GG(\Q_v)$ and $K_v$ a maximal compact subgroup. The  usual Hecke algebra at $v$, with  coefficients in a (say finite) ring $S$,  can be described as    $  \Hom_{S G_v} \left(  S[G_v/K_v],  S[G_v/K_v] \right).$
 If in place of $\Hom$ we use $\Ext$ (see \S \ref{sec:dha} for more details)
 we get a graded extension, which we  may call the ``local derived Hecke algebra'':
\begin{eqnarray} \label{extdef} \mathscr{H}_{v,S} :=  \bigoplus_i \Ext^i_{S G_v} \left(  S[G_v/K_v],  S[G_v/K_v] \right) 
\end{eqnarray}
 
 Such a construction has been considered by P. Schneider \cite{Schneider} in the context of local representation theory in the case when $S$ has  characteristic $v$.  
In the present paper, however, we are solely interested in the opposite case, when $v$ is invertible on $S$.

For elementary reasons, the higher exts are ``almost'' killed by $q_v-1$, where $q_v$ is the size of the residue field; thus this algebra is of most interest when $q_v=1$ in $S$. 
In that case we have (\S \ref{Satake}, Theorem \ref{DHA:Satake}) a {\em   Satake isomorphism}: if $q_v=1$ in $S$,  then
  $ \mathscr{H}_{v,S}$ is isomorphic to the Weyl-invariants on the corresponding algebra for a torus (and is thus   graded-commutative).  
  
Now $\dHecke_{v,S}$ acts on the cohomology $H^*(Y(K), S)$  {\em in a graded fashion} -- the $\Ext^i$ component shifts degree by $+i$. (See 
 \S \ref{sec:explication}  for an explicit version, \S \ref{Arithmeticmanifolds} for the abstract version.)  In particular,
$\bigotimes_{v} \dHecke_{v,\Z/p^n}$  acts on $H^*(Y(K), \Z/p^n)$ and then
  (by passing to the limit as $n \rightarrow \infty$, \S \ref{limit}) we get a  ring of endomorphisms  
  $$ \gHecke \subset \End   \ H^*(Y(K), \Z_{p}),$$
  the ``global derived Hecke algebra.''  
  The degree zero component of $\gHecke$ is the usual Hecke algebra $\mathbb{T}$ -- i.e. the  subalgebra
  of $\End \ H^*$ generated by all Hecke operators.  Here, and elsewhere in the introduction,
  we will use in the Hecke algebra only ``good'' places $v$ relative to $K$. 
  
  {\color{\changecolor} 
   It seems to us likely that $\gHecke$ is graded commutative
   but we do not know this in general.  However, this has limited impact  for global applications because, firstly, different places commute; and secondly 
   our analysis in \S \ref{Satake} provides a large set of ``nice primes''
   at which the local derived Hecke algebra is graded commutative. In particular, 
   all of our work in this paper really analyzes the subalgebra generated by these ``nice primes''  (see    after  \eqref{VV0} for definition). 
This subalgebra is what is important  for all our applications, as will become clear from the discussion of \S \ref{sec:MV} onward. 
In many cases,  the existence of this sub-algebra is  anyway enough to force all of $\gHecke$ to be graded commutative (see Theorem \ref{maintheoremtriv} part (ii) and Proposition \ref{P853}).  }

  As we have mentioned, if we decompose $H^j(Y(K), \Z_{p})$ into eigencharacters for $\mathbb{T}$,
  one finds the same eigencharacters occuring in several different degrees $j$.  See \cite{takagi} for an elementary introduction to this phenomenon. 
    We want to see that  $\mathbb{T}$ is rich enough to   account for this. 
    
    One way of formalizing  ``rich enough'' is  to  complete the cohomology at a given character $\chi: \mathbb{T} \rightarrow \mathbf{F}_{p}$ of the usual Hecke algebra,
    and ask that $H^*(Y(K), \Z_{p})_{\chi}$  be generated over $\gHecke$ in minimal degree. In other words, we should like to check surjectivity of
    the map \begin{equation} \label{bigenough}\gHecke \otimes H^{\BWq}(Y(K), \Z_{p})_{\chi} \rightarrow H^*(Y(K), \Z_{p})_{\chi}\end{equation}      where $\BWq$ \index{$\BWq$, minimal degree of tempered cohomology} is the minimal degree where $H^*(Y(K), \Z_{p})_{\chi}$ is nonvanishing.\footnote{{\color{\changecolor} Note that, if $\gHecke$ is not known to be graded commutative,
    it is not a priori clear it preserves the $\chi$-eigenspace. However, under very mild assumptions it does, and this will be the
    situation in the cases we analyze.  Alternately one
    can switch to the graded commutative subalgebra described above. For the purpose of the introduction, then, the reader may either assume
    that $\gHecke$ preserves $H^*_{\chi}$ or assume that $\gHecke$ is replaced by the ``strict global derived Hecke algebra''   defined in \S \ref{Hecke:resplace}.
    }}
    

   In Theorem \ref{maintheoremtriv} and  Theorem  \ref{mindegreeproof} we prove this in two different cases (in both cases,
    we require the prime $p$ to be large enough): 
    
 Theorem \ref{maintheoremtriv}, proved in  \S \ref{Quillen}, studies the case
   when  
  $\mathbf{G}$ is (the $\Q$-group corresponding to) an inner form
   of $\mathrm{SL}_n$ over an imaginary quadratic field, and   $\chi$ is the character  $T \mapsto \deg(T)$ that sends any Hecke operator $T$ to
   its degree.  The main point is that, in this case, $H^*(Y(K), \Z_{p})_{\chi}$ can be described in terms of algebraic $K$-theory. 

 Theorem \ref{mindegreeproof}, proved in   \S \ref{Patching} and \S \ref{Patching2}, treats the case of $\chi$ associated to a tempered cohomological automorphic form, {\em assuming 
  the existence of Galois representations attached to cohomology classes on $Y(K)$,
   satisfying the expected properties (see \S \ref{GaloisAss}). } \footnote{Our assumptions are similar to Calegari--Geraghty \cite{CG}; however we do not need
   the assumptions on vanishing of cohomology because we allow ourselves to discard small $p$. 
     {\color{\changecolor}   Since the first version of this paper, remarkable progress has been made in analyzing
    the Galois representations attached to torsion classes. In particular, the paper
    \cite[Theorem 3.11]{ManyAuthorPaper} of Allen, Calegari, Caraiani, Gee, Helm, Le Hung, Newton, Scholze, Taylor and Thorne
    establishes, among many other results,  the key local-global compatibility needed for Taylor--Wiles patching.
    Moreover, the paper \cite{SevenAuthorPaper} of Caraiani, Gulotta, Hsu, Johansson, Mocz, Reinecke, and Shih
      eliminates the use of nilpotent ideals in the
    Galois representations originally constructed by Scholze. While these do not precisely match with the inputs
    needed for our setup of the argument, they appear to  address the key  issues,  and so to me it seems very likely that 
    one could produce an unconditional version of our analysis in the near future.}}   In this analysis we also impose some assumptions on $\chi$ for our convenience -- e.g. ``minimal level,'' and excluding congruences with other forms --  we have not attempted
   to be general. Here, the main tool of the proof is a very striking interaction between the derived Hecke algebra and the Taylor--Wiles method. We discuss this interaction   further in \S \ref{sec:explication}.

The proofs of \S \ref{Quillen} and \S \ref{Patching}--\ref{Patching2} are quite different, but they have an     
  an interesting   feature in common. In both cases, we use the derived Hecke algebra at primes $q$
such that restriction to $G_{\Q_q}$ kills certain classes in global Galois cohomology.
These classes live inside a certain dual Selmer group (specifically, the right hand side
of \eqref{Letalef} below).

That this particular dual Selmer group arises is quite striking, because
it seems to be a $p$-adic avatar of a certain motivic cohomology group;
and this same motivic cohomology group is suggested, in \cite{PV},
to act on the rational cohomology of $Y(K)$. 
 This brings us to the core motivation of this paper:
the derived Hecke algebra allows one to construct a $p$-adic realization of the 
operations on rational cohomology proposed in \cite{PV}.  
Therefore, we   digress to describe the conjectures of \cite{PV}. 
 We return to describe the remainder of  the current paper in 
\S \ref{reindexing}. 
%
%

 \subsection{Motivic cohomology}  \label{sec:MV} 
 
 This section is solely motivational, and so we will  
 freely assume various standard conjectures without giving complete references. We shall also allow ourselves to be slightly imprecise  
 in the interest of keeping the exposition brief. We refer to the paper \cite{PV} for full details. 
 
 Let $\chi: \mathbb{T} \rightarrow \Q$ be a character of the usual Hecke algebra, now with $\Q$ values. 
  We will suppose that $\chi$ is tempered and cuspidal. By this, we mean that 
  there is a collection $\pi_1, \dots, \pi_r$ of cuspidal automorphic representations, each tempered at $\infty$, 
  such that the generalized $\chi$-eigenspace $H^*(Y(K), \C)_{\chi}$  is exactly equal to the subspace of cohomology associated to the $\pi_i$s.

Consider now this generalized eigenspace with rational coefficients $$H^*(Y(K), \Q)_{\chi} \subset H^*(Y(K), \Q).$$
 One can understand its dimension data completely. To do so
we introduce some numerical invariants: let  $\delta, \BWq$ be defined such that
   \begin{equation} \delta  = \mathrm{rank} \  \mathbf{G}(\R) - \mathrm{rank} \ K_{\infty}. \end{equation} 
     \begin{equation} \label{qdef0} 2 \BWq +\delta = \dim Y(K). \end{equation} 
  Then  we have 
 (see \cite[Theorem III.5.1]{BW}; also \cite[Cor. 5.5]{Borel2})
\begin{equation} \label{dimeq} \dim H^{\BWq+i} (Y(K), \Q)_{\chi} = {\delta \choose i} \dim H^{\BWq}(Y(K), \Q)_{\chi}.\end{equation}
 
  In \cite{PV} a conjectural explanation for this numerology is proposed. Namely, 
 we construct a $\delta$-dimensional $\Q$-vector space  and 
 suggest that its exterior algebra acts on $H^*(Y(K), \Q)_{\chi}$. To define the vector space requires
a discussion first of the motive associated to $\chi$, and then of its motivic cohomology. 

\subsubsection{The Galois representation and the motive associated to $\chi$}
It is conjectured (and in some cases proven \cite{Scholze}) 
that to such $\chi$ there is, for every prime $p$,  a Galois representation 
$\rho_{\chi}: \Gal(\overline{\Q}/\Q) \rightarrow {}  \LG(\Q_{p})$, where $\LG$ is the Langlands dual group.\footnote{A priori, this takes  $\overline{\Q_{p}}$
coefficients; we will, for simplicity, assume that it can actually be defined over $\Q_{p}$. Moreover,  in general \cite[\S 3.4]{BG}  one has to
replace $\LG$ by a slightly different group to define $\rho_{\chi}$, but the foregoing discussion goes through with no change. 
 } 
 
 We shall suppose that $p$ is a good prime,  not dividing the level of the original arithmetic manifold $Y(K)$ (for the precise meaning of ``level,''
 see after \eqref{prodformula}). 
In particular, this means that $\rho_{\chi}$ should be crystalline upon restriction to $G_{\Q_p}$.

Now we shall compose $\rho$ with the co-adjoint representation $\LG \rightarrow \Aut(\widehat{\mathfrak{g}}^*)$
of $\LG$ on the dual of its own Lie algebra (here $\widehat{\mathfrak{g}}$ is the  Lie algebra of the dual group to $\G$, considered as a $\Q$-group, 
and $\widehat{\mathfrak{g}}^*$ is its $\Q$-linear dual).
The result is  
$$  \Ad^* \rho_{\chi}: \Gal (\overline{\Q}/\Q) \longrightarrow \Aut(\widehat{\mathfrak{g}}^* \otimes \Q_{p})$$

It is also conjectured that  $\Ad^* \rho_{\chi}$ should be motivic.   In other words, there should exist a weight zero  motive   $M_{\coad}$ over $\Q$, the ``coadjoint motive for $\chi$,'' whose Galois realization
is isomorphic to $\Ad^* \rho_{\chi}$:
\begin{equation} \label{Galoisrealization} H^*_{\mathrm{et}}(M_{\coad} \times_{\Q} \overline{\Q}, \Q_{p}) \simeq \Ad^* \rho_{\chi} \mbox{ (in cohomological degree $0$)}.\end{equation}
For simplicity we shall assume that $M_{\coad}$ can be taken to be a Chow motive, and    will suppose that the coefficient field of $M_{\coad}$ is equal to $\Q$.   
\footnote{Since \eqref{Galoisrealization} only determines the {\'e}tale realization, it is more natural to consider $M_{\coad}$ as a homological motive. Assuming standard conjectures,
\cite[\S 7.3 Remark 3.bis]{murre}  this can be promoted (non-canonically) to a Chow motive.  The independence of 
of the constructions that follow requires a further conjecture, e.g. the existence of the Bloch--Beilinson filtration on $K$-theory.}
 \medskip

\subsubsection{Motivic cohomology groups associated to $\chi$} \label{motivic chi} 
For such an $M_{\coad}$, and indeed for any Chow motive $M$, we can define (after Voevodsky) a bigraded family of motivic cohomology groups  $H^{a}_{\mot}(M ,\Q(q))$; the indexing is chosen so that this  admits a comparison map to the corresponding absolute {\'e}tale cohomology group $H^a_{\mathrm{et}}(M, \Q_p(q))$.

We will be solely interested in the motivic cohomology group with $a=q=1$; 
in this case, with the coadjoint motive,  the comparison with {\'e}tale cohomology gives
\begin{equation} \label{saxophone} H^1_{\mot}(M_{\coad}, \Q(1))  \otimes_{\Q} \Q_q \rightarrow H^1(G_{\Q}, \Ad^* \rho_{\chi}(1)).\end{equation}

Now Scholl \cite[Theorem 1.1.6]{scholl} has shown that one can define  (again for any Chow motive $M$ over $\Q$) 
a natural subspace  $H^a_{\mot}(M_{\Z}, \Q(q)) \subset H^a_{\mot}(M, \Q(q))$ of its motivic cohomology, informally speaking ``those classes that extend to a good model over $\mathbf{Z}$.''
Conjecturally, the analogue of the above map should now take values 
  inside the $f$-cohomology of Bloch and Kato \cite{BK}; in the case of interest  the analog of \eqref{saxophone} is now 
    $$H^1_{\mot}((M_{\coad})_{\Z}, \Q(1))  \otimes_{\Q} \Q_q \rightarrow H^1_f(G_{\Q}, \Ad^* \rho_{\chi}(1)).$$
Moreover, this map is conjecturally (\cite[5.3(ii)]{BK})  an isomorphism. 

It may be helpful to note that  Beilinson's conjecture relates this particular motivic cohomology to the value of the $L$-function for $\Ad^* \rho_{\chi}$  at the edge of the critical strip. In particular,
Beilinson's conjectures imply
that $$\dim_{\Q} H^1_{\mot}((M_{\coad})_{\Z}, \Q(1)) = \mbox{ order of vanishing  of $L(s, \Ad^* \rho_{\chi})$  at $s=0$.} $$
 A routine computation with $\Gamma$-factors shows that the right-hand side should indeed be equal to $\delta$.

To keep typography simple, we will denote the group $H^1_{\mot}((M_{\coad})_{\Z},\Q(1))$ simply by $L$: 
\begin{equation} \label{LDef} L := H^1((M_{\coad})_{\Z}, \Q(1)). \end{equation}
so that our discussion above says that, granting standard conjectures, $L$ is a $\Q$-vector space of dimension $\delta$,  and it comes with a map
\begin{equation} \label{Letalef} L \otimes \Q_p \rightarrow \underbrace{ H^1_f(G_{\Q}, \Ad^* \rho_{\chi}(1))}_{:= L_{\Q_p}}\end{equation}

\subsubsection{The complex regulator on $L$ and the conjectures of \cite{PV}}
 There is a complex analogue to \eqref{Letalef}: a complex regulator map on $L$, with target in a certain Deligne cohomology group. Since the details are not important for us,
 we just  call the target of this map $L_{\C}$
 and let  $L_{\C}^*$  be its
$\C$-linear  dual.

 In \cite{PV} we construct an action of $L_{\C}^*$ on $H^*(Y(K), \C)_{\chi}$ by degree $1$ endomorphisms, inducing
\begin{equation} \label{PVc}  H^{\BWq}(Y(K), \C)_{\chi}  \otimes \bigwedge \nolimits^i L_{\C}^* \stackrel{\sim}{\longrightarrow}
 H^{\BWq+i}(Y(K), \C)_{\chi} \end{equation}  
The main conjecture of \cite{PV} is that this action preserves rational structures, i.e.  the $\Q$-linear dual $L^*$  of $L$ carries
$H^*(Y(K), \Q)_{\chi}$  to itself.  In particular, this  means that 
\begin{equation} \label{PC333} \mbox{There is a natural graded  action of $\wedge^* L^*$ on $H^*(Y(K), \Q)_{\chi}$.}
\end{equation}

Therefore, if  one accepts the conjecture 
of \eqref{PC333},  and also believes that \eqref{Letalef} is an isomorphism, it should be possible to 
define a ``natural''
action of \begin{equation} \label{Denzel} \bigwedge \nolimits^* L^* \otimes \Q_p = \bigwedge \nolimits^* H^1_f(G_{\Q}, \Ad^* \rho_{\chi}(1))^*
\acts H^*(Y(K), \Q_p)_{\chi}.
\end{equation}
 Now there is no explicit mention of motivic cohomology, and 
this is where the current paper comes into the story:
in   \S \ref{reciprocity}, we shall explain how the derived Hecke algebra 
   can be used  to produce
such an action.   

This concludes our review of \cite{PV}; we now explain \eqref{Denzel}  a little bit more.

  \subsection{The derived Hecke algebra and Galois cohomology} \label{reindexing}
  
  The main result  of \S \ref{reciprocity} is Theorem \ref{maintheorem}, which constructs  
  an action of 
$   \bigwedge \nolimits^* H^1_f(G_{\Q}, \Ad^* \rho_{\chi}(1))^*$
   on $H^*(Y(K), \Q_p)$.   This is characterized
  in terms of the action of explicit derived Hecke operators. 
   More precisely,    we construct in \S \ref{obstructions} an isomorphism
\begin{equation} \label{dHeckeGalois} \gHecke \otimes \Q_p \simeq \mathbb{T} \otimes_{\Q_p} 
 \bigwedge \nolimits^* H^1_f(G_{\Q}, \Ad^* \rho_{\chi}(1))^*\end{equation}
 (actually, we do this in a case when $\mathbb{T}=\Z_p$, but  in general the argument  should yield the above result). 
  Informally, \eqref{dHeckeGalois}
gives an  ``indexing''  of derived Hecke operators by Galois cohomology. We will describe it concretely in a moment,  see \S \ref{Galoisindexing}. It can be viewed as a ``reciprocity law,''
because it relates the action of the (derived) Hecke algebra to the Galois representation in a direct way.

  To go further, let us assume that the map \eqref{Letalef} is indeed an isomorphism. Denote by $L^*$ the $\Q$-linear dual of $L$,
by $L_{\Q_p}^*$ the $\Q_p$-linear dual of $L_{\Q_p}$; we get also an isomorphism $L^* \otimes \Q_p \simeq L^*_{\Q_p}$.

Thus the derived Hecke algebra gives rise to an action of $\wedge^* L^*_{\Q_p}$
 on $H^*(Y(K), \Q_p)$. The fundamental conjecture, formulated precisely  as Conjecture
\ref{mainconjecture},  is then the following: 
 \begin{quote} 
{\em Let $\wedge^* L^*_{\Q_p}$ act
 on $H^*(Y(K), \Q_p)$  as described above. Then $\wedge^* L ^*\subset \wedge^* L^*_{\Q_p}$ preserves rational cohomology $H^*(Y(K), \Q) \subset H^*(Y(K), \Q_p)$.  }
\end{quote}

 The main point of this paper was  to get to the point where we can make this conjecture!  
What it says is that there is a hidden action of $L^*$ on the $\Q$-cohomology of $Y(K)$,
which can be computed, after tensoring with $\Q_p$, using the derived Hecke algebra. 

Here is the current status of evidence for this conjecture:

 \begin{itemize}
 \item[(i)] The most direct evidence (as of the time of writing) will be given in the paper \cite{HV}, which is joint work with
 Michael Harris. 
 There we develop an analog of the derived Hecke algebra in the setting of coherent cohomology, and formulate
 an analog of the conjecture in this setting. The advantage of this is we are actually able to carry out a numerical test
 (in the case of classical weight one modular forms) and it indeed works.
 
\item[(ii)]    As we have already mentioned,  the conjecture  should be seen as a a $p$-adic analog of the conjecture of \cite{PV} (which tells the archimedean story).
In the archimedean case, we  are able to give substantive evidence for the conjecture by other methods (periods of automorphic forms, and analytic torsion). 

 \item[(iii)] Suitably phrased, the computations of \S \ref{Quillen} can be seen as supporting a modified version of the conjecture. 
It is also easy to verify that the conjecture holds for tori, as we shall discuss in  \S \ref{pisspoor} of this paper.
\end{itemize}

 \begin{remark} 
 Note that, because of our fairly strong  assumptions, \eqref{dHeckeGalois} is even true {\em integrally} in the setting of \S \ref{reciprocity},
 i.e. the global derived Hecke algebra is an exterior algebra over $\Z_p$. I don't expect this to be true in general;
 however, the rational statement  \eqref{dHeckeGalois} should remain valid. One might imagine that  the derived deformation ring of 
 \cite{GV} will have better integral properties than the derived Hecke algebra. 
 \end{remark}

 \subsection{Explication, Koszul duality, Taylor-Wiles}  \label{sec:explication}
 
 We now explain the action of the derived Hecke algebra, and its relationship to Galois cohomology, 
 as explicitly as possible, in the case  when $Y(K)$ is an arithmetic hyperbolic $3$-manifold. Besides explicating
 the foregoing abstract discussion, this will also have the advantage that it allows us to explain the relationship
 between the derived Hecke algebra and the Taylor--Wiles method. 
 
   Suppose  $\GG$ arises from $\PGL_2$ over an imaginary quadratic field $F$, i.e. $\GG = \mathrm{Res}_{F/\Q} \PGL_2$. 
   Let  $\OO$ be the ring of integers of $F$. 
   Therefore the associated manifold $Y(K)$ (see \eqref{YKdef}) is a finite union of hyperbolic $3$-orbifolds. 
Let us suppose, for simplicity, that the class number of $F$ is odd; then, at full level, the associated arithmetic manifold
is simply the quotient of hyperbolic $3$-space $\mathbb{H}^3$ by $\PGL_2(\OO)$. 

In what follows, we fix a prime $p$ and will work always with cohomology with $\Z/p^n$ coefficients. 

Let $\mathfrak{q}$ be a prime ideal of $\OO$, relatively prime to $p$,  and let  $\F_{\mathfrak{q}} =  \OO/\mathfrak{q} $ the residue field. 
Let $$\alpha:   \mathbf{F}_{\mathfrak{q}} ^*  \longrightarrow \Z/p^n$$
be a homomorphism.   By means
of the natural homomorphism
$$\Gamma_0(\mathfrak{q}) \rightarrow  \F_{\mathfrak{q}}^*$$
sending $\left( \begin{array}{cc} a & b\\ c & d \end{array}  \right) \mapsto a/d$, 
we may regard $\alpha$ as a cohomology class $\langle \alpha \rangle \in H^1(\Gamma_0(\mathfrak{q}), \Z/p^n)$. 
Here, $\Gamma_0(\mathfrak{q})$ is as usual defined by the condition that $c \in \mathfrak{q}$. 
\medskip

 Then a typical ``derived Hecke operator'' 
of degree $+1$ is the following (see \S \ref{sec:concrete} for more): 
\begin{equation} \label{concrete}  T_{\mathfrak{q},\alpha}: H^1 \PGL_2(\OO) \stackrel{\pi_1^*}{\leftarrow}  H^1  \Gamma_0(\mathfrak{q}) \stackrel{\cup \langle \alpha \rangle }{\longrightarrow} H^2 \Gamma_0(\mathfrak{q}) \stackrel{\pi_{2*}}{ \rightarrow}  H^2 \PGL_2(\OO).  \end{equation}
Here  $\pi_1, \pi_2$ are the two natural maps $\Gamma_0(\mathfrak{q}) \rightarrow \PGL_2(\OO)$. 
 
In words,
we pull back to level $\Gamma_0(\mathfrak{q})$, cup with $\langle \alpha \rangle$,  and push back (the ``other way'') to level $1$. 
If we omitted the cup product, we would have the usual Hecke operator $T_{\mathfrak{q}}$.
The class $\alpha$ itself is rather uninteresting -- it is a ``congruence class'' in the terminology of \cite{CaVe}, i.e. it becomes trivial on a  congruence subgroup -- but nonetheless this operation seems to be new even in this case. 

The  role of torsion coefficients is vital:
If we took the coefficient ring above to be $\Z$,  there are no homomorphisms
$\mathbf{F}_{\mathfrak{q}}^* \rightarrow \Z$; more generally, 
  $\mathbf{F}_{\mathfrak{q}}^*$ has only torsion cohomology in  positive degree. 
  In fact, even to obtain ``interesting'' operations with $\Z/p$ coefficients,
we need at least that $p$ divide $\mathrm{N}(\mathfrak{q})-1$ (i.e. that $N(\mathfrak{q}) =1$ in the coefficient ring $\Z/p$,
as we mentioned in \S \ref{introdha}).

 What that means is that  elements of $\gHecke $ in characteristic zero
 necessarily arise in a very indirect way:
  as a limit of operations from $\mathscr{H}_{\mathfrak{q}}$ for {\em larger and larger primes $\mathfrak{q}$.}  
This situation is perhaps reminiscent of the Taylor-Wiles method,
and indeed one   miracle of the story is that, although the definition of $\gHecke$
is completely natural, it interacts in a rich way with the Taylor-Wiles method (not merely with
its output, e.g. $R=T$ theorems, but with the internal structure of the method itself).   

To see why this is so let us 
 examine  \eqref{concrete}: the Taylor-Wiles method studies the action of 
$\mathbf{F}_{\mathfrak{q}}^*$ on the cohomology of $\Gamma_1(\mathfrak{q})$ (these are the ``diamond operators'').   On the other hand
when we study $T_{\mathfrak{q},\alpha}$, 
 we are studying the action of $H^*(\mathbf{F}_p^*)$ on $H^*(\Gamma_0(\mathfrak{q}))$. In both settings
 it is vital that  $N \mathfrak{q}-1$ be divisible by high powers of $p$.

But these two  actions just mentioned are very closely related. More generally, 
if a group $G$ acts on a space $X$, the action of $G$ on  homology of $X$
and the action of its group cohomology $H^*(G)$ on the equivariant cohomology $H_G^*(X)$
are closely related: when $G$ is a compact torus, 
for example, this relationship is just Koszul duality \cite{GKM}.  This is  just the situation we are in, with $G = \mathbf{F}_{\mathfrak{q}}^*$,
and $X$ the  classifying space of $\Gamma_1(\mathfrak{q})$.

\subsection{Relationship to Galois cohomology: the ``reciprocity law''} \label{Galoisindexing} 

  Continuing our discussion from \S \ref{sec:explication},   let us describe explicitly how the operator $T_{\mathfrak{q},\alpha}$ is related to Galois cohomology.  
  Said differently, we are explicating the indexing of derived Hecke operators by Galois cohomology that
  is implicit in \eqref{dHeckeGalois}.
The result could be considered to be a reciprocity law, in the same sense as the usual relationship between Hecke operators and Frobenius eigenvalues.
  
  This discussion is (probably inevitably) a bit more technical.  We must again localize our story to a given Hecke eigenclass and  also make some further assumptions on the prime $\mathfrak{q}$.  
For a more precise discussion and proofs, see \S \ref{Galoisindexingproofs} of the main text. 

Fix now a character $\chi: \mathbb{T} \rightarrow \Z_{p}$
of the usual Hecke algebra at level $Y(K)$. 
 Let $$\rho: \Gal(\bar{F}/F) \rightarrow \GL_2(\Z_{p})$$ be the  Galois representation 
 conjecturally associated to $\chi$, 
and let $\rho_m$ be its reduction modulo $p^m$.
We shall assume that $\rho$ is crystalline at all primes above $p$, and also that $p >2$.  Let $\ramprimes$ be the set of finite primes  at which $\rho$ is ramified
(necessarily including all primes above $p$).

Let $\Ad \rho$ be the composite of $\rho$ with the adjoint
representation of $\PGL_2$; 
we will think of the underlying space
of $\Ad \rho$ as the space of $2 \times 2$ matrices with trace zero and entries in $\Z_{p}$. 
Also, let 
  $\Ad^* \rho$  be the $\Z_{p}$-linear dual to $\Ad \rho$ (this is identified with $\Ad \rho$ as a Galois module, so long as $p \neq 2$,
   but we prefer to try to keep them conceptually separate). Finally, $\Ad^* \rho(1)$
will be the Tate-twist of $\Ad^* \rho$. 
 
 Let $\mathfrak{q} \notin \ramprimes$ be a prime of $F$ and let 
 $F_{\mathfrak{q}}$ be the completion of $F$ at $\mathfrak{q}$.  Embed
\begin{equation} \label{Weird-Embed} \Z_{p} \mbox{ with trivial $\Gal(\overline{F_{\mathfrak{q}}}/F_{\mathfrak{q}})$ action}  \hookrightarrow \Ad  \ \rho|_{\Gal(\overline{F_{\mathfrak{q}}}/F_{\mathfrak{q}})}\end{equation}
$$ 1 \mapsto  2 \rho(\mathrm{Frob}_{\mathfrak{q}})-   \mathrm{trace} \rho(\mathrm{Frob}_{\mathfrak{q}}). $$
Explicitly, $\rho(\mathrm{Frob}_{\mathfrak{q}})$ is a $2 \times 2$ matrix over $\Z_{p}$, and the right-hand side above
is a $2\times 2$ matrix over $\Z_{p}$ with trace zero. This eccentric looking formula is a special case of a construction that makes sense for all groups, 
see \S \ref{Galoisindexingproofs}.

The map \eqref{Weird-Embed} gives rise to a similar embedding $\Z/p^m \rightarrow \Ad \ \rho_m|_{G_{F_{\mathfrak{q}}}}$, and thus a pairing of $G_{F_{\mathfrak{q}}}$-modules: 
$$\Z/p^m \times \Ad^* \ \rho(1) \rightarrow \mu_{p^m}. $$
 Thus  by local reciprocity we get a pairing
$$ H^1(F_{\mathfrak{q}},\Z/p^m) \times H^1(F_{\mathfrak{q}}, \Ad^* \ \rho(1)) \rightarrow \Z/p^m \Z,$$
and then (by restricting the second argument to $F_{\mathfrak{q}}$) 
$$ H^1(F_{\mathfrak{q}},\Z/p^m) \times H^1_f(\mathcal{O}[\frac{1}{S}],  \Ad^* \ \rho(1)) \rightarrow \Z/p^m \Z,$$
Here $H^1(\mathcal{O}[\frac{1}{S}], -) $ denotes the subspace of classes in $H^1(F, -)$ that
are unramified outside $S$, and the $f$ subscript means that we restrict further
to classes that are crystalline at $p$. 
 
Now take, as in
\S \ref{sec:explication}, an
 element $\alpha:   \F_{\mathfrak{q}}^* \rightarrow \Z/p^m$
indexing the derived Hecke operator $T_{\mathfrak{q},\alpha}$, 
and make an arbitrary extension  to a homomorphism $\tilde{\alpha}: F_{\mathfrak{q}}^*/(1+\mathfrak{q}) \rightarrow \Z/p^m$. 
This defines a class $\tilde{\alpha} \in H^1(F_{\mathfrak{q}},\Z/p^m)$, well defined up to unramified classes.
The pairing of $\tilde{\alpha}$ with $H^1_f(\mathcal{O}[\frac{1}{S}], \Ad^* \ \rho(1))$, as above, is easily seen to be independent of choice; thus from
a prime ideal $\mathfrak{q}$ and  a homomorphism $ \F_{\mathfrak{q}}^* \rightarrow \Z/p^m$ 
 we have  obtained a homomorphism: 
\begin{equation} \label{padef} [\mathfrak{q}, \alpha ] : H^1_f(\mathcal{O}[\frac{1}{\ramprimes}],  \Ad^* \rho \ (1)) \longrightarrow  \Z/p^m.\end{equation}
 In the main text of this paper (Lemma \ref{reciprocity law})
  we prove  a variant  of the following statement\footnote{Namely, we only work with simply connected groups  -- i.e., we prove
 an analogous result for $\SL_2$ rather than $\PGL_2$ -- and we impose various local conditions on the residual representations. In the introduction, we have stuck with $\PGL_2$ because it's more familiar. For example,
for $\SL_2$, we would need to use only the squares of the usual Hecke operators.  }
 under some further  local hypotheses on the representation $\rho$:

 \begin{quote}   {\em Claim:} There exists $N_0(m)$ such that
 for each pair of prime ideals $\mathfrak{q},\mathfrak{q}'$ satisfying 
\begin{itemize}
\item[(a)] $\mathrm{Norm}(\mathfrak{q}) \equiv \mathrm{Norm}(\mathfrak{q}' ) \equiv 1$ modulo $p^{N_0(m)}$ 
 \item[(b)] 
  the eigenvalues of $\rho(\mathrm{Frob}_{\mathfrak{q}})$, mod $p$,  are distinct elements of $\Z/p\Z$, and the same for $\mathfrak{q}'$; 
\item[(c)] 
 $ [ \mathfrak{q},\alpha ] = [\mathfrak{q}', \alpha ']$  in the notation of \eqref{padef}
 \end{itemize} 
the actions of $T_{\mathfrak{q},\alpha}$ and $T_{\mathfrak{q}', \alpha'}$
 on $H^*(Y(K), \Z/p^m)$ coincide.    
 \end{quote} 
 
 This is a ``reciprocity law,'' of the same nature as the reciprocity law relating Frobenius and Hecke eigenvalues.
 It is the basis for \eqref{dHeckeGalois}. 
 
 It is not as precise as one would like, because  of the annoying extra conditions on $\mathfrak{q},\mathfrak{q}'$ but it is good enough to get \eqref{dHeckeGalois}.
 It is certainly natural to believe that $$[\mathfrak{q},\alpha] = [\mathfrak{q}',\alpha']  \stackrel{?}{\implies} T_{\mathfrak{q},\alpha} = T_{\mathfrak{q}', \alpha'}$$
 (where the equality on the right is a an equality of endomorphisms of $\Z/p^m$-cohomology), without imposing condition (a)  or (b) above. 
  It would be good to prove not only this  but a version that gives information at bad places.  
 Such a formulation is presumably related to a derived deformation ring, as we  describe in the next \S.

\subsection{Further discussion and problems} \label{fdp}

It is not really surprising in retrospect that such cohomology operations should exist.  It took me a  long time to find   them because
of their subtle feature of being patched together from torsion levels.   There is a relatively simple   archimedean analog made via differential forms, see \cite{PV}.

 Here are some metaphors and problems:
\begin{itemize}
\item[(a)]
In the ``Shimura'' context a corresponding structure is provided by ``Lefschetz operators'' (although these act nontrivially
only for nontempered representations). 
But the derived Hecke algebra operators do not recover this structure. Indeed, for weight reasons, 
one expects that the higher degrees of the derived Hecke algebra act trivially in the Shimura case. 
The example of  $\GL_2$ over a field with both real and complex places shows a mixture of features, which would be interesting to study further.

 \item[(b)] The theory of completed cohomology of Calegari--Emerton \cite{CE} already predicts that, if we pass up a congruence tower,
cohomology becomes (under certain conditions) concentrated in a single degree. Said another way,
all the degrees of cohomology have ``the same source,'' and thus one expects to be able to pass from one to another. 

For this reason, it will be interesting to study the action of the {\em mod $p$ derived Hecke algebra of a {\em $p$-adic group}}; but
we stay away from this in the current paper.   (Our results and a global-to-local argument suggest that 
 this derived Hecke algebra might have a nice structure theory.  As mentioned this is studied in \cite{Schneider, SchneiderOllivier}; there is also recent work of Ronchetti \cite{Ronchetti}.)

 \item[(c)] There is also a   story of ``derived deformation rings,'' developed in \cite{GV};
 there is a pro-simplicial ring $\tilde{R}$ that represents deformations of Galois representations with coefficients in 
 simplicial rings.    
The precise definition of $\tilde{R}$, and -- assuming similar conjectures
to those assumed here --   a construction of its action on integral homology, are given in the paper \cite{GV}.
 
 However, the relationship between $\gHecke$ and $\tilde{R}$
is not one of equality: the former acts on cohomology, raising cohomological degree,
and the latter naturally acts on homology, raising homological degree. See the final section of \cite{GV} for  
a formulation of the relationship between the two actions. 

Our expection is that $\tilde{R}$ will have better integral properties than $\gHecke$, in general.

   \item[(d)] Numerical invariants: We can use $\gHecke$ to shift a class from degree $q$ to the complementary degree $\dim \ Y(K) - q$  and then cup
 the resulting classes. This gives an analog of the ``Petersson norm'' which makes sense for a torsion class (or a $p$-adic class).
What is the meaning of the resulting numerical invariants?

 \end{itemize}
 
 \subsection{Acknowledgements}

  Gunnar Carlsson  pointed out to me that my original definition (which was the one presented below in \S \ref{desc2}) should be equivalent
to the much more familiar definition with $\Ext$-groups given in the introduction.   The definition with 
$\mathrm{Ext}$-groups, or rather a differential graded version thereof,
was already defined  by Schneider around 2008 and published in \cite{Schneider}; I was unaware of Schneider's work at the time. In any case
there is little overlap between our paper and \cite{Schneider}. (See also \cite{SchneiderOllivier}). 

  Frank Calegari explained many ideas related to his paper \cite{CG} with Geraghty,  and, more generally, 
  taught me (over several years!) about  Galois representations and their deformations.  He also pointed out several typos and mistakes in the manuscript.    
The joint paper \cite{CaVe} influenced the ideas here,  e.g.  \S \ref{cong-class}.

    I thank David Treumann who explained Koszul duality to me many years ago, and more recently explained to me
 Smith theory and torus localization, which were helpful in the proof of the Satake isomorphism. 
  
  I had a very helpful discussion with Soren Galatius and Craig Westerland.

  I am grateful to Toby Gee for writing extremely clear lecture notes on modularity lifting,   without which I am not sure I would have achieved even my current modest understanding
  of  the Taylor--Wiles method, and to both Shekhar Khare and Michael Harris for taking an interest and for many helpful discussions.
     
     The referee read the paper very thoroughly and made many very helpful suggestions about exposition and clarity. I thank him or her for the substantial effort they made. 
      
     \subsection{Notation} \label{sec:notn}
     
We try to adhere to using $\ell$ or $p$ for the characteristic of coefficient rings,
     and using $q$ or $v$ for the  residue field size of nonarchimedean fields. Thus we may talk about the ``$\ell$-adic Hecke algebra of a $v$-adic group.''
     
     $\mathbf{G}$ will denote a   reductive algebraic group over a number field $F$.  \index{$\mathbf{G}$} \index{$F$} 
     In the local part of our paper --  \S \ref{sec:dha}, \S \ref{Satake}, \S \ref{IHA}  -- we shall work
     over   the  completion of such an $F$ at an arbitrary finite place. 
     In our global applications we will be more specific (just for ease of notation, e.g. not worrying about multiple primes above the residue characteristic):  \S \ref{Quillen} we take $F$ quadratic imaginary, and in \S \ref{Patching} onward we take $F=\Q$.

          It will be convenient at many points to assume that $\mathbf{G}$ is split, and then to fix  a maximal split torus $\mathbf{A}$ inside $\mathbf{G}$, and
     also a Borel subgroup $\mathbf{B}$ containing $\mathbf{A}$.   This endows the cocharacter lattice $X_*(\mathbf{A}) = \Hom(\mathbb{G}_m, \mathbf{A})$ \index{$X_*$}
     with a positive cone  $X_*(\mathbf{A})^+ \subset X_*(\mathbf{A})$, the dual to the cone spanned by the roots of $\mathbf{A}$ on $\mathbf{B}$.
     We will denote by $r = \dim X_*(\mathbf{A})$ the rank of $\mathbf{G}$. 
      
     For $v$ a place of $F$ we let $F_v$ be the completion of $F$ at $v$, $\mathcal{O}_v \subset F_v$ the integer ring, 
  $\mathbf{F}_v$ the residue field and write  $q_v$ for  the cardinality of $\mathbf{F}_v$.  
  We also put  $$G_v = \G(F_v). $$ \index{$G_v$} 
   Attached to $\mathbf{G}$ and a choice of open compact subgroup $K \subset \mathbf{G}(\adele_{F, f})$
     (the finite adele-points of $\mathbf{G}$) there is attached an ``arithmetic manifold'' $Y(K)$, which is a finite union of locally symmetric spaces: \index{$Y(K)$} 
\begin{equation} \label{YKdef} Y(K) =  \G(F) \backslash ( S_{\infty} \times \G(\Afinite))/K,\end{equation} 
where $S_{\infty}$ is the  ``disconnected symmetric space'' for $\G(F \otimes \R)$ -- the quotient of $\G(F \otimes \R)$ by a maximal compact connected subgroup. 
Although it is a minor point, we will take $Y(K)$ as an orbifold, not a manifold, and always compute its cohomology in this sense. 

As before, we introduce the integer invariants $\BWq, \delta$: 
\begin{equation} \label{qdef} \delta = \mathrm{rank}(\mathbf{G}(F \otimes \R)) - \mathrm{rank}(\mbox{maximal compact of $\mathbf{G}(F \otimes \R)$}),\end{equation}
and define $\BWq$ so that $2 \BWq+\delta = \dim Y(K)$. These have   the same significance as described in \eqref{dimeq},
at least if we suppose that the center of $\mathbf{G}$ is anisotropic over $F$.

We  will  will work only with open compact subgroups  with a product structure, i.e. 
\begin{equation} \label{prodformula} K = \prod K_v \end{equation} where
$K_v \subset \G(F_v)$ is an open compact subgroup, and $K_v$ is a hyperspecial maximal compact of $G_v$\index{$K_v$} 
for all but finitely many primes $v$.   A prime $v$ will be ``good'' for $K$ when $K_v$ is hyperspecial
and the ambient group $\G$ is quasi-split at $v$. \index{level of $K$}
The ``level of $K$'' will be, by definition, the (finite) set of all primes  $v$  which are not good.

     $\mathbf{G}$ has a dual group $G^{\vee}$, which we will regard as a {\em split Chevalley group over $\Z$;} in particular, its Lie algebra is defined over $\Z$,
     and its points are defined over any ring $R$.   We regard it as equipped with a maximal torus $T^{\vee}$ inside a Borel subgroup $B^{\vee}$.
     \index{$T^{\vee}$, dual torus} \index{$B^{\vee}$, dual Borel } 
     
       In the discussion of the Taylor--Wiles method, which takes place in \S \ref{Patching} and \S \ref{reciprocity},
     it is convenient to additionally assume:
     \begin{quote} $\mathbf{G}$ is simply connected and $G^{\vee}$ is adjoint.
     \end{quote}
     This is a minor issue, to avoid the usual difficulties of ``square roots.'' One could (better) replace $G^{\vee}$ by some version of the $c$-group of \cite{BG}.

      When we discuss Galois cohomology, we will follow the usual convention
that, for a module $M$ under the absolute Galois group $\mathrm{Gal}(\overline{L}/L)$ of a field $M$,
we denote by $H^*(L,M)$ the  continuous cohomology of the profinite group    $\mathrm{Gal}(\overline{L}/L)$  with coefficients in $M$. 
  For $L$ a number field, with ring of integers $\mathcal{O} \subset L$, we denote by $H^1(\mathcal{O}[\frac{1}{S}], M) \subset H^1(L, M)$ the subset of classes
  that are unramified outside $S$ and $H^1_f(\mathcal{O}[\frac{1}{S}], M) \subset H^1(\mathcal{O}[\frac{1}{S}], M)$ the  classes that are, moreover, crystalline at $p$.

{\color{\changecolor} Warning: In the literature, the subscript $f$  is often used to mean classes that are crystalline at $p$ and unramified at {\em all} other places,
 not merely at places outside $S$.   Indeed, we implicitly used this notation in \S \ref{sec:MV} and \S \ref{reindexing} when we wrote $H^1_f(G_{\Q}, -)$. 
 However, in the remainder of the text, we will not use this convention. To avoid any confusion, on the one occasion (in \S \ref{pisspoor})
 we wish to refer to classes that are crystalline at $p$ and unramified at all other places, we will use the notation $H^1_{f, \mathrm{ur}}$. This notation will be reprised when it is used so the reader need not remember it now.}
    
  \section{Derived Hecke algebra} \label{sec:dha}
     
We introduce the derived Hecke algebra (Definition \ref{Extdefn}) and then give
two equivalent descriptions in \S \ref{desc2} and  \S \ref{doublecoset}.   The model given in \S \ref{desc2} is by far the most useful. 
We shall then  describe the action of the derived Hecke algebra on the cohomology of an arithmetic group in 
       \S \ref{Arithmeticmanifolds}, and then make it a bit more concrete in \S \ref{sec:concrete}. 
       Finally, \S \ref{DHA:oddsandends1} discusses some minor points to do with change of coefficient ring,
       and \S \ref{limit} discusses some other minor points about passage between $\Z/\ell^n$ coefficients and $\Z_{\ell}$ coefficients. 
       
       Appendix \ref{remedial} expands on various points of homological algebra that are used in the current section.

     \subsection{}
       As in \S \ref{sec:notn},  we fix a prime $v$ of $F$, with residue field  $\mathbf{F}_v$ of characteristic $p_v$ and size $q_v$,  and set 
       $ G_v = \G(F_v)$. We denote by $U_v$ an open compact subgroup of $G_v$.
              Eventually, we will use only the case of $U_v$ being either a maximal compact subgroup or an Iwahori subgroup, but there is no need to impose this. When we are working strictly in a local setting, we will abbreviate these simply to $G$ and $U$:
       $$ G = G_v, \ \ U = U_v.$$

       It will also be convenient to fix 
\begin{equation} \label{Vvdef} V_v = \mbox{a pro-$p_v$, normal, finite index subgroup of $U_v$,}\end{equation}
which we  again abbreviate to $V$ when it will cause no confusion.

        Let $S$ be a finite coefficient ring in which $q_v$ is invertible. 
In what follows, by ``$G$-module'' we mean a module $M$ 
under the group algebra $SG$ with the property that every  $m \in M$
has open stabilizer in $G$.   
       The category of $G$-modules is an abelian category and it has enough projective objects (see \S \ref{enough_proj}).

        The usual Hecke algebra  for the pair $(G, U)$ can be defined as the endomorphism ring $\Hom_{SG}(S[X], S[X])$, where $X = G/U$ and 
$S[X]$ denotes the free $S$-module on a set $X$.  Motivated by this, we define:

       \begin{Definition} \label{Extdefn} \index{derived Hecke algebra}\index{$\mathscr{H}(G, U)$}
          The derived Hecke algebra for $(G, U)$  with coefficients in $S$ is the graded algebra \begin{equation} \label{firstdef} \mathscr{H}(G,U)_{S}:= \Ext^*(S [G/U], S [G/U]),\end{equation}
where the $\Ext$-group is taken inside the category of  $G$-modules.  

 \end{Definition}

 Let us record some variants on the notation:
 \begin{itemize}
 \item[-] We will write simply $\mathscr{H}(G, U)$ when the coefficients are understood to be $S$;  
\item[-] We write $\mathscr{H}^{j}(G, U)$ or $\mathscr{H}^{(j)}(G, U)$ for the component in degree $j$,  i.e.  the $\Ext^j$ summand on the right.  \index{$\mathscr{H}^{(j)}$}\index{$\mathscr{H}^j$}
\item[-]  In  global situations where we have  fixed a level structure $K_v \leqslant G_v$ for all $v$, or for almost all $v$, 
we will often write simply $\mathscr{H}_{v,S}$  for the corresponding derived Hecke algebra $\mathscr{H}(G_v, K_v)$.  Again we will write simply $\mathscr{H}_v$\index{$\mathscr{H}_v$}
if the coefficients are understood to be $S$. 
\end{itemize} 
 
If we choose a projective resolution $\mathbf{P}$ of $S[G/U]$
as $G$-module,  then $\mathscr{H}(G, U)$ is identified with
the  cohomology of the differential graded algebra
$\Hom_{SG}(\mathbf{P}, \mathbf{P})$.  
It will be convenient for later use to  make an explicit choice of $\mathbf{P}$:  
Let $\mathbf{Q}$ be a free resolution of the trivial module $S$ in 
the category of $S [U/V]$-modules. We may take $\mathbf{P}$ to be  
  the compact induction  (from $U$ to $G$) of $\mathbf{Q}$. Observe that all the groups $\mathbf{P}_i$
of the resulting resolution are free $S$-modules.

\subsection{Description in terms of invariant functions}  \label{desc2} 
We may also describe $\mathscr{H}(G, U)$ as the algebra of
``$G$-equivariant cohomology classes on $G/U \times G/U$ with finite support modulo $G$.'' 
We now spell out carefully what this means; 
an explicit isomorphism
between this description and Definition \ref{Extdefn} is constructed in Appendix \S \ref{remedial}. 

First some notation: for $x,y \in G/U$, 
we denote by $G_{xy}$  the pointwise stabilizer of $(x,y)$ inside $G$; it is a profinite group.
We denote by $H^*(G_{xy}, S)$  the continuous cohomology of $G_{xy}$ with coefficients in $S$ (discretely topologized).

In this model, an element of $\mathscr{H}(G, U)$ 
is an assignment $h$ that takes  as input $(x,y) \in G/U \times G/U$
and produces as output\  $h(x,y) \in H^*(G_{xy}, S)$, subject to the following conditions:
\begin{itemize}
\item $h$ is $G$-invariant, that is to say, $[g]^* h(gx, gy) =   h(x,y)$,
where $[g]^*: H^*(G_{gx,gy}) \rightarrow H^*(G_{xy})$ is pullback by $\Ad(g)$. 
 \item  $h$ has finite support modulo $G$, i.e. there is a finite subset $T \subset G/U \times G/U$  such that $h(x,y)  = 0$ if $(x,y) $ does not lie in the $G$-orbit of $T$. 
\end{itemize}

The addition and $S$-module structure on $\mathscr{H}(G, U)$ is defined pointwise. 
The product is given by the rule 
 \begin{equation} \label{explication} h_1 * h_2(x,z) = \sum_{y \in G/U} \underbrace{h_1(x,y) }_{H^*(G_{xy})} \cup \underbrace{h_2(y,z)}_{H^*(G_{yz})} \end{equation} 
 where we give the right-hand side the following meaning: 
  The cup product  on the right makes sense in $H^*(G_{xyz},S)$, i.e. first restrict $h_1$ and $h_2$ to $H^*(G_{xyz},S)$, 
  and take the cup product there. 
  Now split $G/U$ as a disjoint union $\coprod O_i$ of orbits under $G_{xz}$; let $O$ be one such orbit. 
 We regard 
\begin{equation} \label{glug} \sum_{y \in O}  h_1(x,y)  \cup h_2(y,z)  := \mathrm{Cores}^{G_{x y_0 z}}_{G_{xz}} \ \left( h_1(x, y_0) \cup  h_2(y_0, z) \right)\end{equation}
 where $y_0 \in O$ is any representative, and the ``trace'' or corestriction is taken
 from $G_{x y_0 z} $ to $G_{xz}$; note that the right-hand side of  \eqref{glug} is independent of choice of $y_0 \in O$. 
 Adding up over orbits $O$ gives the meaning of the right-hand side of \eqref{explication}.

  \begin{remark} 
  Suppose that $\Delta$ is a compact subgroup of $G$ that stabilizes every point of $G/U$. 
  In this case, we can restrict $h$ to get a function $h_{\Delta}: G/U \times G/U \rightarrow H^*(\Delta)$.
  We also have $(h h')_{\Delta}=  h_{\Delta} h'_{\Delta}$, where the right-hand multiplication
is the more familiar
\begin{equation} \label{easier for me} h_{\Delta} h'_{\Delta}(x,z) = \sum_{y \in G/U} h_{\Delta}(x,y) \cup h_{\Delta}'(y,z).\end{equation}
  \end{remark}

     \subsection{Double coset description} \label{doublecoset}
     Finally, we can describe $\mathscr{H}(G, U)$ in terms of double cosets 
     $U \backslash G/U$. 
 For $x \in G/U$ let   $$U_x = U \cap \Ad(g_x) U$$ where $g_x \in G$
represents $x$ (that is to say, $x=g_x U)$. Then $U_x$ is the stabilizer of $x$ in $U$. 

Fix a set of representatives $[U \backslash G/U] \subset G/U$
for the  left $U$-orbits on $G/U$. 
Then we have an isomorphism of $S$-modules
\begin{equation} \label{TTspin} \bigoplus_{x \in [U \backslash G/U]}   H^*(U_x, S) \stackrel{\sim}{\rightarrow} \mathscr{H}(G,U) \end{equation} 
 given thus:   Fix a class $z \in   [U \backslash G/U]$,   and a representative $g_z\in G$ for $z$ -- thus $z = g_z U$. Let $\alpha \in H^*(U_z, S)$. 
  Then the class of $\alpha \in H^*(U_z, S)$, considered as an element of the left-hand side of \eqref{TTspin},
  is carried to the function $h_{z,\alpha} $ on   $G/U \times G/U$ characterized by the following properties:
  \begin{itemize}
  \item[(i)]$h_{z,\alpha}(x,y) = 0$ unless    $(x, y) $ belongs to the $G$-orbit of $(z, eU)$. 
  \item[(ii)]  $h_{z,\alpha}$ sends   $ (z, eU) $ to $\alpha \in H^*(U_z,S)$ -- note that $U_z$ is exactly the common stabilizer of $z$ and $eU$.
  \end{itemize}
  
   This gives another description of $\mathscr{H}(G, U)$. It is harder to directly describe the multiplication rule in this presentation, and we use instead the isomorphism to the previous description. 
  Later on we'll describe explicitly the action of $h_{z,\alpha}$ on the cohomology of an arithmetic manifold. 
 
  Now let us examine the ``size'' of $\mathscr{H}(G, U)$; this discussion is  really only motivational, and so we will be a little informal. Suppose, for example,
 that $G$ is split and $U$ is hyperspecial. In this case, the quotient $U \backslash G/U$
 is parameterized by a dominant chamber  $X_*(\mathbf{A})^+$ inside the co-character group $X_*(\mathbf{A})$ of a maximal split torus $\mathbf{A}$. 
 Moreover, if $x \in G/U$  is a representative for a double coset parameterized by $\lambda \in X_*(\mathbf{A})$, then the group $U_x$ is, modulo a pro-$p$-subgroup, 
 the $\mathbf{F}_v$-points $M_{\lambda}(\mathbf{F}_v)$ of the   centralizer $M_{\lambda}$ of $\lambda$. 
 Thus we obtain an isomorphism of $S$-modules:
 $$\mathscr{H}(G, U) := \bigoplus_{\lambda \in X_*(\mathbf{A})^+}  H^*(M_{\lambda}(\mathbf{F}_v), S)$$
For ``generic'' $\lambda$ -- i.e., away from the walls of $X_*(\mathbf{A})^+$ --
the group
$M_{\lambda}$ is a split torus; the order of its $\mathbf{F}_v$-points is a power   of $(q_v-1)$.  Thus if $(q_v-1)$ is invertible on $S$, 
all the terms of $\mathscr{H}(G, U)$ corresponding to dominant $\lambda$ vanish.

In this paper we will be primarily concerned with the case when $q_v-1=0$ inside $S$.  Although it is certainly interesting to study $\mathscr{H}(G, U)$ in general,
the preceding discussion shows that this case (i.e. $q_v=1$ in $S$) is where $\mathscr{H}(G, U)$ is ``largest.''

\subsection{Derived invariants} \label{sec:derivedinvariants}

 If $M$ is any  complex of $G$-modules,  we may form the derived invariants
$$  \mbox{derived $U$-invariants on $M$} := \underline{\mathrm{Hom}}_{SG}(S [G/U], M) \in \mathbf{D}(\mathrm{Mod}_S)$$
where $\underline{\mathrm{Hom}}$ is now  derived $\mathrm{Hom}$ in the derived category of $G$-modules,
taking values in the derived category of $S$-modules.

Then the derived Hecke algebra automatically acts on the cohomology 
of the derived invariants:
\begin{equation} \label{HGUderived} \mathscr{H}(G, U) \acts  H^*(\mbox{derived $U$-invariants on $M$}).\end{equation}
Indeed, the derived invariants are represented by the complex
$\Hom_{SG}(\mathbf{P}, M)$, where $\mathbf{P}$ is as before any projective resolution of $S[G/U]$. 
The action of $\Hom_{SG}(\mathbf{P}, \mathbf{P})$ on this complex
furnishes the desired (right) action of $\mathscr{H}(G, U)$.

Let us describe the derived invariants in more familiar terms.  Let $V$ be as in \eqref{Vvdef}, and  consider the explicit projective resolution $\mathbf{P}$ discussed in \S \ref{Extdefn}; we see 
that the derived $U$-invariants are computed by the complex $\Hom_{SU}(\mathbf{Q}, M)$.
This coincides with $U/V$-homomorphisms from $\mathbf{Q}$ to the termwise invariants $M^{V}$; since $\mathbf{Q}$ is a projective resolution of $S$
in the category of $U/V$-modules, we see that 
 \begin{equation} \label{twosteps}  \mbox{derived $U$-invariants on $M$} \simeq    \underline{\mathrm{Hom}}_{S U/V}(S, M^{V}) \in  \mathbf{D}(\mathrm{Mod}_S).  \end{equation} 
where the right hand side is derived homomorphisms, in the derived category of $U/V$-modules.
 In other words,  there is an identification
 $$ H^*(\mbox{derived $U$-invariants on $M$}) \simeq \mathbb{H}^*(U/V, M^{V}),$$ 
the group hypercohomology  of the finite group $U/V$ acting on the  complex of termwise invariants $M^{V}$.

\subsection{Arithmetic manifolds} \label{Arithmeticmanifolds}

In the remainder of this section, we describe how the derived Hecke algebras act on the cohomology of arithmetic manifolds. 

We follow the notation of \S \ref{sec:notn}.  In particular, we fix
$K \subset \G(\Afinite)$ an open compact subgroup, 
which we are supposing to have a product structure $K = \prod_{w} K_w$;
let us split this as
$$ K =  K^{(v)} \times K_v$$ where $K^{(v)} = \prod_{w \neq v} K_w$ is the structure ``away from $v$.'' 
Associated to this is an arithmetic manifold $Y(K)$, as in \eqref{YKdef}. 

We will construct an action of the derived Hecke algebra $\mathscr{H}_v = \mathscr{H}_v(G_v, K_v)$   on the cohomology of $Y(K)$.
 To do so, we will exhibit $Y(K)$ as the derived $K_v$-invariants on a suitable $G_v$-module, and then apply \eqref{HGUderived}. 

For $U_v$ any open compact subgroup of $G_v$, let us abridge:
\begin{equation} \label{CUdef} C^*(U_v) = \mbox{cochain complex of  $Y(K^{(v)} \times U_v)$ with $S$ coefficients.}\end{equation}

Now let 
$$ M= \varinjlim_{U_v} C^* (U_v).$$ 
where the limit is taken over  open compact subgroups $U_v \leqslant G_v$. 
 Visibly, $M$ is a complex of smooth $G_v$-modules.  
 Choosing $V_v \subset K_v$ as in \eqref{Vvdef},    we have
$$M^{V_v} \simeq C^* (V_v),$$
since we may interchange  invariants and the direct limit; and then 
 for a finite cover $X \rightarrow Y$ with Galois group $D$
we have an isomorphism $C^*(Y) \stackrel{\sim}{\rightarrow} C^*(X)^{D}$. 
However, the   
 derived invariants of $K_v/V_v$ on $C^*(V_v)$
``coincide with'' the cohomology of $Y(K)$: the natural map
$$   C^*(Y(K)) =   C^*(V_v)^{K_v/V_v} \rightarrow  \underline{\mathrm{Hom}}_{S[K_v/V_v]}(S, C^*(V_v)) $$
is a quasi-isomorphism, in the derived category of $S$-modules. This follows from the fact that the terms
$C^*(V_v)$ have no higher cohomology as $K_v/V_v$-modules, because each $C^j(V_v)$  
  is the module of $S$-valued
 functions on a free $K_v/V_v$-set and is in particular induced from a representation of the trivial group.

We have exhibited
a quasi-isomorphism
$$C^*(Y(K)) \simeq \mbox{derived $K_v$-invariants on $M$}$$
between $C^*(Y(K))$ and a complex that represents the derived invariants of $K_v$ acting on $M$.  Thus \eqref{HGUderived} 
gives a natural right action of $\mathscr{H}(G_v, K_v)$ on the   cohomology  of $Y(K)$.  

\begin{remark}
Although this is strictly a right action,  we will often write it on the left,  which conforms more to the usual notation for Hecke operators;
the reader should therefore remember that the multiplication needs to be appropriately switched at times, but
this will cause almost no issue because the derived Hecke algebra will prove to be graded-commutative at all the places we use. 
\end{remark}

Of course, this description is totally incomprehensible; thus we now work on translating it to something more usable.

\subsection{Digression: pullback from a congruence quotient} \label{cong-class}
 We first need a brief digression to construct certain cohomology classes on $Y(K)$. These are called ``congruence classes''  
 in \cite{CaVe}, because they capitulate in congruence covers of $Y(K)$. 
 
 There is a natural map
\begin{equation} \label{group to manifold} H^*(K_v, S) \longrightarrow H^*(Y(K), S),\end{equation}
where, on the left, $H^*(K_v, S)$ is the continuous cohomology of the profinite group $K$ with coefficients in (discretely topologized) $S$. 
Indeed, any cohomology class for $H^*(K_v, S)$ is inflated from a quotient $K_v/K_{v,1}$. 
Let $K_1$ be the preimage of $K_{v,1}$ in $K$. 
The covering $Y(K_1) \rightarrow Y(K)$ has deck group $K_v/K_{v,1}$, and thus gives rise to a map, well-defined up to homotopy,
\begin{equation} \label{classifyingspacemap} Y(K) \longrightarrow \mbox{classifying space of $K_v/K_{v,1}$.}\end{equation}
We may then pull back cohomology classes along this map to get \eqref{group to manifold}.

 These ``congruence'' cohomology classes have a very simple behavior under Hecke operators:   \begin{lemma} \label{lem:HTriv}
 Let $h$ be in the image of the map \eqref{group to manifold}. 
For any prime $w$ of $F$ that does not divide the level of $K$ or the size of $S$, such that $\mathbf{G}(F_w)$ is split, and any usual Hecke operator 
$T$ supported at $w$, we have
 $$T h= \deg(T) h.$$
 \end{lemma}
 We will give a direct proof, but let us note that one can also deduce the result
from the commutativity of the Hecke algebra at $w$ (which is proved, under mild restrictions on $w$, in 
\S \ref{Satake}).   It is also likely that the assumption that $\mathbf{G}(F_v)$ is split is unnecessary (since $w$ doesn't divide the level of $K$, it is automatically quasi-split by our definitions, which should be enough
for the argument below to go through). 
 
 \proof
It is easy to verify this if $w \neq v$, so we examine only the case $w=v$.

By the assumptions, we may suppose that $K_v = \mathcal{G}(\mathcal{O}_v)$, for a split reductive $\mathcal{G}$ over $\mathcal{O}_v$. 
 Suppose that $T$ arises from the double coset $K_v a_v K_v$, where, without loss,
$a$ lies in a maximal split torus $\mathbf{A}(F_v)$ that is in good position relative to $K_v$ -- i.e.
it extends to a maximal split torus of    $\mathcal{G}$.   

We will show that  $h$ has the same pullback under the two natural maps
$$\pi_1, \pi_2: Y(K \cap \Ad(a_v) K) \rightarrow Y(K),$$
namely, the natural map, and the map induced by multiplication by $a_v$; this implies the Lemma.

There is an isomorphism $X_*(\mathbf{A}) \simeq A_v/ (A_v \cap K_v)$;
let $\mathcal{M}$ be the centralizer  in $\mathcal{G}$ of the co-character in $X_*(\mathbf{A})$
that corresponds to the class of $a$. 
Let $K_2$ be the preimage, under $K_v \rightarrow \mathcal{G}(\mathcal{O}_v/\varpi_v^D )$,
of $\mathcal{M}(\mathcal{O}_v/\varpi_v^D)$; here $D$ is a large enough  integer, and $\varpi_v$ a uniformizer.

Then, on the one hand, the inclusion $K_2 \hookrightarrow K \cap \Ad(a_v) K$ has  index
equal to a power of $q_v$. In particular, it induces an injection on $H^*(-, S)$, so it is enough 
to verify that $\pi_1^* h = \pi_2^* h$ after pullback under $Y(K_2) \rightarrow Y(K \cap \Ad(a_v) K)$. 

However,  the pullback $H^*(\mathbf{G}(\mathbf{F}_v), S) \rightarrow H^*(K_v, S)$ is an isomorphism. 
The class $h$ is therefore actually pulled back from $\mathbf{G}(\mathbf{F}_v)$. Our assertion then follows from the fact that 
the
 natural maps
$K_2 \longrightarrow  \mathbf{G}(\mathbf{F}_v)$ -- namely, the reduction map, 
and the conjugate of the reduction map by $a_v$ --  actually  coincide. 
 This proves that $\pi_1^* h = \pi_2^* h$, and concludes the proof of the Lemma. 
 \qed

 This motivates the following definition:  
 
 \begin{Definition}  \label{HTriv}
 We say a class 
$ h \in H^*(Y(K), S)$ is {\em Hecke-trivial} if,
for all places $v$ that do not divide the level of $K$ and  with residue characteristic invertible on $S$, 
 and all Hecke operators $T$  supported at $v$,
$$(T - \mathrm{deg}(T))^{n} h = 0.$$ 
for a sufficiently large integer $n=n(T)$. 
We denote by $H^*(Y(K), S)_{\triv}$ the submodule of Hecke--trivial classes.    \end{Definition}

\subsection{Concrete expression for the action of $\mathscr{H}_v$ on $H^*(Y(K), S)$.}\label{sec:concrete}

Let us now give a more down-to-earth description of the action of $\mathscr{H}_{v,S}$ on $H^*(Y(K), S)$, with notation as above.
In particular, we will show that the action of elements $h_{z,\alpha}$ can be described in a fashion that is very close to the usual definition of Hecke operators. 

From  $z = g_z K_v \in G_v/K_v$,
and $\alpha \in H^*(K_v \cap \Ad(g_z) K_v)$,  we obtain a class $h_{z, \alpha} \in \mathscr{H}(G_v, K_v)$, by the recipe of \S \ref{doublecoset}.  Then:

\begin{lemma}
Write  $$K_z = K \cap \Ad(g_z) K,  \ \ \  K_z' = K \cap \Ad(g_z^{-1}) K.$$
Also, let $\langle \alpha \rangle$ be the image of $\alpha$ under  $H^*(K_v \cap \Ad(g_z) K_v)  \stackrel{\eqref{group to manifold}}{\rightarrow} H^*(Y(K_z), S)$.   
Then the action of $h_{z,\alpha}$ on $H^*(Y(K),S)$ coincides with the following composite 

 \begin{equation} \label{action-explicit} H^*(Y(K))  \stackrel{[g_z]^*}{\rightarrow} H^*( Y(K_z))  \stackrel{\cup \langle \alpha \rangle}{\rightarrow} H^*(Y(K_z))  \rightarrow H^*(Y(K)),\end{equation}

where all cohomology is taken with $S$ coefficients, 
and  the arrows are (respectively) pullback by the map $Y(K_z) \rightarrow Y(K)$ induced
by the map $g \mapsto g g_z$, 
  cup with $\langle \alpha \rangle$, and push-forward for the standard map $Y(K_z) \rightarrow Y(K)$.  
 \end{lemma}
 
  Note that this is almost exactly the same as a usual Hecke operator; we have just inserted the operation of $\cup \langle \alpha \rangle$ on the way. 
  The fact that $\langle \alpha \rangle$ is Hecke-trivial, in the sense of Definition \ref{HTriv},  is the key point that makes this operation commute with usual Hecke operators.

  \proof
Routine but extremely tedious; see \S \ref{remedial}.
\qed

 \medskip

\noindent {\bf Remark. }  Note also the following trivial case: taking $g_z = 1$, we see that  the operation of ``cup with $\alpha \in H^*(K_v, S)$''
always belongs to the  derived Hecke algebra. 

\medskip
 
{\color{\changecolor}  \noindent {\bf Remark.}  As an example let us write out the argument that derived Hecke operators at different places always commute with one another.   
 Fix places $v \neq w$, elements $g_v \in G_v, g_w \in G_w$,
 and classes $\alpha_v \in H^*(K_v \cap \Ad(g_v) K_v)$ and $\alpha_w \in H^*(K_w \cap \Ad(g_w) K_w)$. 
We claim that the composite of the two associated derived Hecke operators
can be described in the following  way, which makes graded commutativity clear: Push-pull aong
 $$Y(K) \leftarrow Y(K_{g_v} \cap K_{g_w}) \rightarrow Y(K),$$ 
 but cup in the middle with the class of $\alpha_v \cup \alpha_w$.  To verify this claim, examine the following diagram: 
  {\small \begin{equation}
 \xymatrix{
 && Y(K_{g_v} \cap K_{g_w})  \ar[ld]^{\times g_w} \ar[rd] && \\
 & Y(K_{g_v})  \ar[ld]^{\times g_v}  \ar[rd] && Y(K_{g_w}) \ar[ld]^{\times g_w} \ar[rd] & \\ 
 Y(K) && Y(K) && Y(K)
 }
 \end{equation}}
When we write (e.g.) $\times g_v$ we mean that the map is induced by right multiplication by $g_v$. 
The composite of the derived Hecke operators is, by definition, obtained by going along the bottom two rows.
However, the  
 middle diagram is a fiber product square,
 and so the two ways of going from $Y(K_{g_v})$ to $Y(K_{g_w})$, via ``push-pull'' or ``pull-push,'' coincide. 
 To prove the desired claim, then, it suffices to show that the two pullbacks of the class
 $\langle \alpha_v \rangle \in H^*(Y(K_{g_v}))$  to
 $H^*(Y(K_{g_v} \cap K_{g_w}))$ -- via the natural map, and via the map $\times g_w$ -- actually coincide.  
 
 Equivalently, the classes obtained from $\alpha_v$, in the natural way
 on $Y(K_{g_v} \cap g_w K g_w^{-1})$ and on $Y(K_{g_v} \cap g_w^{-1} K  g_w)$
 are in fact compatible, under the map $\times g_w$ from one space to the other. 
 However, these cohomology classes are obtained from a certain covering of the spaces,
 obtained by adding extra level at $v$, and the compatibility follows from the fact that $\times g_w$
 lifts to these coverings. 
%
%
   }

\subsection{Change of coefficients} \label{DHA:oddsandends1}
 
Let us examine what happens under a change of rings $S \rightarrow S'$.  
The description of \S \ref{doublecoset} and the explicit action of \S \ref{sec:concrete}
means that there is a map of Hecke algebra $\mathscr{H}_{v,S} \rightarrow \mathscr{H}_{v, S'}$
compatible with the actions on $H^*(Y(K), S) \rightarrow H^*(Y(K), S')$. However, this does not make quite clear that the change
of rings map is an algebra homomorphism. For completeness let us explain this now, since
we will want to freely pass between $\Z/\ell^n$ coefficients for various $n$s.

The  tensor product $ \otimes_{S} S'$  is a right exact functor from $SG$-modules to $S'G $ modules 
and so it can be derived to a map 
of derived categories.  Note that this carries projectives to projectives since $\Hom_{S'G_v}(P \otimes_{S} S', -) = \Hom_{SG_v} (P, -)$.

This derived tensor product (let us write it as $\underline{\otimes}$) ``carries $S[G_v/K_v]$ to $S'[G_v/K_v]$:''
if we choose a projective replacement $\mathbf{P} \rightarrow S[G_v/K_v]$ the natural map
$$\mathbf{P} \otimes_{S} S' \longrightarrow S'[G_v/K_v]$$
is a quasi-isomorphism. 
Indeed it is possible to choose $\mathbf{P}$
so that each term of $\mathbf{P}$ is free as an $S$-module (see the explicit resolution after \eqref{firstdef}). Then $\mathbf{P} \otimes_{S} S'$
has no cohomology in higher degree (since this complex computes the $\mathrm{Tor}_S(S[G_v/K_v], S')$ and
the former is free) and thus it is a resolution of $S'[G_v/K_v]$.

This yields at once a map
$$ \mathscr{H}_{v,S} \rightarrow \mathscr{H}_{v, S'}.$$ 
from the Hecke algebra with $S$ coefficients, to the same with $S'$ coefficients.  
Explicitly, the left-hand side is represented by the cohomology of the differential graded algebra $\Hom_{SG}(\mathbf{P}, \mathbf{P})$, and this 
dga maps to $\Hom_{S'G_v}(\mathbf{P}\otimes_S S', \mathbf{P} \otimes_{S} S')$, whose cohomology computes $\mathscr{H}_{v,S'}$.  This is the desired algebra map ``change of coefficients.''

 Consider now the  obvious map
 $$ \iota: \varinjlim_{U_v} C^*(U_v)   \rightarrow \varinjlim_{U_v} C^*(U_v; S'),$$
 where the notation is as in  \eqref{CUdef}, and the right-hand side is defined the same way but with $S'$ coefficients. 
This induces 
 $$\iota':  \Hom_{SG}(\mathbf{P},  \varinjlim_{U_v} C^*(U_v)    ) \rightarrow \Hom_{S'G} \left( \mathbf{P} \otimes_S S',  \left(  \varinjlim_{U_v} C^*(U_v; S') \right) \right)$$
 wherein we compose with $\iota$ and extend by $S$-linearity.  There are compatible actions of $\mathscr{H}_{v,S}$
 and $\mathscr{H}_{v,S'}$ on the left and right sides.  On the other hand,  the map $\iota'$
 induces on cohomology the natural map $H^*(Y(K), S) \rightarrow H^*(Y(K), S')$. 
    
To summarize: the actions of $\mathscr{H}_{v,S}$ on $H^*(Y(K), S)$ and
  $\mathscr{H}_{v, S'}$ on $H^*(Y(K), S')$ are compatible, with respect to the natural maps $\mathscr{H}_{v,S} \rightarrow \mathscr{H}_{v,S'}$
  and $H^*(Y(K), S) \rightarrow H^*(Y(K), S')$.

\subsection{Passage from mod $\ell^n$ to $\ell$-adic; the global derived Hecke algebra} \label{limit}

We now write out in grotesque detail certain minor details of the  passage from mod $\ell^n$ to $\ell$-adic coefficients, which will be used without comment in our later proofs.  
 This section should probably be skipped by the reader and consulted only as needed.  
 
 {\color{\changecolor} In what follows, we will describe the straightforward version of the global derived Hecke algebra
using all good primes for $K$ that do not divide $\ell$;  here ``good'' is defined  after \eqref{prodformula}.
{\em When we refer to the global derived Hecke
algebra without any further remark, we are always referring to this version. }
We will remark, after the construction, how to make a definition
using a restricted set of primes.  }

 \index{global derived Hecke algebra} \index{$\gHecke$} 
The action of the derived Hecke algebra gives an algebra of endomorphisms $\widetilde{\mathbb{T}}_n \subset \mathrm{End}( H^*(Y(K), \Z/\ell^n))$, namely
the algebra of endomorphisms generated by all  the derived Hecke algebras $\mathscr{H}_{v, \Z/\ell^n}$
where $v$ varies over good primes of $K$ that are not above $\ell$.

Now we have
$$H^*(Y(K), \Z_{\ell}) = \varprojlim H^*(Y(K), \Z/\ell^n)$$ 
and we {\em define} the global derived algebra (relative to the fixed set $V$ of places)
to be 
\begin{equation} \label{gHeckedef} \gHecke \subset \mathrm{End}(H^*(Y(K), \Z_{\ell}))\end{equation}
to be those endomorphisms of the form $\varprojlim t_n$,   for some compatible system $t_n \in \gHecke_n$, i.e.
$t_n$ ``reduces to $t_{m}$'' for $n > m$ in
the sense that the following diagram should commute:

 \begin{equation}
 \xymatrix{
H^*(Y(K), \Z/\ell^n)  \ar[r]^{t_n} \ar[d] & H^*(Y(K), \Z/\ell^n)  \ar[d] \\
H^*(Y(K), \Z/\ell^{m})  \ar[r]^{t_{m}} & H^*(Y(K), \Z/\ell^{m}) . 
 }
\end{equation}

Let $\gHecke_{n}^*$ be the  systems of elements  
$(t_n , t_{n-1}, \dots, t_1)$, where
$t_r \in \widetilde{\mathbb{T}}_r$ for $ r  \leq n$ are all compatible in the sense that the above diagram should commute for each $t_{r}, t_{r'}$. 
In particular, $\gHecke_{n}^*$ acts on $\Z/\ell^r$-valued cohomology for each $  r \leq n$. 
The inverse limit
$\varprojlim \gHecke_n^*$ acts on $H^*(Y(K), \Z_{\ell})$, and its
image in $\End H^*(Y(K), \Z_{\ell})$ is precisely the global derived Hecke algebra $\gHecke$.

 Fix $m$. For $n \geq m$  consider
 the map
 $$\gHecke_n^* \rightarrow \tilde{\mathbb{T}}_m^*.$$ 
For increasing $n$ and fixed $m$, the image of this map  gives a decreasing sequence of subsets of  the finite set $\tilde{\mathbb{T}}_m^*$. 
This sequence must stabilize. Call this stabilization
$\gHecke_{\infty,m}$; it is a subring of $\widetilde{\mathbb{T}}_m^*$ and thus acts by 
  endomorphisms of $H^*(Y(K), \Z/\ell^m)$; 
  also,  there exists 
$N_m$ such that $\gHecke_{\infty,m}$
coincides with the image of $ \gHecke_{N_m}^*$
in $\tilde{\mathbb{T}}_m^*$.  
  
The natural map
\begin{equation} \label{ontoness} \varprojlim \gHecke_n^* \rightarrow \gHecke_{\infty,m}\end{equation}
is onto, since we're dealing with an inverse system of finite sets.

Let $\dHecke_{v, \Z/\ell^n}$ be the local derived Hecke algebra at $v$ with $\Z/\ell^n$-coefficients.
We show later (\S \ref{surj}) that, if $\ell^n$ divides $q_v -1$, then the natural map
$\dHecke_{v, \Z/\ell^n} \rightarrow \dHecke_{v, \Z/\ell^m}$ is surjective.  
 It follows that if $q_v-1$ is divisible by $\ell^{N_m}$, then  the image of 
  $ \gHecke_{N_m}^*$ acting on $H^*(Y(K), \Z/\ell^m)$
contains the image of $\dHecke_{v, \Z/\ell^{m}}$.  
Therefore, the image of
  $\gHecke_{\infty,m}$ acting on $H^*(Y(K), \Z/\ell^m)$
contains the image of $\dHecke_{v, \Z/\ell^{m}}$.  

In practice, we will establish ``bigness'' results of the following type: 
 \begin{quote}
(*) \label{quotestar} For all $ m \leq n $, there exists sets of primes $Q_n = \{q_1, \dots, q_r\}$
such that $\ell^n$ divides $q_i-1$ and 
the image of $\otimes_i \dHecke_{q_i, \Z/\ell^m}$ acting on $H^*(Y(K), \Z/\ell^m)$
is ``large:'' $H^*(Y(K), \Z/\ell^m)$ is generated over $\otimes_i \dHecke_{q_i, \Z/\ell^m}$    by elements  of some fixed degree $D$. \end{quote} 
  When we prove such results, it will not be for the full cohomology of $Y(K)$ but rather for its localization at some ideal of the Hecke algebra,
  but we suppress that for the current discussion.

  Let us prove that, under this assumption (*), $H^*(Y(K), \Z_{\ell})$ is generated over $\gHecke$ by elements  of degree $D$.  
    The assumption implies (by the previous discussion, with $n=N_m$) that $H^*(Y(K), \Z/\ell^m)$ is generated over
 $\gHecke_{\infty,m}$ by elements of degree $D$; by \eqref{ontoness}, it is also generated over $\varprojlim \gHecke_n^*$   by elements of degree $D$.
 That is to say,
 $$  (\varprojlim \gHecke_n^*) \otimes H^{D}(Y(K), \Z/\ell^m)$$
 surjects onto $H^*(Y(K), \Z/\ell^m)$ for every $m$.  
 By a compactness argument
  the same assertion holds with $\Z_{\ell}$-coefficients.

 More generally, the same type of argument allows us to show that various types of ``largeness''
  can be passed from $\Z/\ell^m$ coefficients to $\Z_{\ell}$.   
    {\color{\changecolor} 

  \subsection{Restricting places and the strict global derived Hecke algebra.} \label{Hecke:resplace} 
   We can restrict the primes and the powers of $\ell$ used in the above construction.
It is convenient to index these restrictions by a function
  $$V: \mbox{primes} \longrightarrow \{0, 1, 2, \dots, \} \cup \{\infty\}$$
  where primes $v$ with $V(v) = 0$ will be not used at all in the definition.
  
   In the above construction, replace $\widetilde{\mathbb{T}}_n$
   by the algebra generated by $\dHecke_{v, \Z/\ell^n}$ where $n \leq V(v)$. 
   Proceeding as above, then,  we obtain a  restricted global Hecke algebra $\gHecke^{(V)}$
   acting on cohomology with $\Z_{\ell}$ coefficients. Informally, 
  each prime $v$ can be used at most at torsion level $\Z/\ell^{V(v)}$. 
   
   For example, taking
\begin{equation} \label{VV0} V_0(v) = \begin{cases} 0, v \mbox{ not good}; \\ 
   \mbox{largest power of $\ell$ dividing $q_v-1$}, \mbox{else}     \end{cases},\end{equation}
      the resulting algebra $\gHecke^{(V_0)}$ has the advantage that
   it will be graded commutative by the results of       \S \ref{Satake},
   at least if $\ell$ doesn't divide the order of the Weyl group. \index{$\gHecke^{(V)}$} \index{$\gHecke^{(V_0)}$}

    It will be sometimes convenient to enlarge this by the usual Hecke algebra, i.e. defining the ``strict'' global derived Hecke algebra\index{strict global Hecke}
   $$ \gHecke' :=  \mbox{algebra generated by $\gHecke^{(V_0)}$ and all  underived Hecke operators at good places prime to $\ell$.}$$
  
    However, by default, when we write $\gHecke$, we mean the ``full'' version using $V=\infty$ for all good primes not above $\ell$,
   and $V=0$ at all other primes. 
   Thus we have inclusions of algebras, each inside endomorphisms of cohomology:   
   $$\mbox{usual (underived) Hecke algebra $\mathbb{T}$} \subset  \gHecke' \subset \gHecke.$$

  The advantage of $\gHecke'$ is that it is clearly graded commutative. 

  Thus, for example, if $\mathfrak{m}$ is a maximal ideal of $\mathbb{T}$, 
  the strict global derived Hecke algebra
  $\gHecke'$ induces an algebra of endomorphisms of the $\mathfrak{m}$-completion
  $H^*(Y(K),\Z_{\ell})_{\mathfrak{m}}$:
  $$     H^*(Y(K),\Z_{\ell})_{\mathfrak{m}} \otimes \gHecke' \longrightarrow H^*(Y(K),\Z_{\ell})_{\mathfrak{m}}.$$
  While {\em a priori} we do not know
  that the full $\gHecke$ preserves $H^*(Y(K), \Z_{\ell})_{\mathfrak{m}}$, this 
  is true under a mild additional assumption: For each good place $w$ not equal to $\ell$, 
  let $\mathbb{T}^{(w)}$ be the prime-to-$w$ usual Hecke algebra,
  and $\mathfrak{m}^{(w)}$ the induced maximal ideal. 
  Suppose that the natural map
\begin{equation} \label{removing w} H^*(Y(K), \Z_{\ell})_{\mathfrak{m}^{(w)}} \rightarrow H^*(Y(K), \Z_{\ell})_{\mathfrak{m}} \end{equation}
  is an isomorphism. This is true, for example, if there exists Galois representations associated to cohomology classes
  (by the argument of \cite[Lemma 6.20]{KT}).    
In this case, the local derived Hecke algebra at $w$ clearly  preserves the left-hand side, 
and so also preserves the right-hand  side.  Since this is true for all good $w$ not dividing $\ell$, the full $\gHecke$
also preserves the $\mathfrak{m}$-completion of cohomology.

   }

  \section{Torus localization and Satake isomorphism} \label{Satake} 
  
   Our main goal here is to prove a version of the Satake isomorphism that applies to the derived Hecke algebra. Namely, 
take $m =\ell^r$  a prime power.   Suppose $q \equiv 1$ modulo $\ell^r$. We show  
   (see \eqref{dSatake} for the precise statement) 
\begin{eqnarray*} \mbox{derived Hecke algebra for  split $q$-adic group with $\Z/\ell^r$-coefficients} \\  \cong  \left(\mbox{derived Hecke algebra for maximal torus with $\Z/\ell^r$ coefficients}\right)^{W} \end{eqnarray*}
where the $W$ superscript means Weyl-fixed, and we also require that $\ell$ does not divide the order of $W$. 

For example,   if $q \equiv 1$ modulo $\ell$, the derived Hecke algebra of $\PGL_2(\Q_q)$ with coefficients in $\Z/\ell$  is isomorphic to 
   $$ \Z/\ell[x_0, x_0^{-1},  y_1, z_2]^{\Z/2}$$
  where $x_0, y_1, z_2$ have (respectively) degree $0,1,2$, and the action of $\Z/2$ switches $x_0^{\pm 1}$ and negates $y_1, z_2$. 
    
 A consequence of our results is that (under our assumptions on $q, m$)  the derived Hecke algebra is graded commutative. We do not know
 if this is valid without any assumption on $q$ and  the coefficient ring $S$.  Recall, however, that $q \equiv 1$  in $S$ is precisely the case where the Hecke algebra is largest, by 
 the discussion of \S \ref{doublecoset}, and understanding this case will be enough for our global analysis.

  \subsection{}

  It is a curious fact that, in characteristic dividing $q_v-1$, 
  the Iwahori-Hecke algebra of a split $F_v$-group is isomorphic to the group algebra of  its affine  Weyl group.   
  A related interesting phenomenon is that, under the same assumptions,
  the Satake isomorphism
  $$ \mbox{ Hecke algebra} \longrightarrow \mbox{Hecke algebra of torus}$$
  is given simply by {\em restriction} (!) 
  
These points  can be explained by
  ``torus localization,'' as we now explain.  Using that method we will derive our 
 Satake isomorphism  below.  Of  course this is a little bit cheap,
 but it turns out to be exactly what we need anyway.

  I am very grateful to David Treumann for conversations about this material.  In particular,
I learned about localization in the context of local geometric Langlands   from his paper {\em Smith theory and geometric Hecke algebras} \cite{Treumann}.

  \subsection{} \label{BasicNotn}
  In this section and the next, $\mathbf{G}$ will be a {\em split group over
  the nonarchimedean local field $F_v$}.
  The coefficient ring for all our Hecke algebras will be taken to be 
  $S= \Z/\ell^r$,  for a prime $\ell$ and $r \geq 1$. 
We shall  suppose that $\ell^r$ divides $(q_v - 1)$, where $q_v$ is the cardinality
  of the residue field $\mathbf{F}_v$.      We also assume that $\ell$ is relatively prime to the order of the Weyl group of $\GG$. 

  We fix other notations as follows: Let $\mathcal{G}$ be a split group over $\mathcal{O}_v$ \index{$\mathcal{G}$} \index{$\mathcal{A}$} \index{$\mathcal{B}$}
  whose generic fiber is identified with $\mathbf{G}$. Let $K_v = \mathcal{G}(\mathcal{O}_v)$,
  a maximal compact subgroup of $G_v = \mathbf{G}(F_v)$. 
  Let $\mathbf{A}$ be a maximal torus in $\mathbf{G}$, and $\mathbf{B}$ a Borel subgroup of $\mathbf{G}$ containing $\mathbf{A}$; we suppose them to\index{$G_v$} \index{$A_v$} \index{$B_v$} 
 extend to a torus $\mathcal{A}$ and Borel $\mathcal{B}$  inside $\mathcal{G}$.  We write $A_v, B_v$ for the $F_v$-points
  of $\mathbf{A},\mathbf{B}$.    We shall use the notation $\mathbf{A}(F_v)^{\circ}$ for the   maximal compact subgroup of $\mathbf{A}(F_v)$,
  and similar notation whenever the maximal compact subgroup is unique.  \index{$\mathbf{A}(F_v)^{\circ}$}
  
  Let $W$ be the Weyl group for $\mathbf{A}$.  
We write $X_*= X_*(\mathbf{A})$ for the co-character group of $\mathbf{A}$. \index{$X_*$} 
  We identify $X_*$ with $A_v/A_v\cap K_v \subset G_v/K_v$  by means of the map
\begin{equation} \label{IWidentification} \chi \in X_* = \Hom(\mathbb{G}_m, \mathbf{A}) \longrightarrow \chi(\varpi_v),\end{equation}
  with $\varpi_v$ a uniformizer. 
 
 We write for short
 $ T =\mathcal{A}(\mathbf{F}_v).$
  The reduction map $A_v \cap K_v \rightarrow  T$ splits uniquely, and so we obtain
  a ``Teichm{\"u}ller'' lift
\begin{equation} \label{teich}  T \hookrightarrow A_v \cap K_v. \end{equation}
This induces a cohomology isomorphism, with $\Z/\ell^r$ coefficients.     \index{$T$}
   
   We have a Cartan decomposition
$$G_v = K_v \cdot A_v \cdot K_v  $$
The $A_v$ component of this decomposition is 
  unique up to the action of the Weyl group $W$.
  
  \begin{theorem} \label{DHA:Satake}  
Let notations be as above; in particular the coefficient ring is always $S=\Z/\ell^r$, where $\ell^r$ divides $q_v-1$,
and $\ell$ does not divide the order of the Weyl group.

Then restriction (in the model of \S  \ref{desc2}) defines an isomorphism
   \begin{equation} \label{dSatake} \mbox{ derived Hecke algebra for $(G_v, K_v)$} 
\stackrel{\sim}{\longrightarrow} \mbox{ derived Hecke algebra for $(A_v, A_v \cap K_v)$ }^W.
  \end{equation}

  \end{theorem}
  
  Let us explicate what we mean by ``restriction.''
As per \S \ref{desc2}, an  element $h$ of the left-hand side  is an association:
$$ (x,y) \in (G_v/K_v)^2 \rightsquigarrow h(x,y) \in H^*(G_{xy}, S),$$
and  its image  $h'$ on the right-hand side is obtained
  by restricting to $A_v/(A_v \cap K_v) \hookrightarrow G_v/K_v$
  and pulling back cohomology classes under the inclusion $A_{xy} \hookrightarrow G_{xy}$. 
  The element $h'$ is clearly $A_v$-invariant, and it is also $W$-invariant:
  $$  [w]^* h'(wx, wy) =   h(x,y)$$
because of the $G$-invariance of $h$. 
 
Because $A_{xy} = A_v \cap K_v$ for each $x,y$, and the  the (Teichm{\"u}ller) inclusion \eqref{teich}
$T \hookrightarrow A_v \cap K_v$ induces a cohomology isomorphism,
we can regard $h'$  as a function $X_* \times X_* \longrightarrow H^*(T,S)$.
We will often regard $h'$ as such without explicit comment.  
The multiplication in this model is usual convolution in the $X_*$ variable, together with multiplication in $H^*(T, S)$. 
 We may therefore identify the right-hand side of \eqref{dSatake} with   \begin{equation}  \label{alterpres} \left( S[X_*] \otimes H^*(T, S) \right)^W,\end{equation}
  just as in \eqref{TTspin}, i.e. restrict to $X_* \times \{0\}$ and identify functions on $X_*$ with the group algebra in the obvious way.
  
    \subsection{} \label{surj}
  
  A useful corollary to this result is the following (although even easier, as it does not use the algebra structure):
  The induced map
$$\mbox{derived Hecke algebra  over $\Z/\ell^n$} \rightarrow \mbox{ derived Hecke algebra  over  $\Z/\ell^m$ }$$
  is a {\em surjection} for $n > m$, under our assumption that $\ell^n$ divides $q_v-1$. 
  (We used this in the discussion  of \S \ref{limit}). 
  
  In fact, we're reduced to checking the same fact when $C$ is a cyclic group of order divisible by $\ell^n$, i.e.
  $$H^*(C, \Z/\ell^n) \rightarrow H^*(C, \Z/\ell^m)$$
  is surjective.  This follows from a straightforward computation.

\subsection{Some useful Lemmas}
 
\begin{lemma} \label{nontriv}
Any nontrivial root  $\alpha$ of $\mathbf{A}$ on $\mathbf{G}$ is nontrivial
on the $\ell$-Sylow of $A_v \cap K_v$. In particular, $\alpha$ induces a nontrivial
character  $\mathcal{A}(\mathbf{F}_v) \rightarrow \mathbf{F}_v^{\times}$. 
\end{lemma}
\proof
This is just a matter of checking the residue characteristic is forced to be big enough:
if the claim is not true,  the root $\alpha$ would be divisible 
 by $\ell^r$ in $X^*(\mathbf{A})$;  but roots are divisible at most by $2$ because $\langle \alpha, \alpha^{\vee} \rangle =2$, 
 and $\ell > 2$ because it's prime to the order of the Weyl group.
 \qed

 \begin{lemma}\label{finitefields} 
Use notation as above; in particular $\mathbf{F}_v$ is a finite field of cardinality $q_v \equiv 1$ modulo $\ell^r$,
and the order of the Weyl group is not divisible by $\ell$. 

Then the restriction map from the cohomology $H^*(\mathcal{G}(\mathbf{F}_v), \Z/\ell^r)$    to Weyl-fixed cohomology of the torus $H^*(\mathcal{A}(\mathbf{F}_v), \Z/\ell^r)^W$
is an isomorphism.
\end{lemma}

\proof   Write for short  (and just for this proof) $G, A,B$ for the $\mathbf{F}_v$-points of $\mathcal{G}, \mathcal{A}, \mathcal{B}$.

 Consider the composite of restrictions 
\begin{equation} \label{sequence} H^*(G,  \Z/\ell^r)  \rightarrow H^*(B, \Z/\ell^r) \stackrel{\sim}{ \rightarrow} H^*(A, \Z/\ell^r).\end{equation}
The second map is an isomorphism and its inverse is specified by corestriction. 
Therefore we can transport the $W$-action on $H^*(A)$  to a $W$-action on $H^*(B, \Z/\ell^r)$; explicitly the action of $w$ is
\begin{equation} \label{w} \mathrm{Cores}^B_A \circ  [w]_A \circ \mathrm{Res}^{B}_A .\end{equation} 
where $[w]_A$ is pullback of cohomology classes under $\Ad(w^{-1}): A \rightarrow A$. 

We will now show that $\mathrm{Res}^G_B \circ \mathrm{Cores}^{G}_B = \sum_{w \in W} w$,
where the $w$-action on $H^*(B)$ is that just defined. 
Since $\mathrm{Cores}^{G}_B \mathrm{Res}^G_B = |W|$, which is invertible in $\Z/\ell^r$,
we see that $\mathrm{Res}^G_B$ is injective and $\mathrm{Cores}^G_B$ is surjective; so 
  $\mathrm{Res}^G_B$ is an isomorphism onto the $W$-invariants on $H^*(B, \Z/\ell^r)$,
  which implies the Lemma.

By the  usual formula \cite[Proposition 9.5]{Brown}, using the Weyl group $W$ as a system of representatives for double cosets,
the composite equals
$$\sum_{w \in W}  \mathrm{Cores}_{w  B w^{-1} \cap B}^B  \cdot \Ad(w^{-1})^* \cdot \mathrm{Res}^{B}_{B \cap w^{-1} B w}  $$
But $w^{-1} B w \cap B$ contains  $A$, and  $[w B w^{-1} \cap B: A] =1$ modulo $\ell^r$.  So we can rewrite the $w$-term as 
$$ \mathrm{Cores}_{w  B w^{-1} \cap B}^B  \cdot \Ad(w^{-1})^* \cdot  \mathrm{Cores}_{A}^{B \cap w^{-1} B w }  \mathrm{Res}^{B}_{A} = \Cores_{A}^{B} \cdot [w]_A \cdot \Res^B_A  $$
which is exactly the $W$-action on $H^*(B)$, by  \eqref{w}. \qed 

\begin{lemma}\label{ff2}  Let $G_1, G_2$ be finite groups. Suppose that $G_1 \hookrightarrow G_1 \times G_2$ is the natural inclusion, and $M$ is a module for $G_1 \times G_2$
 killed by  the order $\# G_2$ of $G_2$. Then the corestriction map $H^*(G_1, M) \rightarrow H^*(G_1 \times G_2, M)$
 is zero. 
 \end{lemma}
 \proof Indeed, the composite $H^*(G_1 \times G_2) \stackrel{Res}{\rightarrow} H^*(G_1) \stackrel{Cores}{\rightarrow} H^*(G_1 \times G_2)$
 is multiplication by the order  of $G_2$, and is therefore zero with $M$ coefficients; but the first $\mathrm{Res}$ is surjective
 because $G_1 \rightarrow G_1 \times G_2$ is split. 
\qed

\begin{lemma}  \label{doublecentralizer} 
 Let $\Gamma \subset G_v$ be a finite $\ell$-subgroup.
 Let $\mathbf{S}$ be the double centralizer of $\Gamma$,
 considered as an algebraic $F_v$-subgroup of $\mathbf{G}$.
 Then:
   \begin{itemize}
  \item[(a)] $\mathbf{S}$ 
 has component group of prime-to-$\ell$ order, 
 \item[(b)]   the maximal compact subgroup $\mathbf{S}(F_v)^{\circ}$ of its $F_v$-points fixes
 every point of $G_v/K_v$ that is fixed by $\Gamma$. 
 \end{itemize} 
  
  \end{lemma}
\proof  Let $x \in G_v/K_v$ be fixed by $\Gamma$. Conjugating $\Gamma$ by $G_v$ we may suppose that $x=K_v$, the identity coset
 in $G_v/K_v$. 

 Now, the quotient of  the orders of $\mathcal{G}$ and $\mathcal{A}$
 over the finite field $\mathbf{F}_v$ is congruent to $|W|$ modulo $\ell$, because of our assumption $\ell$ divides $q-1$. Therefore there is an $\ell$-Sylow of $K_v$ contained in $A_v \cap K_v$.
Thus, further conjugating $\Gamma$ by $K_v$ we can further assume that 
 that $\Gamma \subset \mathbf{A}(F_v) \cap K_v$.   The centralizer $\mathbf{Z}(\Gamma)$ of $\Gamma$ is then a subgroup containing $\mathbf{A}$.
  The double centralizer $\mathbf{S}$ is  thus contained in $\mathbf{A}$ and, of course, it contains $\Gamma$.

Because $\mathbf{S} \subset \mathbf{A}$, the maximal compact subgroup 
 of $\mathbf{S}(F_v)$ is contained in the maximal compact $\mathbf{A}(F_v) \cap K_v$ of $\mathbf{A}$; the latter fixes $x$. This
 proves (b). 
  
  To verify the assertion about the component group of $\mathbf{S}$, we first verify that $\mathbf{Z}(\Gamma)$ is connected. 
  Note that $\mathbf{S}$ is contained in $\mathbf{Z}(\Gamma)$ by the analysis above, so it is in fact the center of $\mathbf{Z}(\Gamma)$. 
  Then we are reduced to the following assertion: 
 for any reductive group $\mathbf{Z}$, the component group of the center of $\mathbf{Z}$
  is only divisible by primes dividing the order of the Weyl group. Replacing $\mathbf{Z}$ by  
  its quotient by the connected center, we can check the same assertion for $\mathbf{Z}$ semisimple;
  so it is enough to check for $\mathbf{Z}$ simply connected semisimple. There it is obvious case by case. 
   
  To see that $\mathbf{Z}=\mathbf{Z}(\Gamma)  $ is connected, we can reason as follows:  $\mathbf{A}$ is a maximal split
  torus  within  $\mathbf{Z}$, so any element of $\mathbf{Z}/\mathbf{Z}^0$ has a representative in $\mathbf{Z}$  that belongs to the normalizer of $\mathbf{A}$.
  Here $\mathbf{Z}^0$ denotes the connected component.
  So it is enough to show that any  $\overline{F_v}$-point $n$ in the normalizer of $\mathbf{A}$ that belongs to $\ZZ$ actually belongs to $\ZZ^0$.
  Let $w \in \Aut(\mathbf{A})$ be the  element of the Weyl group of $\mathbf{A}$ corresponding to such an $n$. 
  Fix $\gamma \in \Gamma$. Since $n$ centralizes $\Gamma$, we see that 
  $w$ fixes $\gamma$. 
Write $N$ for the $\ell$-part of $q-1$. So $\gamma \in \mathbf{A}[N] \simeq X_* \otimes \mu_{N}$;
  fixing a primitive $N$th root, we can identify $\mathbf{A}[N]$ with $X_*/N$. 
  Since the order of $w$ is relatively prime to $\ell$, we see -- by taking invariants in $X_* \rightarrow X_* \rightarrow X_*/N$ --
  that $\gamma$ actually lies in the image of some $w$-fixed character $\mathbb{G}_m \rightarrow \TT$.   
 
  Applying this reasoning for each $\gamma \in \Gamma$,
  we see that $w$ actually  centralizes a subtorus of $\TT$ containing $\Gamma$. But the centralizer of that torus
  is a connected group, thus contained in $\mathbf{Z}^0$.  We conclude that $\mathbf{Z}$ is connected, as we claimed. 
  \qed

 \begin{lemma} \label{splitness}  Suppose $z \in G_v/K_v$ does not belong to the image of $X_*$
 (where the map  $X_* {\hookrightarrow}  G_v/K_v$ was  defined in
\eqref{IWidentification}). 
 
 Let $\Gamma$ be an $\ell$-Sylow of $A_v \cap K_v$; 
 thus $\Gamma$ is an $\ell$-Sylow subgroup of $G_{xy}$.
 Let $\Gamma_z$ be the stabilizer of $z$ in $\Gamma$.  
 
 Then the corestriction map  $H^*(\Gamma_z) \rightarrow H^*(\Gamma)$ is zero with $\Z/\ell^r$ coefficients.
\end{lemma}

 \proof

 Note that the centralizer  and so also the double centralizer of $\Gamma$ is simply $\mathbf{A}$. (Any root of $\mathbf{A}$ is nontrivial on the $\ell$-Sylow of $A_v \cap K_v$,  by Lemma \ref{nontriv}, so the connected centralizer is $\mathbf{A}$; 
the centralizer cannot be larger than $\mathbf{A}$ because any element of the Weyl group acts nontrivially on $\Gamma \simeq X_*(\mathbf{A})/\ell^r$).

Let $\mathbf{S}$ be the double centralizer of   $\Gamma_z$. 
  Since $\Gamma_z \subset \Gamma$
we also have $\mathbf{S} \subset \mathbf{A}$. 

 Let $\mathbf{S}^0$ be the identity component of $\mathbf{S}$; it is a split torus.  Because ((a) of \S \ref{doublecentralizer}) the component group of $\mathbf{S}$ is prime-to-$\ell$, we see that 
$\Gamma_z$ lies inside $\mathbf{S}^0$, and thus inside the maximal compact subgroup of $\mathbf{S}^0(F_v)$.
Let $\Gamma_z^*$ be the $\ell$-Sylow of $\mathbf{S}^0(F_v)$. Thus $\Gamma_z \subset \Gamma_z^*$. 

Choose a complement $\mathbf{S}' \subset \mathbf{A}$ to $\mathbf{S}^0$, i.e. a subtorus  
with the property that $\mathbf{S}^0 \times \mathbf{S}' \rightarrow \mathbf{A}$ is an isomorphism.  
Now $\Gamma \subset \mathbf{A}(F_v)^{\circ}$ is an $\ell$-Sylow by computation of orders, so
therefore $$\Gamma = \Gamma_z^* \times \Gamma'$$
where $\Gamma'$ is the $\ell$-Sylow of $\mathbf{S}'(F_v)^{\circ}$. 

If $\mathbf{S}'$ were trivial, then $\mathbf{S}^0 = \mathbf{A}$;   in that case,  by (b) of \S \ref{doublecentralizer}, 
$z$ lies in the fixed set of $\mathbf{A}(F_v)^{\circ}$, which is none other\footnote{ Let $\mathbf{N}$ be the unipotent radical of the Borel $\mathbf{B} \supset \mathbf{A}$. 
If $\mathbf{A}(F_v)^{\circ} x K_v = x K_v$, we have  $x^{-1} \mathbf{A}(F_v)^{\circ} x \subset K_v$.
By using the Iwasawa decomposition, it is enough to check that if this inclusion holds
 for some $x=n \in \mathbf{N}(F_v)$, then in fact $n \in K_v$.   
 In that case we have $n^{-1}  a  n \in K_v$   for all $a \in \mathbf{A}(F_v)^{\circ}$,
 and in particular $n^{-1} ( \Ad(a) n) \in K_v$ for all such $a$.
 
Choose a generic positive element $\lambda \in X_*(\mathbf{A})$, giving an enumeration of the positive roots $\alpha_1, \dots, \alpha_s$
so that $\langle \alpha_i, \lambda \rangle$ is increasing. 
For each such root
 we have a root subspace $u_{i}: \mathbf{G}_a \rightarrow
 N$, and the product map  $u_s(x_s) u_{s-1}(x_{s-1}) \dots u_1(x_1)$, 
 from $\mathbf{G}_a^s \rightarrow \mathbf{N}$, 
extends to an isomorphism
 of  schemes    over $\mathcal{O}_v$. In this ordering, the commutator
$[u_{i}, u_j]$ involves only $u_{k}$ with $k > \max(i,j)$. 
 
Let $x_1$ be the $\alpha_1$ coordinate of $n$. We have $(\alpha_1(a)-1) x_1 \in \mathcal{O}_v$
 for all $a \in \mathbf{A}(F_v)^{\circ}$, which implies $x_1 \in \mathcal{O}_v$, cf.  second paragraph of the proof of Lemma 
\ref{splitness}. Adjust $n$ on the right by $u_{1}(-x_1)$
 to arrange that $x_1$ is trivial. Now proceed the same way for the $\alpha_2, \alpha_3, \dots$ coordinate. 
 } than 
$$X_* \subset G_v/K_v,$$
which contradicts our assumption. Therefore, $\mathbf{S}'$ is nontrivial. We see at once that the order of $\Gamma'$ is divisible
by the $\ell$-part of $q-1$.

Thus, by Lemma \ref{ff2}, the corestriction from $\Gamma_z^*$ to $\Gamma$ is zero with $\Z/\ell^r$ coefficients. 
The corestriction from $\Gamma_z$ to $\Gamma$ factors through this one, so it is zero too.

  \qed

\subsection{Proof of Theorem   \ref{DHA:Satake}}\label{Satake algebra proof}
Throughout the proof, as in the statement of the theorem, we take coefficients in $S = \Z/\ell^r$. Accordingly, we drop explicit
mention of the coefficients from the notation.

 Recall the explicit description of the Satake map, using the identification
 \eqref{alterpres} of the toral derived Hecke algebra: 
 
 Given an assignment
 $(x,y) \in G_v/K_v \mapsto h(x,y) \in H^*(G_{xy})$, we associate
 to it the  function  $X_* \times X_* \longrightarrow H^*(T)$, given by 
 $$h': (x,y) \in (A_v/A_v \cap K_v)^2 \mapsto \mathrm{Res}^{G_{xy}}_{T} h(x,y) \in H^*(T)$$

 We must show that  the rule
 $h \mapsto h'$ gives  an isomorphism
   \begin{equation} \mbox{ derived Hecke algebra for $(G_v, K_v)$} 
  \simeq \mbox{ derived Hecke algebra for $(A_v, A_v \cap K_v)$ }^W.
  \end{equation}

We  first verify that $h \mapsto h'$ is bijective. Each element of  the derived Hecke algebra  for $(G_v, K_v)$ is  uniquely of the form $\sum h_{a, \alpha}$
  where $a \in X_*^+$ and $\alpha \in H^*(K_v \cap \Ad(a) K_v)$, with notation as in \S \ref{doublecoset}. 
The intersection of $K_v a K_v$ with $X_*$
  is precisely given by the $W$-orbit of $a$ by uniqueness of the Cartan decomposition. 
  So the map $h \mapsto h'$ sends $h_{a, \alpha}$ to the  function $h'_{a, \alpha}$ on $X_* \times X_*$
  characterized by $W$-invariance and:  
   \begin{itemize}
  \item[(i)]$h_{z,\alpha}(x,e) = 0$ unless    $x \in W a$; 
  \item[(ii)]  $h_{z,\alpha}$ sends   $ (a, e) $ to  the image of $\alpha \in H^*(K_v \cap \Ad(a) K_v) \rightarrow H^*(T)$. 
  \end{itemize}
    It is enough, then, to show that each element of  $\mbox{ derived Hecke algebra for $(A_v, A_v \cap K_v)$ }^W$
  is uniquely a sum of such elements $h'_{a,\alpha}$.  This comes down to the fact that the map
\begin{equation} \label{isogn} H^*(K_v \cap \Ad(a) K_v) \longrightarrow  H^*(T)^{W_a} \end{equation}
  is an isomorphism, where $W_a$ is the stabilizer of $a$ in the Weyl group. 
  But, if we write $M$ for the Levi subgroup of $\mathcal{G}$
  that centralizes $a$, then $K_v \cap \Ad(a) K_v$ is, modulo a pro-$p$-group,
  the $k$-points $M(k)$, and $W_a$ is identified with the Weyl group of $M$. 
  So \eqref{isogn} follows from Lemma \ref{finitefields}.

To show that  $h \mapsto h'$ preserves multiplication, we compute $(h_1 h_2)'(x,z)$; it equals  the restriction, from $G_{xz}$ to $T$, of
 $$ \sum_{O \subset G_v/K_v} \sum_{y \in O} h_1(x,y) \cup h_2(y,z).$$
the sum being grouped, as before, over orbits $O$ of $G_{xz}$ on such $y$. 
Recall that the inner sum  is understood by computing the cup-product $h_1(x,y) \cup h_2(y,z)$ for a single $y \in O$, 
and then corestricting from $G_{xyz}$ to $G_{xz}$.  
Therefore,

$$(h_1 h_2)'(x,z)  = \sum_{O}  \underbrace{ \mathrm{Res}^{G_{xz}}_T \mathrm{Cores}^{G_{xz}}_{G_{xyz}} h_1(x,y) \cup h_2(y,z) }_{:= H(O)}$$
and as usual we can express $H(O)$ a sum over
 $T$-orbits on $G_{xz}/G_{xyz}$, that is to say, as a sum of $T$ orbits  $O' \subset O$:
\begin{equation} \label{yuk}  H(O) =\sum_{O'} \mathrm{Cores}^{T}_{T_{y'}} \mathrm{Res}^{G_{xyz}}_{T_{y'}} (\dots )\end{equation}
where we have chosen a representative $y' \in O'$ for each $T$-orbit $O'$ upon $O$; and
  the  injection $T_{y'} \rightarrow G_{xyz}$ that defines the restriction map is induced by an element of $G_{xz}$
conjugating $y'$ to $y$. 

We saw  in Lemma \ref{splitness} that the corestriction map vanishes unless $y'$ actually belongs to $X_*$. 
(Indeed, writing $\Gamma$ for the unique $\ell$-Sylow of  the abelian group $T$, 
then $\Gamma_{y'}$ is an $\ell$-Sylow of $T_{y'}$,
and the corestriction map induced by $\Gamma_{y'} \rightarrow T_{y'}$ is surjective on cohomology.)
In the case when $y' \in X_*$, we have  $T_{y'} = T$ in which case $O' = \{y'\}$. We conclude that 
$$H(O) =\sum_{y' \in O \cap X_*}  \mathrm{Res}^{G_{xyz}}_{T}(h_1(x,y) \cup h_2(y,z))%
$$
and finally adding up all $O$ we get
\begin{eqnarray*}(h_1 h_2)'(x,z) &=& \sum_{y \in X_*}  \mathrm{Res}^{G_{xyz}}_{T} (h_1(x,y) \cup h_2(y,z) ) \\ 
 &=& \sum_{y \in X_*}  \mathrm{Res}^{G_{xy}}_{T}  h_1(x,y) \cup  \mathrm{Res}^{G_{yz}}_{T}  h_2(y,z)  =   h_1' h_2' (x,z).\end{eqnarray*}
This concludes the proof of the theorem. 
 
 \section{Iwahori-Hecke algebra}\label{IHA}
 In this section, we collect a few important facts about Iwahori--Hecke algebras. 
In particular, we discuss the structure of the Iwahori--Hecke algebra at a Taylor--Wiles prime (\S \ref{Iwahori}), the
relation between modules over the (usual, i.e. non-derived) Iwahori-Hecke algebra
and modules over the (usual) spherical Hecke algebras (\S \ref{IwahoriHecke})
and finally briefly discuss a localization result for the derived Iwahori-Hecke algebra
(\S \ref{dhLocalization}). 

 These results are presumably well-known to experts
but they help us polish our presentation of the Taylor-Wiles method -- indeed similar ideas appear in the paper of Khare and Thorne \cite{KT}.

\subsection{} 
 
 We continue with the notation of the prior section (\S \ref{BasicNotn}). 
 In particular, $G_v$ is the $F_v$-points of a reductive split group. 
 
 In this section, we will also make use of the affine Weyl group $\tilde{W}$ attached to $\mathbf{G}$; by definition
  \index{$\tilde{W}$} 
 this is the semidirect product $X_*  \rtimes W$  where $X_*$ is the cocharacter group of  the maximal torus $\mathbf{A}$,
 and $W$ is the Weyl group of $\mathbf{A}$. 
 
 Let $S$ be the ring $\Z/\ell^r$, for a prime $\ell$; this will be the coefficient ring for all our Hecke algebras and derived Hecke algebras.    We suppose that $q_v \equiv 1$ modulo $\ell^r$
 and that $\ell$ doesn't divide the order of the Weyl group.  

{\color{\changecolor} Let $I_v$ be an Iwahori subgroup of $G_v$ contained inside $K_v$ and in good position with reference to $\mathbf{A}$.  By this we mean that $I_v$ stabilizes
a chamber of the building that lies inside the apartment defined by $\mathbf{A}$. 
 An explicit choice of such an $I_v$ can be obtained from an integral model $\mathcal{B}$ of a Borel subgroup containing $\mathbf{A}$:
$$ I_v = \mbox{preimage of $\mathcal{B}(\mathbf{F}_v)$ inside $\mathcal{G}(\mathcal{O}_v)$,}$$
and the other such $I_v$s are $W$-conjugate to this one. }

It will be helpful to keep in mind that the index $[K_v:I_v] \equiv |W|$ modulo $\ell^r$, in particular, this index is invertible in $S$.
  Take the Haar measure on $G_v$ which assigns $I_v$ mass $1$.

 \subsection{The structure of the Iwahori algebra} \label{Iwahori}   
   {\color{\changecolor} Let $\mathrm{H}_I$ be the Hecke algebra for $I_v$.
   We will understand this to be defined as
   $$ \mathrm{H}_I  := \Hom_{SG_v}(S[G_v/I_v], S[G_v/I_v]).$$

This is identified with the 
set of
 $S$-valued and finitely supported functions $f$  on $I_v \backslash G_v/I_v$. Namely, 
 identifying such functions with measures (multiplying by the Haar measure on $G_v$, thought
 of as valued in $S$), each such 
function $f$ acts by right convolution on $S[G_v/I_v]$, and therefore defines an element of $\mathrm{H}_I$. 
\index{$\mathrm{H}_K$} \index{$\mathrm{H}_{IK}$} \index{$\mathrm{H}_{KI}$} \index{$\mathrm{H}_{II}$} 
Therefore, in the text, we will often produce elements of $\mathrm{H}_I$ by describing
the associated bi-invariant function. 

{\em Warning:} the resulting identification
\begin{equation} \label{left right sequel} \mathrm{H}_I \simeq \left(\mbox{functions on $I_v \backslash G_v/I_v$} \right)\end{equation}
is an anti-isomorphism of algebras if we equip the right-hand side with the convolution product.  
When we multiply elements of $\mathrm{H}_I$, we will always understand the multiplication to be that of $\mathrm{H}_I$,
not the multiplication induced from the right hand side of \eqref{left right sequel}. 

We similarly define 
\begin{eqnarray*} \label{HKIdef} \mathrm{H}_K  =  \Hom_{SG_v}(S[G_v/K_v], S[G_v/K_v]), \\
 \mathrm{H}_{IK} = \Hom_{SG_v} (S[G_v/K_v], S[G_v/I_v]), \\ 
  \mathrm{H}_{IK} = \Hom_{SG_v} (S[G_v/I_v], S[G_v/K_v]).
\end{eqnarray*}

If $V$ is any $G$-representation, the algebra  $ \Hom_{SG_v}(S[G_v/K_v], S[G_v/K_v])$ 
acts on the right on $V^{K_v} = \Hom_{G_v}(S[G_v/K_v], V)$. Similarly, elements $\mathrm{H}_{IK}$
induces endomorphisms  $V^{I_v} \rightarrow V^{K_v}$. 
Indeed a useful mnemonic for    the subscript ``$IK$''  is that,  acting as explained above,   $\mathrm{H}_{IK}$ goes from $I$-invariants to $K$-invariants, and so on.
Also   $\mathrm{H}_{IK}, \mathrm{H}_{KI}$  are bimodules for $\mathrm{H}_K$ and $\mathrm{H}_I$.

As before, each of these can be identified with a space of functions. Thus, for example, 
 \begin{eqnarray} \label{leftright2} \mathrm{H}_{IK} \simeq \mbox{functions on $K_v \backslash G_v/I_v$},    \end{eqnarray}
 and similarly for $\mathrm{H}_{IK}, \mathrm{H}_K$. 
As before, these identifications  arise  from the right convolution action of the functions
  on $S[G_v/K_v]$ or $S[G_v/I_v]$.  
  
 Note that, somewhat contrary to what the notation might suggest,  an element of $\mathrm{H}_{IK}$
  considered as a function is {\em left} $K_v$-invariant and {\em right} $I_v$-invariant.         Again, we will use these identifications such as \eqref{leftright2} without comment to produce
  elements inside the various Hecke algebras. The same warning applies here: the identifications do not preserve multiplication; the order
  must be switched, just as \eqref{left right sequel} is an anti-automorphism.  To avoid confusion here, our convention is that products are always to be understood
  in the sense of the $\mathrm{H}_I, \mathrm{H}_K$, etc., and {\em not} via convolution of functions.


It is useful later to define  
\begin{equation} \label{eK definition} e_K = \frac{1_{K_v}}{\mathrm{measure}(K_v)}.\end{equation}
Considered as an element of $\mathrm{H}_I$ it  is idempotent. 
When  considered as an element of $\mathrm{H}_{IK}$, it carries the identity coset of $S[G_v/K_v]$ to 
  $\sum_{k \in K/I} k I_v \in S[G_v/I_v]$, and when considered as an element of $\mathrm{H}_{KI}$, carries the identity coset of $S[G_v/I_v]$
  to $\frac{1}{[K_v:I_v]} eK_v \in S[G_v/K_v]$. In particular, $e_K \in \mathrm{H}_{IK}$ induces the co-restriction  map $V^{I_v} \rightarrow V^{K_v}$,
  and $e_K \in \mathrm{H}_{KI}$ induces the map $V^{K_v} \rightarrow V^{I_v}$ which is the natural inclusion divided by the index $[K_v:I_v]$. 
   }

Because $q$ is congruent to $1$ modulo $\ell^r$, the structure of $\mathrm{H}_I$ is very simple.
 It is isomorphic simply to the group algebra 
 of the affine Weyl group:
 \begin{equation} \label{affineweyl} \mathrm{H}_I \simeq S[\tilde{W}] \end{equation}
An explicit anti-isomorphism sends the characteristic function of $I_v w I_v$ to the
 element $w$, for $w \in \tilde{W}$; in particular $e_K$ is sent
 to $\frac{1}{|W|} \sum_{w \in W} w$, the sum over the usual Weyl group. 
 
This follows from the standard presentation of the Iwahori-Hecke algebra (for a reference
with complex coefficients, see \cite[Theorem 4.2]{ChrissGinzburg}); the key point
is that the relation $(T_s-q)(T_s+1)$ simplifies to $T_s^2=1$ when $q=1$ in the coefficients.  
 Actually it is also possible to verify $\mathrm{H}_I$ is isomorphic to $S[\tilde{W}]$ by using torus localization,
although we omit the details.

\subsection{Central element and discriminant} \label{Zdis} 
Every element of $S[X_*]^W$ is central  in $S[\tilde{W}]$. Therefore, \eqref{affineweyl} yields 
a  natural map (indeed an isomorphism) from $Z := S[X_*]^W$ to the center of $\mathrm{H}_I$.  
 
Then $\mathrm{H}_I, \mathrm{H}_K$ have structures of $Z$-algebra and $\mathrm{H}_{IK}, \mathrm{H}_{KI}$
have structure of $Z$-module, all of which are compatible in the obvious way. 

For example, the  ring homomorphism $Z \rightarrow \mathrm{H}_K$
is given by $z \mapsto e_K z e_K = e_K z$. In fact this is a ring {\em isomorphism}, as follows easily from the explicit presentation.
Then e.g. $\mathrm{H}_{IK}$ has two structures of $Z$-module, one via $Z \rightarrow \mathrm{H}_I$ and one via
$Z \rightarrow \mathrm{H}_K$, and the ``compatibility'' is that these two structures coincide.

\subsection{Iwahori Hecke algebra} \label{IwahoriHecke}
  
 {\color{\changecolor}  The next statement 
 asserts that the bimodules $\mathrm{H}_{KI}, \mathrm{H}_{IK}$ give
 equivalences of categories, at least over the open subset of $\Spec(Z)$
 where the ``$W$-covering'' $\Spec \ S[X_*] \rightarrow \Spec  \ S[X_*]^W$ is {\'e}tale. }
  \footnote{{\color{\changecolor} 
I am grateful to Peter Schneider for correcting an error. In an earlier version of this paper, this Lemma was formulated
over a larger open subset of $\Spec \ Z$, but one step in the proof was not correct in this greater generality.}}
     
\begin{lemma} \label{IHequivalence}
Let $\mathfrak{m}$ be a maximal ideal of $Z$ 
over which the map $Z \rightarrow S[X_*]$ is {\'e}tale,
and write $Z'$ for the localization of $Z$ at $\mathfrak{m}$.

   Write $\mathrm{H}_K' = \mathrm{H}_K \otimes_{Z} Z'$ and define similarly
$\mathrm{H}_I', \mathrm{H}_{IK}', \mathrm{H}_{KI}'$. 
Then the bimodules $\mathrm{H}_{IK}'$ and $\mathrm{H}_{KI}'$
induce inverse equivalences of categories between $\mathrm{H}_K'$ modules and $\mathrm{H}_I'$ modules. \end{lemma}  
This is probably well-known in characteristic zero at least.


\proof   (Sketch).

Let's show, for example, that the natural map induced by multiplication  
\begin{equation} \label{targetmap}
  \mathrm{H}_{IK} \otimes_{\mathrm{H}_K}  \mathrm{H}_{KI}
 \rightarrow \mathrm{H}_I \end{equation}
yields an isomorphism after localization at $\mathfrak{m}$. 
In what follows, we denote such localization with a prime: $S[X_*]' = S[X_*] \otimes_Z Z'$. 

Note that  $\mathrm{H}_I$  is  free of rank $w$ as a right module over $S[X_*]$  (clear from \eqref{affineweyl}) 
and, being finite {\'e}tale, $S[X_*]'$ is
\locally
free as a right module over $Z'$
of rank $w$.

 So $\mathrm{H}_I'$ is \locally free of rank $w^2$ over $Z'$.
Similarly, $\mathrm{H}_{IK}'$ and $\mathrm{H}_{KI}'$ are \locally free of rank $w$
and  $\mathrm{H}_K'$ is  locally free of rank $1$ as a $Z'$-module.  The
obstruction to \eqref{targetmap} being an isomorphism then given   by the vanishing of a suitable determinant;
it is enough, therefore, to  show that
\eqref{targetmap} is onto after reducing modulo the  maximal ideal $\mathfrak{m}'$ of $Z'$.

We can extend
the natural homomorphism $Z' \rightarrow Z'/\mathfrak{m}'$
to a homomorphism $\chi: S[X_*]' \rightarrow k$, 
  with $k$ an algebraically closed field containing the finite field $Z'/\mathfrak{m}'$. 
   Note that $k$ has characteristic $\ell$, and that $\chi$ is not fixed by any element of $W$.
      

Now  a homomorphism from $X_* \rightarrow k^{\times}$
is the same as
an unramified $k$-valued character $\chi$ of the maximal torus $A_v$;
so we may identify $\chi$ to a character of $A_v$. 
We form the corresponding induced representation $V= V_{\chi}$. 
Its elements consist of locally constant $k$-valued functions on $G_{v}$
that transform according to $\chi$ on a Borel subgroup containing $A_v$. 
Now
$V^{I_v}$ is a $k$-vector space of rank $w$, and $V^{K_v}$ is a $k$-vector space of rank $1$,
and $Z$ acts on these spaces via the character $\chi$ (as follows, e.g.
from \eqref{explicit00} below).   

We now show that the natural maps
$$\mathrm{H}_{ab} \otimes_Z  k \rightarrow \Hom(V^a, V^b)$$
are isomorphisms for $a$ and $b$ belonging to $\{I, K\}$; that implies the 
desired claim, that is to say, that \eqref{targetmap} is onto after reducing modulo $\mathfrak{m}$.

Because the two sides have the same rank it is enough to check surjectivity. 
In fact, it's enough to show surjectivity in the case of 
$\mathrm{H}_{II}$ and to show that all the other maps are nonzero (because then, for example, the image for $\mathrm{H}_{IK}$ would be a nonzero subspace
of $\Hom(V^{I_v}, V^{K_v})$ which is stable under $\Hom(V^{I_v}, V^{I_v})$.) 
 The other maps are clearly nonzero: the element $e_K$ induces a nonzero map
in each of the cases $IK, KI, KK$. So we are reduced to seeing that 
\begin{equation} \label{moo} \mathrm{H}_{I} \twoheadrightarrow \Hom(V^{I_v}, V^{I_v}).\end{equation}

But there's a standard basis for $V^{I_v}$ indexed by the Weyl group:  $v_w \ (w \in W)$, whose restriction to $K$ 
is the characteristic function of the Bruhat cell  indexed by $w$. 
 The group algebra of $W$, inside $\mathrm{H}_I$, acts by permuting the elements $v_w$. 
Also the element $\lambda \in X_* $, considered again inside $\mathrm{H}_I$, acts by  \begin{equation} \label{explicit00}  \lambda \cdot v_w = \langle w \chi, \lambda \rangle v_w.\end{equation}
In other words, as a representation of $\tilde{W}$, this is the representation
induced from the generic character $\chi$, and so clearly irreducible.  The surjectivity of \eqref{moo} follows. 
 \qed

 \subsection{Localization for the derived Iwahori--Hecke algebra} \label{dhLocalization}
It will later on be helpful to make use of localization for the {\em derived} Iwahori--Hecke algebra.  

We define the derived Iwahori Hecke algebra as per the recipe of \S \ref{sec:dha}, i.e.
$$ \mathscr{H}_I := \Ext_{SG_v}^*(S[G_v/I_v], S[G_v/I_v]).$$
As before this  is isomorphic to the algebra of functions  $h$ that associate to $(x,y) \in G_v/I_v \times G_v/I_v$
a class $h(x,y) \in H^*(G_{xy}, S)$, with the product as described in \S \ref{desc2}.
  In a similar way, we get derived versions
$\mathscr{H}_{IK},\mathscr{H}_{KI}$ of the bimodules $\mathrm{H}_{KI}, \mathrm{H}_{IK}$ defined earlier.

 Now, we can consider  ``restriction to $\tilde{W}$, '' i.e.
\begin{equation} \label{Tau1}  h_1  \in \mathscr{H}_{II} \longrightarrow  h_1' \in  \mbox{ functions $\tilde{W} \times \tilde{W} \rightarrow H^*(T,S)$} \end{equation}
where $T$ is as in 
\S \ref{BasicNotn}
 and we identify $w \in \tilde{W}$ with $wI \in G_v/I_v$; and 
\begin{equation} \label{Tau2}  h_2 \in \mathscr{H}_{IK} \longrightarrow  h_2' \in \mbox{ functions $\tilde{W} \times X_* \rightarrow H^*(T,S)$} \end{equation} where here we identify $x \in X_*$ with the associated coset $xK_v$; and we used the fact that $T$ stabilizes pointwise both  $\tilde{W} \cdot I_v$ and $X_* \cdot K_v$ to restrict cohomology classes to $T$. 
 Finally we have a similar map for $\mathscr{H}_{KI}$.

 Note that the right-hand side of \eqref{Tau1} has an algebra structure,   at least  restricting to functions supported on finitely many $\tilde{W}$-orbits,
 by means of the formula $H_1 H_2 (x,y)=
 \sum_{z \in \tilde{W}} H_1(x,z) \cup H_2(z,y)$ (where all of $x,y,z$ all belong to $\tilde{W}$). 
  This algebra acts on the right-hand side of \eqref{Tau2}, by means 
  of the same formula, but $H_2$ now belongs to the right-hand side of \eqref{Tau2}, 
  and therefore $y$ is now taken to belong to $X_*$.

 \begin{lemma}
Under our current notation and assumptions,
(see \S \ref{BasicNotn}), 
the map \eqref{Tau1} is an algebra morphism.
 Similarly, the map \eqref{Tau2}  
 is compatible with the  map \eqref{Tau1} and the product $\mathscr{H}_{II} \times \mathscr{H}_{IK} \rightarrow \mathscr{H}_{IK}$;
similarly for $\mathscr{H}_{KI}$.
 \end{lemma}
 
 \proof 
 
 We want to show that (where $h_1,h_2 \in \mathscr{H}_{II}$)
  $$ \mathrm{Res}(h_1 h_2) = h_1'  h_2',$$
 where $\mathrm{Res}$ means to restrict all $G_v/I_v$ arguments to $\tilde{W}$ and restrict cohomology classes to $T$; 
we also want similar statements for the $\mathscr{H}_{II}$-action on $\mathscr{H}_{KI}$ and $\mathscr{H}_{IK}$. 
 
 By precisely the same argument as in \S \ref{Satake algebra proof}, we are reduced to the following claim:

 \begin{quote} 
{\em Claim:}   
Let $y'$ belong to either $G_v/K_v$ or $G_v/I_v$. Let $\Gamma$ be an $\ell$-Sylow of $A_v \cap K_v$. 
Let $\Gamma_{y'}$ be the stabilizer of $y'$ in $\Gamma$.
Then
 the corestriction $H^*(\Gamma_{y'}, S) \rightarrow H^*(\Gamma, S)$ vanishes, 
 unless $y' \in X_* \subset G_v/K_v$ or $y' \in \tilde{W} \subset G_v/I_v$.  \end{quote}

We repeat
the reasoning of Lemma \ref{splitness}:  let $\mathbf{S}$ be the algebraic double centralizer 
 of  $\Gamma_{y'}$.  As before, $\Gamma_{y'} \subset \Gamma$ gives $\mathbf{S} \subset \mathbf{A}$. 
Let $\mathbf{S}^0$ be the identity component of $\mathbf{S}$. By Lemma  \ref{doublecentralizer},  the component group of $\mathbf{S}$ is prime-to-$\ell$. Therefore,  
$\Gamma_{y'}$ lies inside $\mathbf{S}^0$. Let $\Gamma_{y'}^* \subset \mathbf{S}^0(F_v)^{\circ}$ be an $\ell$-Sylow of the maximal compact subgroup.  Thus $\Gamma_{y'} \subset \Gamma_{y'}^*$.

Choose a complement $\mathbf{S}' \subset \mathbf{A}$ to $\mathbf{S}^0$, i.e. a subtorus
with the property that $\mathbf{S}^0 \times \mathbf{S}' \rightarrow \mathbf{A}$ is an isomorphism.  
Then $\Gamma_{y'}^* \subset \mathbf{S}^0(F_v)^{\circ}$ is an $\ell$-Sylow,
and $\Gamma \subset \mathbf{A}(F_v)^{\circ}$ is an $\ell$-Sylow. 
Therefore $$\Gamma = \Gamma_{y'}^* \times \Gamma'$$
where $\Gamma'$ is the $\ell$-Sylow of $\mathbf{S}'(F_v)^{\circ}$. 

If $\mathbf{S}'$ were trivial, then $\mathbf{S}^0 = \mathbf{A}$.
In that case  $y'$ lies in the fixed set of $\mathbf{A}(F_v)^{\circ}$. 
In the case where $y' \in G_v/K_v$ this was proved in
  \S \ref{doublecentralizer}. However, the proof of this assertion also applies word for word
to establish the same assertion in the case $y' \in G_v/I_v$. 
 
The fixed set of  $\mathbf{A}(F_v)^{\circ}$
 on $G_v/K_v$ is $X_*$, as before, and the fixed set  of  $\mathbf{A}(F_v)^{\circ}$
on $G_v/I_v$ is precisely $\tilde{W} I_v \subset G_v/I_v$.\footnote{Here is a proof of the latter claim: 
if $gI_v$ is fixed, then   $g K_v \in X_* K_v$, and modifying $g$ by an element
of $\mathbf{A}(F_v)$, we can suppose $g \in K_v$. We are reduced to computing
the $\mathbf{A}(F_v)^{\circ}$-fixed points on $K_v/I_v$, which amount to the torus fixed points on a flag variety over $\mathbf{F}_v$ --
using Lemma \ref{nontriv} to avoid problems with small residue field, these fixed points are precisely the $wI_v$ with $w \in W$, as desired.}

Otherwise, $\mathbf{S}'$ is not trivial, the corestriction $\Gamma_{y'}^* \rightarrow \Gamma$
vanishes as before, and so the corestriction $\Gamma_{y'} \rightarrow \Gamma$ vanishes too. 
 \qed

\section{The trivial representation} \label{Quillen}

 In this section we give our first piece of global evidence that the derived Hecke algebra
 can account for the ``degree spreading'' of Hecke eigenclasses.

 \subsection{} 
Let $D$ be a division algebra of dimension $d^2$ over an imaginary quadratic field $F$.
Let $\GG$ be the algebraic group of elements of norm $1$ inside $D$. 
Let $Y(K)$ be the arithmetic manifold \eqref{YKdef}  associated to $\GG$ and a level structure $K$.
We shall suppose  $K$ to be contained in the stabilizer of some maximal order $\mathcal{O}_D$. 
Observe that $\dim \  Y(K) = d^2-1$. 
 
 In this section we study the
 derived Hecke action on  ``Hecke-trivial part''  (see Definition \ref{HTriv}) of the cohomology of $Y(K)$. 
  
 Recall that  the global derived Hecke algebra $\gHecke$ consists of all endomorphisms of $H^*(\Z_{\ell})$
that are limits (\S \ref{limit}) , under $H^*(\Z_{\ell}) \simeq \varprojlim H^*(\Z/\ell^n)$, of
endomorphisms that lie in the algebra generated by all $  \mathscr{H}_{v,\Z/\ell^n}$.

\begin{theorem}\label{maintheoremtriv} 
For all but finitely many primes $\ell$, 
{\color{\changecolor} the action of $\gHecke$ on $H^*(Y(K), \Z_{\ell})$ preserves the trivial summand $H^*(Y(K), \Z_{\ell})_{\triv}$, and:}
\begin{itemize} 
\item[(i)] The trivial part of the cohomology
$H^*(Y(K), \Z_{\ell})_{\mathrm{triv}}$ is cyclic over $\gHecke$, generated by the trivial class;
\item[(ii)]  The image $\gHecke_{\triv}$ of $\gHecke$ in $\End \ H^*(Y(K), \Z_{\ell})_{\triv}$
is graded commutative. Also $\gHecke_{\triv} \otimes \Q_{\ell}$ 
coincides with $\Q_{\ell}$-algebra generated by $H^*(Y(K), \Q)_{\triv}$ acting
 on itself by means of the cup product. 
 \end{itemize}
\end{theorem}

 Note the significance of the second part of the statement: inside the $\Q_{\ell}$-derived Hecke algebra there is a natural ``preferred''
 rational structure. Our general conjecture (Conjecture \ref{mainconjecture}) says that this should be true in great generality and the preferred rational structure is related to motivic cohomology.   Certainly the situation that we discuss here is quite easy
 compared to the general case, but nonetheless it has several points of interest. 
 
We also note that the theorem  is almost certainly {\em false} (in the form stated above) if $F$ is not totally imaginary, for reasons related to  (a) of \S \ref{fdp}.

We will deduce the Theorem from the following: 
 
\begin{lemma}  \label{pairs}
Notation as above, so that $\mathbf{G}$ is the algebraic group arising from a division algebra over  the imaginary quadratic field $F$.     
For all sufficiently large $\ell$, the following statement holds: 
\begin{quote} 
 For each integer $n$ there are infinitely many places $v$ of the field $F$, with $q_v \equiv 1 $ modulo $\ell^n$ and where the division algebra is locally split, 
such that the pullback map   of \eqref{group to manifold}
$$H^*(\mathbf{G}(\mathbf{F}_v), \Z/\ell) \rightarrow H^*(Y(K), \Z/\ell)_{\triv} $$
 is surjective.  
\end{quote}
 \end{lemma}
 Note that the map above really does take values in the Hecke-trivial cohomology, by Lemma \ref{lem:HTriv}.

\proof (Summary)
The proof of Lemma \eqref{pairs} occupies  \S \ref{sec:Recollections} 
-- 
\S \ref{sec:Chebo}.  After some initial setup, 
we give in \S \ref{abcd} certain conditions (a), (b), (c), (d) which imply the Lemma;
and then after \S \ref{abcd} we check these conditions can actually be satisfied. 
\qed

First of all, let us explain   why the Lemma implies the theorem:

 \subsection{Lemma \eqref{pairs} implies the theorem}  \label{lemmaimpliestheorem}

 \proof  
 {\color{\changecolor} First of all, $\gHecke$ preserves $H^*(Y(K), \Z_{\ell})_{\triv}$. 
 This follows by the argument described around \eqref{removing w}. Explicitly,
 fixing any prime $w$, we may find  a prime-to-$w$ usual Hecke operator $T$
for which the trivial part of $\Z_{\ell}$-cohomology coincides with the generalized zero eigenspace for $T-\deg(T)$. 
(This can be checked over $\C$, since we are supposing the prime $\ell$ to be large enough. See e.g. discussion after \eqref{raru}).}

By avoiding a finite set of $\ell$, we may clearly suppose    that $H^*(Y(K), \Z)$ has no $\ell$-torsion,
and that $\ell > d$. 
Similarly, we suppose that  
 $$H^*(Y(K), \Z)_{\triv} \otimes_{\Z} \Z/\ell \rightarrow H^*(Y(K), \Z/\ell)_{\triv}$$
 is an isomorphism: see \S \ref{lconstrain} for an explanation.

By  Lemma \ref{finitefields}, the map $H^*(\mathbf{G}(\mathbf{F}_v), \Z/\ell^n) \rightarrow H^*(\mathbf{G}(\mathbf{F}_v), \Z/\ell)$ is surjective 
if $q_v \equiv 1$ modulo $\ell^n$.
It follows from this that the surjectivity assertion of   Lemma \ref{pairs}  continues to hold with coefficients modulo $\ell^n$. 

 Now we can consider an element of $H^*(\mathbf{G}(\mathbf{F}_v), \Z/\ell^n)$
 as an element of the derived Hecke algebra for $\mathbf{G}$ at $v$ (see the Remark of \S \ref{sec:concrete}). 
So the assertion implies that the cup product action of each $h \in H^*(Y(K), \Z/\ell^n)_{\triv}$ is contained
in the action of the derived Hecke algebra; by passage to the limit,
the cup product action of $H^*(Y(K), \Z_{\ell})_{\triv}$ on itself is contained in the action of the derived Hecke algebra.

Let $B_{\ell}$ be the  image of cup product $H^*(Y(K), \Z_{\ell})_{\triv} \rightarrow \End(H^*(Y(K), \Z_{\ell})_{\triv})$. Let
 $\gHecke_{\triv}$ be the image of $\gHecke$ inside $\End(H^*(Y(K), \Z_{\ell})_{\triv})$.
 It remains to show that these are equal. This comes down to the fact that $\gHecke_{\triv}$ is contained in the commutant
 of $B_{\ell}$ and so cannot be larger than it.  In more detail:
 
 {\color{\changecolor} Let $\gHecke_{\triv}^{(V_0)} $ be the subring of $\gHecke_{\triv}$ defined by only using
local derived Hecke algebras $\dHecke_{q, \Z/\ell^n}$ with $\ell^n$ dividing $q-1$, i.e.
the restricted variant of the global derived Hecke algebra defined after \eqref{VV0}. 

We have inclusions
\begin{equation} \label{inclusions BHH} B_{\ell} \subset \gHecke_{\triv}^{(V_0)} \subset \gHecke_{\triv},\end{equation}
where the first inclusion follows from the argument just given.
To conclude the proof, we will show that these are both equalities. 
Note, first of all, that each element of $\gHecke_{\triv}^{(V_0)}$ commutes, in the graded sense, 
with each element of $\gHecke_{\triv}$. Indeed,  if $q \equiv 1$ modulo $\ell^n$,
the global action of $\dHecke_{q, \Z/\ell^n}$ is readily seen to commute with  $\dHecke_{v, \Z/\ell^n}$
for $v \neq q$, and it commutes in the case $v=q$ because $\dHecke_{q, \Z/\ell^n}$ is known to be graded commutative.

Choose now $h \in \gHecke_{\triv}$.  There is $b \in B_{\ell}$ such that $h.1 = b.1$ (here $1$ is the trivial class in $H^0$). 
Then $(h-b). 1 = 0$.  
The same is true then for both the even and odd components of $(h-b)$. 
But, as we just saw, both components commute in the graded sense, with $\gHecke^{(V_0)}_{\triv}$,
as we just saw, so in fact both of these components kill  all of $H^*_{\triv}$.   Thus $h-b= 0$, so
$h \in B_{\ell}$, as required. }
\qed

\subsection{Recollections} \label{sec:Recollections} 
 
Let $N$ be an integer. (We will shortly fix it to be ``large enough.'')

Let $\mathrm{U}_N$ be the standard unitary group, the stabilizer of $\sum_{i=1}^N |z_i|^2$. 
 There are  natural maps  
 \begin{equation} \label{mapN} \mbox{bi-invariant differential forms on $\mathrm{U}_N$}  \stackrel{\sim}{\rightarrow}  H^*(\mathrm{U}_N, \C)   \longrightarrow H^*(\GL_N(\OO_F), \C) \end{equation}
obtained by  the natural identification of  $\GL_N(\C)$-invariant differential forms on $\GL_N(\C)/\mathrm{U}_N$ 
with bi-invariant differential forms on $\mathrm{U}_N$ (and then Hodge theory).     The notation is a little confusing:  $H^*(\mathrm{U}_N)$ above refers to the singular cohomology of $\mathrm{U}_N$ as a topological space,
whereas $H^*(\GL_N(\OO_F))$ refers to the group cohomology of $\GL_N(\OO_F)$. 

Moreover, 
the algebra of invariant differential forms on $\mathrm{U}_N$ is a free exterior algebra with  primitive  generators
$\Omega_1, \Omega_3, \dots, \Omega_{2N-1}$ in degree $1, 3, 5, \dots, 2N-1$;
``primitive'' is taken with respect to the 
coproduct on cohomology, induced by $\mathrm{U}_N\times \mathrm{U}_N \rightarrow \mathrm{U}_N$. 
An  explicit representative for $\Omega_j$ can be taken as
\begin{equation} \label{oinkoink} X_1, \dots, X_j \in \mathrm{Lie}(\mathrm{U}_N) \mapsto \mbox{anti-symmetrization of $\mathrm{trace}(X_1 \cdots X_j) $} \end{equation}
The same symbols $\Omega_i$ will also be used to denote the corresponding invariant differential forms on $\GL_N(\C)/\mathrm{U}_N$. 
For later use, note
that these can be restricted to cohomology classes for $\mathrm{SU}_N$ and also  to invariant differential forms on $\SL_N(\C)/\mathrm{SU}_N$;
these restrictions kill $\Omega_1$.

There are natural inclusions $\mathrm{U}_N \hookrightarrow \mathrm{U}_{N+1}$
and $\GL_N(\OO_F) \hookrightarrow \GL_{N+1}(\OO_F)$. 
For fixed $j$ and large enough $N$, these induce isomorphisms in $H^j(-, \C)$. Moreover,
these isomorphisms are compatible with increasing $N$.   By passage to the inverse limit we get 
 $$H^*(\mathrm{U}_{\infty}, \C)  {\longrightarrow} H^*(\GL_{\infty}(\OO_F), \C)$$
 Here (e.g) $\GL_{\infty}$  means in fact $\varinjlim \GL_N$.

 Both sides here  carry compatible coproducts;  for the right hand-side we can take the coproduct induced by  ``intertwining'' map  (see e.g. \cite[Chapter 2]{Srinivas})
 $\GL_{\infty} \times \GL_{\infty} \rightarrow \GL_{\infty}.$   \footnote{To see that the coproducts are compatible amounts to say that
 the group multiplication and the intertwining map both induce the same homomorphism on the stable cohomology of $U_N$. 
 This follows from the Eckmann-Hilton argument, or more directly as follows: Embed two commuting copies $U_N^{(a)} \times U_N^{(b)} \rightarrow U_{2N}$
 by the intertwining map. The group multiplication $U_{2N} \times U_{2N} \rightarrow U_{2N}$  when
 restricted to $U_N^{(a)} \times U_N^{(b)}$ on the source gives  the intertwining map   $U_N^{(a)} \times U_N^{(b)} \rightarrow U_{2N}$. This fact, together with stability of homology, shows that the two coproducts are compatible, as claimed.}
  
The corresponding Pontryagin product on homology will be denoted by $*$.

In what follows,  we fix $N$ to be divisible by $d^2$ and chosen so large that
\begin{itemize}
\item[-]  the inclusion $\GL_N \hookrightarrow \GL_{\infty}$
induces an isomorphism of integral group homology in degrees $\leq d^2$, both
with entries in $\mathcal{O}_F$ and with entries in any residue field. (This is possible because,
by a theorem of van der Kallen \cite[p. 289]{vdK}, the range of homological stability
can be taken uniformly in these cases; indeed, van der Kallen's bounds for stability involve
only the Krull dimension in the case of commutative rings.)

We will denote the stabilization map  
$$H_i(\GL_N(\OO_F)) \rightarrow H_i(\GL_{\infty}(\OO_F))$$
by $a \mapsto a^{(\infty)}$ and its inverse by $b \mapsto b^{(N)}$. 
We will use this notation for any choice of coefficients, not merely $\Z$.

\item[-]
The map \eqref{mapN} induces an surjection in degrees   $< d^2$. (That this is possible is a consequence of Borel's result \cite[(7.5)]{Borel} asserting that cohomology classes 
in sufficiently low degree are representable by invariant differential forms on the symmetric space; Borel's result is for a semisimple group, but we readily deduce the claimed result by applying it to $\SL_N$.)  We have written ``surjection'' instead of isomorphism
just because of the issue of working with $\GL$ rather than $\SL$: the differential form corresponding to $\Omega_1$
dies under \eqref{mapN}). 
\end{itemize}

In a similar way to \eqref{mapN}, we have  an  isomorphisms
\begin{equation} \label{raru} H^*(\mathrm{SU}_d, \C) \rightarrow  \mbox{$\SL_d(\C)$ invariant diff. forms on $\mathrm{SL}_d(\C)/\mathrm{SU}_d$}  \stackrel{\sim}{\rightarrow} H^*(Y(K), \C)_{\triv}. \end{equation}
For the surjectivity of the final map:  if a differential form
$\omega$  on $Y(K)$ satisfies $T \omega = \mathrm{deg}(T) \omega$ for even one Hecke operator $T$,
then by an easy ``maximum modulus'' argument it must be invariant,\footnote{This uses compactness of $Y(K)$; in the general case the answer is substantially more complicated} i.e. represented by a $\GG(F \otimes\R)$ invariant form on the corresponding symmetric space. In particular, the cohomology $H^*(Y(K), \C)_{\triv}$
is a free exterior algebra, generated in degrees $3,5, \dots, 2d-1$.

 Finally let us recall  (Borel) that the $K$-theory of $\OO_F$ is, modulo torsion, one-dimensional in each odd degree,
 and that (Quillen) for  any finite field $\mathbf{F}$ of size $q$, 
 the even $K$-groups vanish and the odd $K_{2s+1}(\mathbf{F}) \simeq \Z/(q^s-1)$.

\subsection{The constraints on $\ell$} \label{lconstrain}

We impose the following constraints on $\ell$: 
\begin{itemize}
\item[(i)] The cohomology of $Y(K)$ with coefficients in $\ell$ is torsion-free.
\item[(ii)]  $\ell$ doesn't divide $  \gcd_{v} ( (q_v-1) \cdots (q_v^{2d-1}-1))$,
where the $\gcd$ is taken over all  $q_v \geq q_0$ for large enough $q_0$. (This $\gcd$ stabilizes for $q_0$ large enough.)
\item[(iii)]  $\ell$ is relatively prime to the numerator and denominator of the rational number  $M \in \mathbf{Q}^*$ 
 defined in \eqref{Mdef}.
 \item[(iv)]  the cohomology $H^*(Y(K), \Z/\ell)_{\triv}$ is a free-exterior algebra
on generators in degree $3,\dots, 2d-1$.
\item[(v)] $\ell > d^2$ 
\item[(vi)] The cohomology of $\GL_{\infty}(\mathcal{O}_F)$ is  free of $\ell$-torsion in degrees less than $d^2$. 
\end{itemize}

All these assertions are automatically true for $\ell$ big enough. This is obvious for (i), (ii), (iii), (v) and follows from the standard stability results for (vi). We examine (iv): 
We saw  after \eqref{raru} that 
$H^*(Y(K), \Q)_{\triv} $
is a free exterior algebra; fix generators $e_3, e_5, \dots$ that belong to $H^*(Y(K), \Z)$. 
The products of the $e_i$ are linearly independent over $\Q$, so their reductions are also  linearly independent over $\Z/\ell$ for large enough $\ell$.
It remains to show that they span $H^*(Y(K), \Z/\ell)_{\triv}$. 
But that is obvious by counting dimensions: if we fix a Hecke operator $T$, then 
for sufficiently large $\ell$ the generalized zero eigenspace of  $T-\deg(T)$  on $H^*(Y(K), \Z/\ell)$ 
has the same dimension as the generalized zero eigenspace of $T-\deg(T)$ on $H^*(Y(K), \C)$.

 This concludes the proof that all of (i) -- (vi) above are automatically valid for large enough $\ell$. 
It would be interesting to see what happens for ``bad'' $\ell$. 
\subsection{} 
Let $N$ be a large integer, as chosen  in \S \ref{sec:Recollections}.  
Fix an embedding 
\begin{equation} \label{iotadef}  \iota: \GG \hookrightarrow \SL_N\end{equation}
for some large $N$, by taking a sum of many copies of the representation
that arises from   the division algebra  acting on itself. 
 Then (in suitable coordinates) we may suppose that the open compact subgroup $K$ is carried into the standard maximal compact $\prod_{v} \SL_N(\mathcal{O}_{v})$ of $\SL_N(\Afinite)$,  and the arithmetic group $\GG(F) \cap K$ is consequently carried into $\SL_N(\OO_F)$. 

The map $\iota$ gives rise to a map of symmetric spaces, i.e. a map $$\GG(\C)/\mathrm{SU}_d \rightarrow
\SL_N(\C)/\mathrm{SU}_N,$$
where we have chosen a maximal compact for $\GG(\C)$, which is isomorphic to $\mathrm{SU}_d$,
and then we have chosen a hermitian form on $\C^N$ whose stabilizer $\mathrm{SU}_N$ contains $\iota(\mathrm{SU}_d)$. 
Thus we get an  embedding of locally symmetric spaces,   also denoted by $\iota$: 
\begin{equation} \label{def of iota} \iota: Y(K) \rightarrow   \underbrace{ \SL_N(\OO_F) \backslash \SL_N(\C)/\mathrm{SU}_N}_{\simeq B(\SL_N(\OO_F))}.\end{equation}
 We can further compose $\iota$ with the inclusion of $\SL_N(\OO_F)$ to $\GL_N(\OO_F)$ to give
 a map
 $$Y(K) \rightarrow B(\GL_N(\OO_F))$$

\subsection{} \label{abcd} 
In this section, we will formulate four claims (a)--(d) that will imply Lemma \ref{pairs}. We will verify the claims in the remainder of the section.

Let $3 \leq i \leq 2d-1$ be odd and
let $a_i \in K_i(\mathcal{O}_F) $ be  chosen so that it generates $K_i(\mathcal{O}_F)$ modulo torsion. 
Let $[a_i]$ be the image of $a_i$ inside $ H_i(\GL_{\infty}(\mathcal{O}_F), \Z_{\ell})$;
as per our notation above, 
 $[a_i]^{(N)}$ is  its preimage under the isomorphism $H_i(\GL_N(\mathcal{O}_F),\Z_{\ell}) \rightarrow H_i(\GL_{\infty}(\OO_F),\Z_{\ell})$. 

We'll show that (for any $n$) there are  infinitely many places $v$, splitting the division algebra and with $q_v \equiv 1$ modulo $\ell^n$, 
and classes  $\xi_i \in H^i( \GL_N(\mathbf{F}_v),  \Z/\ell)$ with the property that:

\begin{itemize}
\item[(a)]  The image of $a_i$ in $K_i(\mathbf{F}_v)/\ell$ is nonzero. 
Call its image $b_i \in K_i(\mathbf{F}_v)/\ell$; therefore, 
 $b_i$  is a generator of $K_i(\mathbf{F}_v)/\ell$.
 
\item[(b)]  The pairing $ \langle \xi_i, [b_i]^{(N)} \rangle \neq 0$
where $[b_i]$ is  defined similarly to  $[a_i]$: it is the associated homology class
under $K_i(\mathbf{F}_v)  \rightarrow H_i (\GL_{\infty}(\mathbf{F}_v),\Z/\ell)$,
and correspondingly we have $[b_i]^{(N)} \in H_i(\GL_N(\mathbf{F}_v), \Z/\ell)$.

\item[(c)]  $ \langle \xi_3 \cup \dots  \cup \xi_{2d-1} ,  \left(  [b_3]  * \cdots *  [b_{2d-1}] \right)^{(N)} \rangle \neq 0$,
where * is  the  Pontryagin product on the homology of $\GL_{\infty}(\mathbf{F}_v)$. 
 
\item[(d)]  
Let $[Y(K)] \in H_{d^2-1}(Y(K), \Q)$ be the fundamental class of $Y(K)$. Then 
there exists $M \in \Q^*$ such that
the image $\iota_* [Y(K)] \in  H_*(\GL_N(\mathcal{O}_F) ,\Q) $   satisfies:
\begin{equation} \label{Mdef}  \iota_* [Y(K)] = M \cdot   \left(  [a_3] *[a_5] * \dots * [a_{2d-1}] \right)^{(N)}. \end{equation}
where on the right we have the Pontryagin product for $\GL_{\infty}(\mathcal{O}_F)$, and $\iota$
is as in \eqref{def of iota}.
\end{itemize}

Let us first see why (a) -- (d) implies  Lemma \ref{pairs}.
 Recall that we chose $\ell$ to not divide the  
 numerator or denominator of $M$, and also that the cohomology of $\GL_{\infty}(\mathcal{O}_F)$ and so also
 $\GL_N(\mathcal{O}_F)$ is $\ell$-torsion free in degrees $< d^2$;
therefore,  condition (d) implies an equality in $\Z_{\ell}$-homology:
 \begin{equation} \label{i unit} \iota_*[Y(K)] =   \left( \mbox{unit} \right) \cdot  \left(  [a_3] \dots * [a_{2d-1}] \right)^{(N)}.\end{equation}

Let $\pi$ be the projection from $\GL_N(\mathcal{O}_F)$ to $\GL_N(\mathbf{F}_v)$. 
Write $\Xi_i = \pi^* \xi_i \in H^*(\GL_N(\OO_F), \Z/\ell)$ 
and $\eta_i = \iota^* \Xi_i \in H^*(Y(K), \Z/\ell)$. We have then
$$ \langle \Xi_i, [a_i]^{(N)} \rangle  = \langle \pi^* \xi_i, [a_i]^{(N)} \rangle =  \langle \xi_i, \pi_* [a_i]^{(N)} \rangle = \langle \xi_i, [b_i]^{(N)} \rangle \neq 0$$ 
since $\pi_* [a_i] = [b_i]$.  Also, in a similar way, 
\begin{multline}  \langle \Xi_3 \cup \dots \cup \Xi_{2d-1},  \left( [a_3]* \dots  * [a_{2d-1}]  \right)^{(N)} \rangle \\
= \mbox{unit multiple of }  \langle \xi_3 \cup \cdots \cup \xi_{2d-1},   \left( [b_3] * \dots *  [b_{2d-1}] \right)^{(N)}  \rangle \neq 0\end{multline}
because the Pontryagin products and the stabilization maps are compatible with $\pi_*$. 
From this and \eqref{i unit} we get
 $$ \langle \eta_3 \cup \dots \cup \eta_{2d-1} , [Y(K)] \rangle = \langle \Xi_3 \cup \dots \cup \Xi_{2d-1}, \iota_* [Y(K)] \rangle \neq 0.$$

 But the Hecke-trivial cohomology of $Y(K)$ modulo $\ell$ is a free exterior algebra on generators in degrees $ 3, \dots, (2d-1)$. 
 Fix such generators -- call them $\nu_3, \nu_5, \dots$. 
Each $\eta_i$ is  also an element of this Hecke-trivial cohomology by Lemma \ref{lem:HTriv}. It follows that
 $$\eta_i = \mbox{ unit} \cdot \nu_i + \left(  \mbox{ product of $\nu_j$s with $j < i$}\right)$$
 because otherwise the cup product
 $\eta_3 \cup \dots \cup \eta_{2d-1}$ would be trivial.
 We conclude that, in fact, the 
 map
\begin{equation} \label{sarge000} H^*(\mathbf{G}(\mathbf{F}_v), \Z/\ell) \rightarrow H^*(Y(K), \Z/\ell)_{\triv} \end{equation}
 is {\em onto} as required.
 Therefore, to prove Lemma \ref{pairs} it is sufficient to prove (a)--(d) above.

 \subsection{Verification of (d) from \S \ref{abcd}} 
 Since we are supposing \eqref{mapN} to be a surjection in degrees up to $d^2$, it's enough to 
verify that there is a nonzero $M \in \mathbf{C}$ such that 
\begin{equation} \label{Fr}  \langle \iota_* [Y(K)], \omega \rangle =   M  \langle \left(  [a_3]* \dots * [a_{2d-1}] \right)^{(N)},  \omega \rangle\end{equation}
 whenever $\omega$ is an invariant differential form on $\GL_N(\C)$.
We compute in the case  (see \eqref{oinkoink} for the definition):
$$\omega = \Omega_J = \Omega_{j_1} \wedge \dots \wedge \Omega_{j_t}$$
where $J = \{j_1, \dots, j_t\}$, and  show that both
sides   are nonzero  if and only if $J=\{3, 5, \dots, 2d-1\}$. That is enough to prove \eqref{Fr}.

 The right hand side  of \eqref{Fr} equals 
\begin{equation} \label{KPd} \langle  [a_3]* \dots * [a_{2d-1}] , \omega^{(\infty)}  \rangle = \langle  [a_3] \otimes  \dots \otimes [a_{2d-1}], \mathrm{coproduct}(\omega^{(\infty)}) \rangle,\end{equation}
 where we allow ourselves to write $\omega^{(\infty)} \in H^*(\GL_{\infty}(\mathcal{O}_F), \C)$ for the  stabilization of the  cohomology class corresponding to $\omega$. 
 Now $\coproduct(\omega^{(\infty)} ) $ 
 is the product
 of various terms of the shape
  $$  \left( \Omega_{j_1} \otimes 1 \otimes 1 \otimes \cdots 1 + 1 \otimes \Omega_{j_1} \otimes 1 \otimes \dots \otimes 1 + \cdots \right)^{(\infty)} $$
 and from this we see (just for degree reasons) that the term on the right of \eqref{KPd} vanishes if there is even one $j_i$ larger than $2d-1$
 or one $j_i$ equal to $1$.  
 So  $J \subset \{3, \dots, 2d-1\}$. Again, for degree reasons, equality must hold. That shows the right-hand side is zero unless  $J=\{3, 5, \dots, 2d-1\}$.
 When $J=\{3, 5, \dots, 2d-1\}$, the right-hand side becomes 
 $$ \langle a_3, \Omega_3 \rangle \cdot  \langle a_5, \Omega_5 \rangle \dots  \langle a_{2d-1}, \Omega_{2d-1} \rangle$$
 and each factor $\langle a_j, \Omega_j \rangle$ is  nonzero: this is the nontriviality of the Borel regulator.

Now let us examine the left-hand side of \eqref{Fr}, which equals $\langle [Y(K)], \iota^* \Omega_{j_1} \wedge \dots \wedge \iota^* \Omega_{j_d} \rangle$.  It's easy to see that $\iota^* \Omega_1$ vanishes. We also  claim that $\iota^* \Omega_j$ must vanish for $j > 2d-1$.     Indeed we claim that $\iota^* \Omega_j$, which defines an invariant form on $Y(K)$
and thus corresponds by \eqref{raru} to an invariant differential form in the cohomology of $\mathrm{SU}_d$, is {\em primitive} as such.
  For that consider this diagram: 
   \begin{equation}
 \xymatrix{
 \mbox{invariant forms on $\SL_N(\C)/\mathrm{SU}_N$} \ar[r] \ar[d] &  \mbox{invariant forms on $\mathbf{G}(\C)/SU_d$} \ \ar[d] \\
\mbox{$\mathrm{SU}_N$-invariant forms on $i \mathfrak{su}_N$} \ar[r]\ar[d]  &  \mbox{$SU_d$-invariant forms on $i \mathfrak{su}_d$} \ar[d] . \\
H^*(\mathrm{SU}_N) \ar[r]  & H^*(SU_d).
 }
\end{equation}
where the various maps of groups arise from the map $\iota$ of \eqref{iotadef}, and the top vertical
maps arise by restriction to the tangent space of the identity coset. 
In other words, the element of $H^*(\mathrm{SU}_d)$ corresponding to $\iota^* \Omega_j$
is just the pull-back of the element of $H^*(\mathrm{U}_N)$ corresponding to $\Omega_j$
under the group homomorphism $$\varphi: SU_d \rightarrow U_N$$ 
induced by $\iota$. In particular, $\varphi^* \Omega_j$ is also primitive. 

 This shows that the left-hand side of \eqref{Fr} vanishes unless $J=\{3, 5, \dots, 2d-1\}$.  We must still check that it is actually nonvanishing in this case.
For this, we must show that $\iota^*( \Omega_3 \wedge \dots \wedge \Omega_{2d-1}) $
is nonvanishing, equivalently that
$$ \varphi^* \Omega_3 \wedge \dots \wedge \varphi^* \Omega_{2d-1}$$
is a nonvanishing element of the cohomology of $\mathrm{SU}_d$. Since each $\varphi^* \Omega_j$ is primitive it is enough
to see that they are all nonzero.    The $N$-dimensional representation of $\mathrm{SU}_d$ defined by $\varphi$ is 
 isomorphic to the sum of many copies of the standard representation of $\mathrm{SU}_d$. Now one can just compute explicitly with \eqref{oinkoink}.

 \subsection{Verification of  (b) and  (c) from \S \ref{abcd}}
 In words, what we have to do is produce elements $\theta_i \in H^i(\GL_{\infty}(\mathbf{F}_v), \Z/\ell)$
 for each odd $3 \leq i \leq 2d-1$; these $\theta_i$ should  detect (pair nontrivially with) a generator  of $K_i(\mathbf{F}_v)$, and the cup product
 $\theta_3 \cup \dots \cup \theta_{2d-1}$ 
should detect the Pontryagin product of  the homology classes associated to those generators.   Then we 
 may take $\xi_i = \theta_i^{(N)}  \in H^i(\GL_N(\mathbf{F}_v),\Z/\ell)$.

 Quillen shows a natural choice for $\theta_i$:  an equivariant Chern class derived from the standard representation of $\GL_N$.  
 In other words, write $G = \GL_N(\mathbf{F}_v)$ and write $\Gamma$ for the Galois group of $\mathbf{F}_v$.   
 The standard representation of $\GL_N$ can be considered a $G$-equivariant vector bundle on $\mathrm{Spec}(\mathbf{F}_v)$, and thus we get a Chern class
  $$c_{2i} \in H^{2i}_{G, \et}(\Spec \mathbf{F}_v , \Z/\ell(i)) \rightarrow H^1_{\et}(\Spec \mathbf{F}_v,  \Z/\ell(i)) \otimes H^{2i-1}(G, \Z/\ell)$$
  where on the left we have equivariant etale cohomology, and on the right we have usual etale cohomology; the arrow is explicated
  in  \cite[II.1.2, Lemma 1]{Soule}. 
   We may identify the {\'e}tale cohomlogy of $\mathbf{F}_v$ with the (continuous) group cohomology of $\Gamma$, and here $\Gamma$ 
  is acting trivially on $\mu_{\ell}$, and so $H^1(\Spec \ \mathbf{F}_v, \Z/\ell(i)) = H^1(\Gamma, \Z/\ell(i))$ is identified with simply $\mu_{\ell}^{\otimes i}$.
 Fix a generator  $\alpha$ for $\mu_{\ell}$. 
  Thus the image of $c_{2i}$ is of the form  $\alpha^i \otimes \theta_{2i-1}$,
  for some $\theta_{2i-1} \in H^{2i-1}(G, \Z/\ell)$. 
    Similarly we can map $c_{2i}$
  into $H^0(\Gamma, \Z/\ell(i)) \otimes H^{2i}(G, \Z/\ell)$; in that way we get a class $\theta_{2i}' \in H^{2i}(G, \Z/\ell)$ so that the image of $c_{2i}$
  is $\alpha^i \otimes \theta_{2i}'$. 
  The class $c_{2i}$ gives a morphism ( the ``Soul{\'e} Chern class'')
\begin{equation} \label{Soule} s: H_{2i-1}(G, \Z/\ell) \rightarrow H^1(\Gamma, \Z/\ell(i))\end{equation}
  which sends $\lambda \in H_{2i-1}$ to $\langle  \theta_{2i-1}, \lambda \rangle \alpha^{\otimes i}$. 

These constructions are ``stable'' under increasing $N$ -- see \cite[p. 257]{Soule} -- 
and so we can consider $\theta_{2i-1}, \theta_{2i}'$ as classes in $H^*(\GL_{\infty}(\mathbf{F}_v), \Z/\ell)$,
and the Soul{\'e} map as a map $H_{2i-1}(\GL_{\infty}(\mathbf{F}_v), \Z/\ell) \rightarrow H^1(\Gamma, \Z/\ell(i))$.

Now Soul{\'e} Chern class is  known to be surjective when precomposed with $K_{2i-1} \rightarrow H_{2i-1}$ (\cite[Prop. 5, p 284]{Soule})
(so long as $i < \ell$ -- true by assumption on $\ell$). 
This immediately verifies property (b), that is to say if we fix a generator $b_{2i-1}$
for $K_{2i-1}(\mathbf{F}_v)$ with associated homology class $[b_{2i-1}]$  we have
$$ 0 \neq s([b_{2i-1}]) = \langle \theta_{2i-1}, [b_{2i-1}] \rangle \alpha^{\otimes i}$$
so $\langle \theta_{2i-1}, [b_{2i-1}] \rangle \neq 0$ as desired.

    To compute e.g. $\langle \theta_3 \cup  \theta_5 \cup \theta_7,  [b_3]  *  [b_5] * [b_7] \rangle $ we rewrite it as 
\begin{equation} \label{todo}   
= \langle \coproduct(\theta_3) \coproduct(\theta_5) \coproduct(\theta_7), b_3 \otimes b_5 \otimes b_7 \rangle  \end{equation}

Let us note that each $[b_i]$ is primitive in homology -- i.e. $\langle [b_i], \alpha \cup \beta \rangle = 0$ if $\alpha, \beta$ are cohomology classes both in positive degree. 
This is because $[b_i]$ comes from the image of the  Hurewicz map $\pi_i \rightarrow H_i$, so we can just pull back to the sphere $S^i$ and compute.

  Quillen has shown \cite[Proposition 2]{QuillenFp} (see also \cite[Remark 2, p.569]{QuillenFp})   that  the coproduct of (e.g.) $\theta_7$ equals  $$\coproduct( \theta_7)= \underbrace{\theta_0'}_{1} \otimes  \theta_7 + \theta_2' \otimes \theta_5 +  \theta_4' \otimes \theta_3 +\theta_6' \otimes \theta_1 + \mbox{symmetric terms}.$$
Thus, when we take the product $\coproduct(\theta_3) \coproduct(\theta_5) \coproduct(\theta_7)$,
we get a sum of several terms; because of the 
  {\em primitivity} of $[b_i]$ just noted, 
  the only terms that contribute to \eqref{todo} will be
  those coming from
{\small   $$(1 \otimes 1 \otimes \theta_3 + 1 \otimes \theta_3 \otimes 1+ \theta_3 \otimes 1 \otimes 1) (1 \otimes 1 \otimes \theta_5 + 1 \otimes \theta_5 \otimes 1+ \theta_5 \otimes 1 \otimes 1) 
 (1 \otimes 1 \otimes \theta_7 + 1 \otimes \theta_7 \otimes 1+ \theta_7 \otimes 1 \otimes 1) $$}
  and the only term of these with degree $3,5,7$ in the first, second and third factors is $\theta_3 \otimes \theta_5 \otimes \theta_7$. 
Therefore, 
\begin{align*} \langle \theta_3 \cup  \theta_5 \cup \theta_7,  [b_3]  *  [b_5] * [b_7] \rangle &=&  \langle  \theta_3 \otimes \theta_5 \otimes \theta_7, [b_3] \otimes [b_5] \otimes [b_7] \rangle  \\ &=& \prod \langle \theta_{3},[b_3] \rangle \langle \theta_5, [b_5] \rangle \langle \theta_7, [b_7] \rangle \neq 0 \end{align*}
which gives (c) in the case $d=4$; the general case is the same.

   \subsection{Verification of (a) from \S \ref{abcd}} \label{sec:Chebo}
Write $\mathcal{O}' = \mathcal{O}_F[1/\ell]$. 
The  Soul{\'e} maps from \eqref{Soule}  fit in a commutative diagram 
 \begin{equation}
 \xymatrix{
 K_{2i-1}(\mathcal{O}')  \ar[r] \ar[d] & K_{2i-1}(\mathbf{F}_v) \ar[d] \\
 H^1(\Gamma, \Z_{\ell}(i)) \ar[r] & H^1(\mathrm{Gal}(\overline{\mathbf{F}_v}/\mathbf{F}_v), \Z_{\ell}(i))  .
 }
\end{equation}
where $\Gamma$ is  now the Galois group of  the maximal unramified extension  of $\mathcal{O}'$ ; note that  
 $\mathcal{O}_F \rightarrow \mathcal{O}'$ induces an isomorphism on $K_{2i-1}$ for $i > 1$ (see \cite[Theorem 4.6]{Weibel}). 
The right-hand vertical arrow is an isomorphism for $i <\ell$ (Soul{\'e}, {\em loc. cit.}) , and   the left-hand vertical arrow is a surjection (see \cite{Kahn}).

The map $H^1(\Gamma, \Z_{\ell}(i))/\ell \rightarrow H^1(\Gamma, \mathbf{F}_{\ell}(i))$ 
is an injection.  Choose an element $a_{2i-1} \in K_{2i-1}(\mathcal{O}')$ whose image in $ H^1(\Gamma, \mathbf{F}_{\ell}(i))$
is nonzero (note that $H^1(\Gamma, \mathbf{Z}_{\ell}(i))/\ell$ is nonzero, by computing Euler characteristic). 
A nontrivial class in $H^1(\Gamma, \mathbf{F}_{\ell}(i))$ is represented by a nontrivial extension 
\begin{equation} \label{exli} \mathbf{F}_{\ell}(i) \rightarrow M \rightarrow \mathbf{F}_{\ell},\end{equation} in other words, by a homomorphism
$\Gamma \rightarrow \GL_2(\Z/\ell)$ of the form 
 $$ \left( \begin{array}{cc} \omega^i & \rho \\ 0 & 1 \end{array}\right) $$ 
 Note that the image of this homomorphism must have size divisible by $\ell$; 
for otherwise the extension \eqref{exli} splits. By Chebotarev density, we may find infinitely many Frobenius elements $\mathrm{Frob}_v$ which map
 to 
  $$ \left( \begin{array}{cc} 1 & 1 \\ 0 & 1 \end{array}\right) $$ 
 which  implies that the restriction of the extension class to $H^1(\mathbf{F}_v, \mathbf{F}_{\ell}(i))$ is nontrivial.
 
 In other words, the image of $a_{2i-1}$ in $H^1(\mathrm{Gal}(\overline{\mathbf{F}_v}/\mathbf{F}_v), \mathbf{F}_{\ell}(i)) $,
 and so also in $K_{2i-1}/\ell$, 
 is nontrivial for finitely many $v$.  This proves our assertion.

\section{Setup for patching}  \label{Patching} 

We continue in a global setting, but now turn to the study of tempered cohomology. This study 
will occupy most of the remainder of the paper (\S \ref{Patching} -- \S \ref{reciprocity}). 

 In the current section (\S \ref{Patching}) we will set up  the various assumptions needed,
 and in the next section (\S \ref{Patching2}) we will use the relationship between
 the Taylor--Wiles method and the derived Hecke algebra (outlined in \S \ref{sec:explication})
 to prove our target theorem, Theorem \ref{mindegreeproof}: it says the global derived Hecke algebra
actually is big enough to be able to account for the degree spread of cohomology.

In \S \ref{reciprocity} we will explain how to index elements of the global derived Hecke algebra by a Selmer group,
and use this to formulate our main Conjecture \ref{mainconjecture}.

A few apologies are in order:
 \begin{itemize}

\item 
 We   switch notation slightly, working with mod $p^n$ coefficients rather than mod $\ell^n$, to better
make contact with the standard presentations of Galois representations and the Taylor-Wiles method.

\item We have made no attempt to optimize the method for small primes,
and in particular make rather strong assumptions; in particular, we assume both that 
we are in the ``minimal case'' of formally smooth local deformation rings, and that
the Hecke ring at base level is isomorphic to $\Z_p$ (no congruences). 

\end{itemize}

\subsection{Assumptions}  \label{Sec:assumptions}

 \begin{enumerate}
 
 \item Our general notations are  as in \S \ref{sec:notn}, but we now specialize
 to the case that the number field $F$ is $\Q$, and that 
  $\mathbf{G}$ is  a {\em simply connected},  semisimple $\Q$-group.  We assume that $\mathbf{G}$ is split,  and we fix a Borel subgroup $\mathbf{B}$ and a maximal torus $\mathbf{A}$ contained inside $\mathbf{B}$.

\item 
We fix a level structure $K_0 \subset   \mathbf{G}(\Afinite)$, the ``base level.''   \index{$K_0$} \index{$Y(1)$}
   We write $$Y(1)=Y(K_0)$$ for the arithmetic manifold (\eqref{YKdef}) of level $K_0$.
   We sometimes refer to this as the ``level $1$'' arithmetic manifold even though it is not literally so. 

   \item Let $\Pi$ be a tempered  cohomological automorphic cuspidal representation for $\mathbf{G}$, factorizing
as $\Pi = \Pi_{\infty} \otimes \Pi_f$ over archimedean and finite places.  \index{$\Pi$}   
We suppose that $\Pi^{K_0} \neq 0$, so that $\Pi$ actually contributes to the cohomology of $Y(1)$. 
      
   \item We write $\mathbb{T}_{K_0}$ for the Hecke algebra at level $K_0$. It will be convenient to follow the definition of
   \cite{KT} and define this in a derived sense:
   Consider the chain complex of $Y(K_0)$, with $\Z_p$ coefficients,  as an object in the derived category of $\Z_p$-modules; each
  (prime to the level, and to $p$) Hecke operator gives an endomorphism of this object. Define 
   $\mathbb{T}_{K_0}$ to be the  $\Z_p$-algebra generated by such endomorphisms.  
This has the advantage that $\mathbb{T}_{K_0}$ acts on cohomology with coefficients
in any $\Z_p$-module. 

We may similarly form the Hecke algebra $\mathbb{T}_K$ at a deeper level $K \subset K_0$; 
unless specified, it will be generated only by Hecke operators at good primes for $K$. 
   
  \item  We shall suppose that the coefficient field of $\Pi$ is $\Q$, for simplicity -- by this, we mean that  the eigenvalues of Hecke correspondences at 
 all good places for $K_0$ lie in $\Q$,
  or equivalently the underlying representation $\Pi_v$ has a $\Q$-rational structure for  all such $v$. %
Under this assumption,  $\Pi$ gives
   rise to a ring homomorphism 
\begin{equation} \label{oort}  \mathbb{T}_{K_0} \rightarrow \Z_p.\end{equation}

\index{$\Pi$} \index{$R$} 
 
 \item 
Let $T$ be the set of ramified places for $\Pi$, together with
any places at which $K_0$ is not hyperspecial.  \index{$T$ (set of places)}

 \item 
 Write
 $$k =\Z/p\Z$$ and fix an algebraic closure $\overline{k}$ for $k$
 where  $p > 5$ is a prime such that:
  
 \begin{itemize}
 \item[(a)]  $H^*(Y(K_0), \Z)$ is $p$-torsion free.\item[(b)] $p$ does not 
divide the order of the Weyl group of $\GG$, and also $p \notin T$, i.e. $p$ is not a bad 
place for $\Pi$ or $K$. 

      \item[(c)]  
``No congruences between $\Pi$ and other forms at level $K_0$:''

Consider the composite homomorphism
\begin{equation} \label{chi_def}  \chi: \mathbb{T}_{K_0} \rightarrow \Z_p \rightarrow k,\end{equation}
where the first map is the action on $\Pi$, and the second map the obvious one.
Let $\mathfrak{m} = \ker(\chi) $ be its kernel.  \index{$\mathfrak{m}$} 
We  shall require that  the induced map of completions is an {\em isomorphism}:  
\begin{equation} \label{nocong0} \mathbb{T}_{K_0, \mathfrak{m}} \stackrel{\sim}{\longrightarrow} \Z_p\end{equation}
and moreover we shall assume the vanishing of   homology after completion at the maximal ideal $\mathfrak{m}$:
 \begin{equation} \label{nocong2} H_j(Y(1), \Z_p)_{\mathfrak{m}} =  0,   j  \notin  [\BWq, \BWq+\delta].  \end{equation} 
 Observe that in favorable situations \eqref{nocong0} implies \eqref{nocong2}, and both should be true for all large enough $p$ -- 
  see the Remark below
for a further discussion.    Informally, \eqref{nocong0} and \eqref{nocong2} enforce
that there are no congruences, modulo $p$, between $\Pi$ and other cohomological forms at level $K_0$. 

Note that the definitions of $\chi, \mathfrak{m}$ make sense at any level.  Thus we use the notation $\chi, \mathfrak{m}$ sometimes for the corresponding notions for other level structures $Y(K)$,
where $K \subset K_0$.

\end{itemize}
 
\item We put $$\ramprimes  = T \cup \{p\},$$ the collection of all primes that we have to worry about.  \index{$S$ (set of ramified places)}

\end{enumerate}

 \begin{remark}
We expect that  7(c) should be automatically valid for $p$ sufficiently large; in practice, for the purposes of this paper,
it is not an onerous assumption (the minimal level conditions, enforced in (e) of \S \ref{GaloisAss}, is more restrictive).

We give a proof that 7(c) is valid for all large enough $p$,  for $\mathbf{G}$ an  inner form of $\mathrm{SL}_n$ such that $Y(K)$ is compact.
It is likely this can be generalized to other settings, with more work.    In what follows, denote by $\mathbb{T}_{K_0}$
 the Hecke algebra defined as above, but with $\Z$ coefficients; this is easily seen to be finitely generated over $\Z$.   Firstly,  the algebra $\mathbb{T}_{K_0} \otimes \C$ is semisimple,
 because it acts faithfully on $H^*(Y(K_0), \C)$ and this action is semisimple (there is an invariant metric on harmonic forms). 
 Thus, for all large enough primes $p$, $\mathbb{T}_{K_0}$ is {\'e}tale over $\mathbf{Z}$ and thus \eqref{nocong0} must be valid. 
If  the homology in \eqref{nocong2} is nonvanishing, then 
there exists an eigenclass for $\mathbb{T}_{K_0}$  on $H_j$ whose associated character
factors through $\mathbb{T}_{K_0, \mathfrak{m}}$;  by \eqref{nocong0}, this character must coincide with the action of $\mathbb{T}_{K_0}$
 on $\Pi$.  In other words,  the Hecke eigensystem associated to $\Pi$ occurs in degree $j$.  
 This eigensystem corresponds to an automorphic representation $\Pi'$ such that $\Pi'_v = \Pi_v$
 for almost all $v$.  By the strong multiplicity one theorem for $\mathbf{GL}_n$, this implies that $\Pi'_{\infty}$ is tempered cuspidal, and then
 it has nonvanishing $(\mathfrak{g}, K)$ cohomology only in degrees $[\BWq, \BWq+\delta]$. 
  \end{remark}

That's the basic setup; now 
for Galois representations. 

\subsection{Assumptions about Galois representations; the deformation ring $\defring$ } \label{GaloisAss}
 We will make assumptions very close to \cite[Conjecture 6.1]{GV}.
 We briefly summarize them and refer the reader to  \cite{GV} for full details:

Let $K \subset K_0$ be a deeper level structure, $\mathbb{T}_K$ be as above, and let $\mathbb{T}_K \rightarrow k$
be a character. {\color{\changecolor} We require that there exist a   Galois representation
$$\mathrm{Gal}(\overline{\Q}/\Q) \rightarrow G^{\vee}(\bar{k})$$
satisfying the usual unramified compatibility, see (a) below.  }

Moreover, for the specific character $\mathbb{T}_K \rightarrow k$ 
  as 
in  \eqref{chi_def}, i.e. 
the map associated to  the fixed automorphic representation $\Pi$,
reduced modulo $p$, we require more precise statements:
Let $\mathfrak{m}$ be the kernel of $\mathbb{T}_K \rightarrow k$, and $\mathbb{T}_{K,\mathfrak{m}}$ the completion of $\mathbb{T}_K$ at $\mathfrak{m}$. 
We require there to exist a Galois representation\footnote{Note that the assumption that there is a $\mathbb{T}$-valued Galois representation,
rather than a weaker notion such as a determinant, is not reasonable unless one has a condition like ``residual irreducibility.'' In our case, however, we are assuming that
the residual representation $\rhobar$ has very large image anyway -- see (b) below.} 
$$ \tilde{\rho}: \mathrm{Gal}(\overline{\Q}/\Q) \longrightarrow   G^{\vee}(\mathbb{T}_{K, \mathfrak{m}})$$
 with the following properties: 

\begin{itemize}
\item[(a)]
(Unramified compatibility): Fix a representation $\tau$ of $G^{\vee}$.
For all primes $q$ not dividing the level of $K$, 
the representation $\rho$ is unramified at $q$, and  the action of $\mathrm{trace}(\tau \circ \tilde{\rho})(\mathrm{Frob}_q) \in   \mathbb{T}_{K, \mathfrak{m}}$
coincides with the image of  the (Satake)-associated Hecke operator 
$T_{q, \tau}$.

 \item[(b)]  Let $\rhobar$ be the reduction modulo $p$ of $\tilde{\rho}$, so that
$$\rhobar: \Gal(\overline{\Q}/\Q) \rightarrow G^{\vee}(k).$$ 
 Then $\rhobar$ has big image:  when restricted to the Galois
 group of $\Q(\zeta_{p^{\infty}})$, 
 the image of $\rhobar$   contains the image of the $k$-points of the simply connected cover of $G^{\vee}$.
 
 \item[(c)] (Vague version: see \cite[Conjecture 6.1]{GV} for precise formulation): There is a reasonable notion of ``crystalline at $p$'' representation into $G^{\vee}$,
 and the representation $\tilde{\rho}$ is ``crystalline at $p$.'' 
 
 \item[(d)] (Vague version: see \cite[Conjecture 6.1]{GV} and references therein for precise version): The representation
 $\tilde{\rho}$ satisfies the expected local constraints (``local--global compatibility'') when restricted to $\Q_{q}$; 
 here $q$ is a Taylor--Wiles prime (\S \ref{TWprimes stuff}) that divides the level of $K$.  We {\em also} assume a natural version of local global compatibility
 at Iwahori level $Y_0(q)$, formulated  before  Lemma \ref{needaref}. 
 
 \item[(e)] ```All local deformation rings are all formally smooth:''   we suppose that   $$H^0(\Q_q, \Ad \rhobar) = H^2(\Q_q, \Ad \rhobar)=0 \ \ \ \mbox{ for all } q \in \ramprimes = T \cup \{p\}.$$
This means that the local deformation ring of $\rhobar$ at primes $q \in T$  is isomorphic to $\Z_p$,
{\em and} the  local deformation ring of $\rhobar$ at $q=p$ is formally smooth.  This keeps our notation as light as possible. (We also use the formal smoothness
to squeeze the most out of the Taylor--Wiles method, but 
probably one can get something without it.)
 \end{itemize}
 
 In particular, the natural map $\mathbb{T}_{K_0,\mathfrak{m}} \rightarrow \Z_p$  of \eqref{nocong0} gives rise to a Galois representation
 valued in $G^{\vee}(\Z_p)$, which we shall just call $\rho$: \index{$\rho$} 
\begin{equation} \label{Zplift} \rho: \Gal(\overline{\Q}/\Q) \longrightarrow G^{\vee}(\Z_p),\end{equation}
 which of course lifts the residual representation:
 $$\rhobar: \Gal(\overline{\Q}/\Q) \longrightarrow G^{\vee}(k).$$ \index{$\rhobar$} 
 \index{$\defring$} 
 Let $\defring$ be the universal  crystalline deformation ring of $\rhobar$, allowing ramification only at the set $\ramprimes$.    
Good references for deformation rings  are \cite{CHT} 
or \cite{AWS}.

Let us now set up notations for Taylor--Wiles primes.

\subsection{Taylor--Wiles primes and auxiliary level structures} \label{TWprimes stuff}
\index{Taylor--Wiles prime}
A {\em Taylor-Wiles prime} of level $n$ is a prime $q$ (we will also occasionally use the letter $\ell$), not dividing the level of $K_0$, 
such that:
\begin{itemize}
\item[-] $p^n$ divides $q-1$, and 
\item[-] $\rhobar(\Frob_q)$ is conjugate to a strongly regular element of $T^{\vee}(k)$.
\end{itemize}
Here an element $t \in T^{\vee}(k)$ is {\em strongly regular} if its centralizer inside $G^{\vee}$ is equal to $T^{\vee}$. \index{strongly regular elememnt}

We are really interested in systems of such primes, and it is useful to keep track of the  strongly regular element as part of the data.
Fix once and for all a sufficiently large integer $s$. We will work with collections of such primes of cardinality $s$, \index{$s$ (size of Taylor--Wiles datum)} 
which we call Taylor--Wiles data:

\begin{itemize}

\item[-]   A {\em Taylor-Wiles datum} of level $n$ is   a set of primes $Q_n = (q_1, \dots, q_s)$  together with strongly regular elements $(\Frob_{q_1}^T, \dots, \Frob_{q_s}^T) \in T^{\vee}(k)$  \index{$\Frob_q^T$}
 such that 
\begin{itemize}
\item[-]    $p^{n}$ divides $q_i-1$, and 
\item[-]  $\rhobar(\Frob_{q_i})$ is  conjugate to  $\Frob_{q_i}^T$.
\end{itemize} 
  
 We will usually use $q$ to denote a typical element of a Taylor--Wiles set of primes, but occasionally we will also use the letter $\ell$.
Note that the set of possible choices for each $\Frob_{q_i}^T$ has size $|W|$, the order of the Weyl group.

\item[-] 
Let $r$ be the rank of the maximal torus $\TT$. Set \begin{equation}  \label{Tsdef}  R=rs, \ \ \ T^{\vee}_s = (T^{\vee})^s, \ \ \ \  W_s = W^s \end{equation} 
 thus $T^{\vee}_s$ is a torus and $R$ is the rank of $T^{\vee}_s$, and $W_s$ acts on $T^{\vee}_s$.  \index{$W_s$}

\item[-]  Level structures: 
If $q \notin \ramprimes$ is prime, we denote by  \index{$Y_0(q)$} 
\begin{equation} \label{Y0ldef} Y_0(q) \rightarrow Y(1)\end{equation} 
the covering obtained by adding Iwahori level structure at $q$, i.e. we replace $K_0$
by the preimage of a Borel subgroup under $K_0 \rightarrow \GG(\F_{q})$. 
Similarly, we get $Y_1(q)$ by taking the preimage of a unipotent radical of a Borel subgroup under the same mapping.  \index{$Y_1(q)$}
The covering $Y_1(q) \rightarrow Y_0(q)$ is ``Galois,'' with Galois group $ \mathbf{A}(\mathbf{F}_{q}) \simeq (\mathbf{F}_{q}^{\times})^r$.

  Suppose that  $p^{n}$ divides $q-1$. In that case,  define $Y_1(q,n) $ to be  \index{$Y_1(q,n)$} 
the unique subcovering of $Y_1(q) \rightarrow Y_0(q)$ such that
the covering $Y_1(q, n) \rightarrow Y_0(q)$
has covering group $(\Z/p^n)^{r}$. In summary:
\begin{equation} \label{Y1def}  \overbrace{ Y_1(q) \rightarrow  \underbrace{ Y_1(q, n) \rightarrow Y_0(q)}_{\mathbf{A}(\mathbf{F}_{q})/p^n \simeq (\Z/p^n)^{r}}}^{ \mathbf{A}(\mathbf{F}_{q}) \simeq (\mathbf{F}_{q}^{\times})^{r}} \end{equation}

\item[-] For a Taylor--Wiles  datum $Q_n$ of level $n$, we  let   \index{$Y_1^*(Q_n)$} 
 $Y_1^*(Q_n)$ be the fiber-product, over $Y(1)$, of all the coverings $Y_1(q_i,n) \rightarrow Y(1)$. 
 Similarly we define $Y_0(Q_n)$. \index{$Y_0(Q_n)$} 
Therefore, $Y_1^*(Q_n) \rightarrow Y_0(Q_n)$ is Galois;  \index{$T_n$}
we write 
\begin{equation} \label{idontcare} T_n = \left(  \mbox{ Galois group of  $Y_1^*(Q_n) \rightarrow Y(1)$ }  \right) = \prod_{q \in Q_n} \mathbf{A}(\mathbf{F}_q)/p^n,\end{equation}
thus we have (non-canonically) 
\begin{equation} \label{identifT} T_n   \simeq (\Z/p^n)^R.\end{equation}

   As a shorthand we write
$$H_*(Q_n, S) := H_*(Y_1^*(Q_n), S)$$
where $S$ is a ring of coefficients.

 \end{itemize}

\subsection{Rings of diamond operators} \label{sec:diamond} 
We now set up the rings that are generated by ``diamond operators,'' i.e. the deck transformation groups of our various coverings $Y_1^*(Q_n) \rightarrow Y_0(Q_n)$. 
Continuing  with the notation of the prior section, put   \index{$S_n$} 
\begin{equation} \label{sndef} S_n  = \Z/p^n[T_n], \ \ S_n' =  \Z_p [T_n], \end{equation}
so these act on  the $\Z/p^n$ and $\Z_p$-valued chain complex of $Y_1^*(Q_n)$. 
Fixing, as in \eqref{identifT} , an isomorphism of $T_n$ with $(\Z/p^n)^R$, we can identify this ring as follows.
\begin{equation} \label{explicit} S_n \simeq    \Z/p^n[x_1, \dots, x_R]/( (1+x_i)^{p^n}-1). \end{equation} 
where $x_i = [e_i]- 1$, $e_i$ being a generator of the $i$th factor $\Z/p^n$ under the isomorphism \eqref{identifT}.   Recall $R=rs$ as in \eqref{Tsdef}.

Finally we form  a ``limit ring''  \index{$\Sinf$} 
\begin{equation} \label{Sinfdef} \Sinf := \Z_p[[x_1, \dots, x_R]]\end{equation}

The presentation \eqref{explicit} gives rise to obvious maps
 $ \Sinf \twoheadrightarrow S_n$ and  (compatible) augmentations of $\Sinf$ and $S_n$ to $\Z_p$ and $\Z/p^n$ respectively, carrying all the $x_i$ variables to zero; we denote by $\Iinf$ and $I_n$ the corresponding kernels,
so that $\Sinf/\Iinf \simeq \Z_p$ and $S_n/I_n \simeq \Z/p^n$. 
We need the following easy Lemma:
 \begin{lemma*}
Let notations be as above. The natural map
\begin{equation} \label{toinfinity} \Ext^*_{S_n }(\Z/p^n, \Z/p^n) \rightarrow \Ext^*_{\Sinf/p^n}(\Z/p^n, \Z/p^n) \end{equation}
(change of rings) is {\em surjective}. Also,  the natural map
\begin{equation} \label{corfunct} \Ext^*_{\Sinf/p^n}(\Z/p^n, \Z/p^n) \rightarrow \Ext^*_{\Sinf}(\Z_p, \Z/p^n)\end{equation} 
 (change of ring, and functoriality of $\Ext$ in the first argument) is an isomorphism. 
 
\end{lemma*}

 Recall that the change of ring map  $\Ext_B \rightarrow \Ext_A$ induced by a ring map $A \rightarrow B$  can be realized by thinking of $\Ext_B$ in terms of extensions of $B$-modules, and then
 just regarding it as an extension of $A$-modules.  
 
\proof 
\iflongproof
For \eqref{toinfinity}, it is sufficient to check surjectivity on  $\Ext^1$ because the right hand side
is generated by $\Ext^0, \Ext^1$: it can be computed  to be an exterior algebra using a Koszul resolution, see Lemma \ref{KoszulExtComputation}.
  Now ``morally speaking''  this is because it is the pullback on  group $H^1(-, \Z/p^n)$ induced
by $\Z_p^R \rightarrow (\Z/p^n)^R$, but  this is not a real proof  (at least without discussing the relation between cohomology of profinite groups, 
and $\Ext$s over the corresponding ``completed'' group algebra).

So we just compute both sides by using the  compatible (under $\Sinf/p^n \rightarrow S_n$) sequences
$I_n \rightarrow S_n \rightarrow \Z/p^n$ and $\Iinf/p^n \rightarrow \Sinf/p^n \rightarrow \Z/p^n$;
taking homomorphisms into $\Z/p^n$, we get 
\begin{equation} \label{coocoo} \Ext^1_{S_n}(\Z/p^n, \Z/p^n) \simeq \Hom_{S_n}(I_n, \Z/p^n) \left( \simeq \Hom(I_n/I_n^2, \Z/p^n) \right) \end{equation}
and similar.    (Note that the image of $\Hom_{S_n}(S_n, \Z/p^n)$ in $\Hom_{S_n}(I_n, \Z/p^n)$ is zero.) 

So we need to check that the map 
\begin{equation} \label{godot} \Iinf/p^n \rightarrow I_n\end{equation}
induces an isomorphism when we take $\Hom_{\Sinf}(-, \Z/p^n)$.

The homomorphisms $\Hom_{\Sinf}(I_n, \Z/p^n)$ are precisely given by homomorphisms
$\varphi: T_n \rightarrow \Z/p^n$, namely, we send $x_i = [e_i] -1 \in I_n$  to  $\varphi(e_i)$,
using the notation after \eqref{explicit}.  
Since $\Iinf/(p^n, \Iinf^2)$ is a free $\Z/p^n$ module on $x_1, \dots, x_{R}$,
it follows at once that \eqref{godot} induces an isomorphism as desired.

 Now we discuss \eqref{corfunct}, which is not hard but we spell it out, mainly to be clear because 
 I find change of rings confusing.  
 Let $\pi: P \rightarrow \Z_p$ be the Koszul  resolution of $\Z_p$
 as an $\Sinf$-module, and $\overline{P} = P/p^n$ its reduction mod $p^n$, 
 so that $\overline{\pi}: \overline{P} \rightarrow \Z/p^n$ is the Koszul resolution of $\Z/p^n$ as an $\Sinf/p^n$-module.  
 There is an identification $\Ext^*_{\Sinf/p^n}(\Z/p^n, \Z/p^n)$
 with the cohomology of $\Hom(\overline{P}, \Z/p^n)$; it sends a closed element $C \in  \Hom(\overline{P}, \Z/p^n[m])$ to the class 
  $\alpha \in \Ext^m_{\Sinf/p^n}(\Z/p^n, \Z/p^n)$
  represented by   the diagram $\Z/p^n \stackrel{\overline{\pi}}{\longleftarrow} \overline{P}  \stackrel{C}{\rightarrow} \Z/p^n[m]$
 in the derived category of $\Sinf/p^n$-modules (one can invert quasi-isomorphisms in the derived category).  
 
 Consider now the diagram  
 \begin{equation} \label{feeling_groovy0}
 \xymatrix{
\Z/p^n    & \overline{P}   \ar[l]^{\sim}_{ \overline{\pi}} \ar[r]^C  & \Z/p^n[m]   \\
\Z_p    \ar[u]^{A}  & \ar[l]^{\sim}_{\pi} P \ar[l] \ar[u]^B   \\
 }
 \end{equation}
 where $\sim$ means quasi-isomorphism; $A,B$ are the natural projections.

   Now, the image $\alpha'  \in \Ext_{\Sinf}(\Z_p, \Z/p^n)$ of $\alpha$, under the map \eqref{corfunct}, 
is represented by the map $C \circ \overline{\pi}^{-1} \circ A$  inside the derived category of $\Sinf$-modules, equivalently, by the composition
$C \circ B \circ \pi^{-1}$.  
In other words, the image $\alpha'$ of $\alpha$ is represented by the class $C' \in \Hom(P, \Z/p^n[m])$ obtained by pulling back
$ C \in \Hom(\overline{P}, \Z/p^n[m])$ via $P \rightarrow \overline{P}$.  That  pullback induces an isomorphism of complexes 
$$\Hom_{\Sinf/p^n}(\overline{P}, \Z/p^n) \simeq \Hom_{\Sinf}(P, \Z/p^n),$$ and thus, passing to cohomology, the desired isomorphism of $\Ext$-groups. 
\fi
\qed

\subsection{ Adding level $q$ structure: The relationship between homology at level $Y(1)$ and level $Y_0(q)$ } \label{Morita}

Fix a prime $q$ that doesn't divide the level of $K_0$, satisfying
$q \equiv 1$ modulo $p^n$. 
This section and the next two \S \ref{Morita}, \S \ref{HYln},   \S \ref{Lgcompat} all gather some ``standard'' properties of passing  between level $1$ and level $q$,
working with $\Z/p^n$-cohomology. 
We   remind the reader that ``level $1$'' does not literally mean level $1$, but just the base level $K_0$ at which we work.

In the current section, we discuss the relationship between the homology of $Y(1)$ and $Y_0(q)$; recall that $Y_0(q)$ was obtained
by adding Iwahori level at $q$ (see \eqref{Y0ldef}).

We can use the discussion of \S \ref{IwahoriHecke}:  
For $S=\Z/p^n$,  let $\mathrm{H}_I, \mathrm{H}_K$ and so on be the (usual, i.e. no ``derived''!) Hecke algebras for $(G_{q}, K_{q})$, with $S$-coefficients 
 and let $\mathrm{H}_{IK}, \mathrm{H}_{KI}$ be the bimodules previously defined in \S \ref{Iwahori}. 
 There are  natural maps
\begin{equation} \label{juin} H^*(Y(1),S) \otimes_{\mathrm{H}_{K}} \mathrm{H}_{KI} \longrightarrow H^*(Y_0(q),S),\end{equation}
$$   \ \ H^*(Y_0(q),S) \otimes_{\mathrm{H}_{I}} \mathrm{H}_{IK} \longrightarrow H^*(Y(1),S).
$$  
These maps are defined by the ``usual double coset formulas.'' More formally, we may identify  $\mathrm{H}_K$ with $\Hom_{SG_q} (S[G_q/K_q], S[G_q/K_q])$,  
$\mathrm{H}_{KI}$ with $\Hom_{SG_q}(S[G_q/K_q], S[G_q/I_q]$) and so forth;   one may then   proceed as in the discussion of
\S \ref{Arithmeticmanifolds}  to define the maps of \eqref{juin}. 
 
 Before we go further, let us formulate a natural notion of   ``local-global compatibility'' at level $Y_0(q)$: 
  the center of the Iwahori-Hecke algebra $\mathrm{H}_I$ at level $q$,
which is identified (\S \ref{Zdis}) with the Hecke algebra $\mathrm{H}_K$. We shall suppose:  \begin{quote} (Local-global compatibility at level $Y_0(q)$:) 
$\mathrm{H}_K$, identified with the center of $\mathrm{H}_I$ as just explained, acts on $H_*(Y_0(q), k)_{\mathfrak{m}}$ by means of the same (generalized)  eigencharacter $\mathrm{H}_K \rightarrow k$
by which $\mathrm{H}_K$ acts on $\Pi$.  
\end{quote} 

It is feasible that this assumption could be avoided entirely but since it is very likely to be proven along with the other, more essential,
local-global compatibility at Taylor-Wiles primes it seems harmless to assume it. 
 
 \begin{lemma}  \label{needaref}
 Let notation be as above,
and assume the local-global compatibility just mentioned. 
 {\color{\changecolor} Suppose that the prime $q$ is such that $\rhobar(\mathrm{Frob}_q)$ is conjugate
 to a strongly regular element of $T^{\vee}(k)$.}
Then the maps \eqref{juin} are   isomorphisms when  we localize
at a maximal ideal of $\mathrm{H}_K$ induced by $\mathfrak{m}$. 
In particular, if $q$ is part of a Taylor-Wiles datum of level $n$, then \eqref{juin} induces 
\begin{equation} \label{Twloc} H^*(Y_0(q), \Z/p^n)_{\mathfrak{m}}  \stackrel{\sim}{\longleftarrow} H^*(Y(1), \Z/p^n)_{\mathfrak{m}} \otimes_{\mathrm{H}_K} \mathrm{H}_{KI}.\end{equation}
  \end{lemma}

\proof 
The assertion that  \eqref{Twloc} is an isomorphism follows from the assertion that \eqref{juin} is an isomorphism by  ``localization at $\mathfrak{m}$'' (taking a little care because the Hecke algebra for $Y_0(q)$ and $Y(1)$ 
are not quite identically defined, the former omitting the prime above $q$. For a more detailed treatment of this point, please see \cite[Lemma 6.20]{KT}).   

%
%
%
%
%

So it is enough to prove the first assertion. 
Everything will be with  $\Z/p^n$ coefficients.  First  note that  the two natural ways of making a map
$$H^*(Y(1)) \otimes_{\mathrm{H}_K} \mathrm{H}_{KI} \otimes_{\mathrm{H}_I} \mathrm{H}_{IK} \longrightarrow H^*(Y(1)).$$
(i.e., first contracting the first two coordinates, or first contracting the second two coordinates)
both coincide; similarly  the other way around. %
We make the rest of the argument in a more abstract setting.

Suppose $R_{1}, R_2$ are two rings, 
and we are given an $(R_1, R_2)$-bimodule $M_{12}$ and an $(R_2,R_1)$-bimodule $M_{21}$,
giving associated functors $$F(-) =  - \otimes_{R_1} M_{12}, \ \ G(-) = - \otimes_{R_2} M_{21}$$ 
from (right) $R_1$-modules to (right) $R_2$-modules and vice versa, respectively.

We assume that $F,G$ define an equivalence of categories, i.e. there are natural equivalences
from $FG$ to the identity functor and from $GF$ to the identity functor. 
(In our setting above $R_1=\mathrm{H}_K', R_2 = \mathrm{H}_I'$ -- the primes denote localization at the ideal of $\mathrm{H}_K$ induced 
by $\mathfrak{m}$ --   
and the bimodules are $\mathrm{H}_{KI}',  \mathrm{H}_{IK}'$; the assumptions are satisfied by  Lemma \ref{IHequivalence}.)
 
  Next let $X$ be an $R_1$-module, let $Y$ be an $R_2$-module. (In our setting these are given by the  localized homology of $Y(1)$ and $Y_0(q)$, respectively.) 
Suppose also we're given 
maps $\alpha: F(X) \rightarrow Y, \ \beta: G(Y) \rightarrow X$ in such a way that the induced maps
\begin{equation} \label{oq}  G(F(X)) \stackrel{G(\alpha)}{ \longrightarrow }G(Y)  \stackrel{\beta}{\longrightarrow}  X \end{equation} 
\begin{equation} \label{oq2}   F(G(Y))  \stackrel{F(\beta)}{\longrightarrow}  F(X)  \stackrel{\alpha}{\longrightarrow }  Y \end{equation} 
arise from the specified maps $GF \rightarrow \mathrm{id}$ and $FG \rightarrow \mathrm{id}$.

Then the maps $F(X) \rightarrow Y$ and $G(Y) \rightarrow X$ must be  surjections
(by inspection of \eqref{oq} and  \eqref{oq2}); then in the diagram
$ G(F(X)) \rightarrow G(Y) \rightarrow X$
we have a composite of surjections giving an isomorphism,  so both are isomorphisms; in particular,
$G(Y) \rightarrow X$ and $F(X) \rightarrow Y$ must be isomorphisms.    

Applied to our original context, this concludes the proof of the first assertion of the Lemma.
\qed 

\medskip

 Note also that inside $\mathrm{H}_I$ we have a copy of the monoid algebra  
 $k[X_*^+]$, namely, the action of ``$U_{q}$-operators'' $I_q \chi I_q$
 for $\chi \in X_*^+$ (see \eqref{IWidentification} for identification of $X_*$ with a coset space; $X_*^+$
 is the positive cone corresponding to the Borel subgroup $\mathbf{B}$). 
 Each element $t \in T^{\vee}(k)$ defines a character  $\chi_t: k[X_*^+] \rightarrow k$ of this monoid algebra: 
 $t$ corresponds 
  (by ``Local Langlands'', easy in this split case)
 to an unramified character $\mathbf{A}(\Q_{q}) \rightarrow k^{\times}$, i.e. to a homomorphism $X_* \rightarrow k^{\times}$
 and then we just take the linear extension
\begin{equation} \label{chitdef} \chi_{t}: k[X_*^+] \rightarrow k.\end{equation}  \index{$\chi_t$} 
 Using the  previous Lemma, it is easy to compute the action of this monoid algebra:
  \begin{corollary} \label{MoritaCor}   
  Suppose that $q$ is part of a Taylor-Wiles datum. Then  the generalized  eigenvalues of $k[X_*^+]$ acting on $H^*(Y_0(q), k)_{\mathfrak{m}}$
are all of the  form $\chi_{\Frob^T}$, \index{$\chi_{\Frob^T}$}
where    $\Frob^T \in T^{\vee}(k)$ is conjugate to the Frobenius at $q$, and the notation is defined in \eqref{chitdef}. 
\end{corollary} 
The corollary asserts just the usual 
 relationship between ``$T_{q}$ eigenvalues at level $1$ and $U_{q}$ eigenvalues at level $q$.''  We will omit the proof; it follows
 from the prior Lemma and a straightforward computaiton.

\begin{remark}  Let $H^*(Y_0(q), \Z/p^n)_{\mathfrak{m}, \chi_{\Frob^T}}$ be the summand of $H^*(Y_0(q), \Z/p^n)_{\mathfrak{m}}$
corresponding to the $\chi_{\Frob^T}$-eigenspace of $k[X_*^+]$. 
For later use, we note that there is an isomorphism:
\begin{equation}\label{doubleprojremark} H^*(Y_0(q), \Z/p^n)_{\mathfrak{m}, \chi_{\Frob^T}} \cong H^*(Y(1), \Z/p^n)_{\mathfrak{m}}\end{equation}
where  the map from left to right is the  pushforward $\pi_*$, and an inverse in the other direction is given   by the pullback  $\pi^*$ together with projection to the
$\chi_{\Frob^T}$ eigenspace.  

 To see these are inverses we compute in $\mathrm{H}_{I}$: with reference to the actions of \eqref{juin}, the forward (pushforward) map corresponds  
to $e_K \in \mathrm{H}_{IK}$, and the reverse map corresponds to
 $|W| e_K q  \in \mathrm{H}_{KI}$ where $q \in S[X_*^+] \hookrightarrow \mathrm{H}_I$
is chosen to realize to the projection on the $\chi_{\Frob^T}$ eigenspace,
and $e_K q$ means the product of $e_K \in \mathrm{H}_{KI}$ with $q \in \mathrm{H}_I$.
 (See remarks after \eqref{eK definition}). 

 When we compose them we get
$|W| e_K q \in \mathrm{H}_{I}$ or $|W| e_K q e_K \in \mathrm{H}_K$; to see, for example, that    the former acts as the identity endomorphism
on $H^*(Y_0(q), \Z/p^n)_{\mathfrak{m}, \chi_{\Frob^T}}$, 
observe that it  can be written as $ \sum_{w \in W} e_w q$ with $e_w = I w I$,
 and then we just use the fact that $q$ annihilates all the  $k[X_*^+]$-eigenspaces  except the one indexed by $\Frob_T$,
 whereas the $e_w$ permutes the various eigenspaces. 
\end{remark} 
   
   \subsection{Adding level $q$ structure continued: The homology at level $Y_1(q, n)$} \label{HYln} 
   
   This section together with the previous and subsequent ones -- \S \ref{Morita}, \S \ref{HYln},   \S \ref{Lgcompat} --  all gather some ``standard'' properties of passing between level $1$ and level $q$. So  continue with $q$ as in the prior \S, i.e., part of a Taylor--Wiles datum of level $n$. 
   
   We now consider the  homology of $Y_1(q, n)$ -- as defined in \eqref{Y1def} -- and the action of its  Hecke algebra at level $Y_1(q, n)$. 
   Now this homology is basically glued from the homology of $Y_0(q)$:
   if we let  
   $\mathcal{F}$ be the push-forward of  the constant sheaf $k$
 from $Y_1(q,n)$ to $Y_0(q)$, then $\mathcal{F}$ is  a successive extension of copies of the constant sheaf  $k$,
 in a Hecke equivariant way. 
   
    We state  some consequences of this more formally:
   
   \begin{lemma}  \label{vanishing lemma BW} Assume that $q$ is part of a Taylor--Wiles datum of level $n$. 
   Assume local-global compatibility for $Y_0(q)$, in the sense described after \eqref{juin};
   let other assumptions be as in  \S \ref{Sec:assumptions}.
   Let  $\BWq,\delta$ be as in   \eqref{qdef}.  
   Then the homology $H_j(Y_1(q,n), \Z_p)_{\mathfrak{m}}$
  vanishes for $j \notin [\BWq, \BWq+\delta]$.   \end{lemma} 
 
\proof
Clearly we can replace the role of $\Z_p$ by $\mathbf{F}_p = k$, and then 
by the remark before the proof, it is enough to prove the same for  $Y_0(q)$.  
By \eqref{Twloc}, it suffices to prove the same vanishing statement for $Y(1)$
with $\mathbf{F}_p = k$ coefficients.  But this is part of our assumption (7(a) from 
 \S \ref{Sec:assumptions}.) \qed

   Now  we want to say that this relationship between the homology of $Y_1(q, n)$ and $Y_0(q)$ is equivariant
   for Hecke operators at $q$.  The full Hecke algebra is somewhat complicated and we just deal with its ``positive, commutative subalgebra.''  \index{$\widetilde{\mathcal{X}_*}$}
 Set $\widetilde{\mathcal{X}_*}$ to be the quotient of $\mathbf{A}(\Q_{q})$
   by the subgroup  $p^n \mathbf{A}(\Z_{q})$.       Denote by $\Delta_{q}$ the quotient of $\Q_{q}^{\times}$ by the subgroup of $\Z_{q}^{\times}$ of index $p^n$. Therefore $\Delta_{q}$ and  $\widetilde{\mathcal{X}_*}$ depend on $n$ but we will suppress that from the notation for simplicity.  We may identify
      $$\widetilde{\mathcal{X}_*} = X_*(\mathbf{A}) \otimes \Delta_{q},$$
      and we think of $\widetilde{\mathcal{X}_*}$ as a thickened version of the character lattice $X_*(\mathbf{A})$.  
 From the valuation $\Delta_{q} \rightarrow \Z$ we get
\begin{equation} \label{reduction} \widetilde{\mathcal{X}_*} \rightarrow  X_* \end{equation}
and we can define the ``positive cone'' $\widetilde{\mathcal{X}_*} ^+  \subset \widetilde{\mathcal{X}_*}$ as the preimage of $X_*^+$.

Let   $I$ be an Iwahori subgroup of $\G(\Q_q)$, and $I' \vartriangleleft I$
 the subgroup corresponding to the covering $Y_1(q, n)$, i.e. $I/I' \simeq (\Z/p^n)^r$. 
Then there is a natural map from $\Z_p[\widetilde{\mathcal{X}}_*^+]$  
to the Iwahori-Hecke algebra  at level $I'$ 
sending $\chi \in \widetilde{\mathcal{X}_*}$ to the coset $I'  \chi   I'$. 
In particular, we get an action of $\Z_p[\widetilde{\mathcal{X}}_*^+]$
on the homology of $Y_1(q, n)$.

\begin{lemma}
Any generalized eigenvalue of $k[\widetilde{\mathcal{X}}_*^+]$
acting on $H_*(Y_1(q, n),k)_{\mathfrak{m}} $ is also
a generalized eigenvalue of $k[\widetilde{\mathcal{X}}_*^+] $ acting on
$H_*(Y_0(q),k)_{\mathfrak{m}}$ via the map $k[\widetilde{\mathcal{X}}_*^+]  \rightarrow k[X_*^+]$ induced by \eqref{reduction}. 
 \end{lemma}

 Thus, by  Corollary
\ref{MoritaCor}, we get a splitting  
 \begin{equation} \label{bigsplit} H_*(Y_1(q, n),k)_{\mathfrak{m}} = \bigoplus_{\stackrel{ \Frob^T \in T^{\vee}(k)}{\Frob^T \sim \rhobar(\Frob_q)}} \left( H_* (Y_1(q,n), k) \right)_{\mathfrak{m},\chi_{\Frob^T}}\end{equation} 
   into   the  sum of generalized eigenspaces associated to the characters $\chi_{\Frob^T}: k[\widetilde{\mathcal{X}}_*^+] \rightarrow k$. (In the subscript, $\sim$ means ``is conjugate to.'') 
 Again, this is nothing but a fancy way of talking about the decomposition into ``$U_{q}$-eigenspaces.''  The only point to note is that the decomposition is
canonically  indexed by elements of $T^{\vee}(k)$ conjugate to the Frobenius. 

 \proof   
  \iflongproof Write for short $G = \G(\Q_{q})$; let $I$ be an Iwahori subgroup of $G$, and $I' \vartriangleleft I$
 the subgroup corresponding to the covering $Y_1(q, n)$, i.e. $I/I' \simeq (\Z/p^n)^r$. 
 We prove the same statement in cohomology and without the $\mathfrak{m}$; the desired
 statement follows by dualizing and localizing.
  
 By the discussion of \S \ref{Arithmeticmanifolds} we can identify
 $$H^*(Y_1(q, n), k) \simeq \Ext^*_{kG} (k[G/I'], M)$$
 $$ H^*(Y_0(q), k) \simeq \Ext^*_{kG}(k[G/I], M)$$
 
 where $M$ is the direct limit of co-chain complexes of  a family of coverings, obtained
 by adding more and more level structure at $q$.  
 
 {\em Claim:} We may filter $k[G/I']$ by $G$-submodules 
 $F^0 \subset F^1 \subset F^2 \subset \dots $
such that: 
 \begin{itemize}
 \item[(i)] each successive quotient $F^{i+1}/F^i$ is isomorphic as $G$-module to a sum of copies of $k[G/I]$,
 \item[(ii)] For every $\chi \in \widetilde{\mathcal{X}_*^+}$,  the action of $I' \chi I'$ preserves the filtration, and  the action on the quotients coincides with the action of $I \chi I$. 
 \end{itemize}
 
Assuming the existence of this filtration, the result follows easily: we get  long exact sequences  of the form $$ \Ext^j_{kG}(F^{i-1}, M) \rightarrow \Ext^j_{kG}(F^i, M) \rightarrow \Ext^j( k[G/I], M )^{\bigoplus e}  \rightarrow $$
 and these are equivariant for the action of $I' \chi I'$, which acts by via $I \chi I$ on the right-hand summand.
 The result follows immediately by a descending induction.

To construct the desired filtration regard $k[G/I']$ as the  compact induction from $I$ to $G$ of the $I$-representation $k[I/I']$, i.e
 $$k[G]\otimes_{k[I]} k[I/I'] $$
 In particular, $\Delta=I/I'$ acts  by $G$-endomorphisms on $k[G/I']$, which is the action ``by right multiplication.'' This action 
 of $\delta \in \Delta$ coincides with the action of $I' \delta I'$ and in particular 
 commutes with the action of $k[\widetilde{\mathcal{X}}_*^+]$. 
 
 Now we filter $k[G/I']$ by the kernels  $k[G/I'] \langle \mathfrak{m}^j \rangle$ of successive powers  $\mathfrak{m}^j$ of the maximal 
ideal  $\mathfrak{m}$  in $k[\Delta]$.    This filtration  is stable
for  $k[\widetilde{\mathcal{X}}_*^+]$ because the actions of $\Delta, k[\widetilde{\mathcal{X}}_*^+]$ commute. 
Also 
the $j$th term of the resulting filtration is thus
$$F^j= k[G] \otimes_{k[I]}  k[\Delta]\langle \mathfrak{m}^j \rangle$$
and the $j$th graded is just $k[G] \otimes_{k[I]} \frac{\langle  \mathfrak{m}^j \rangle}{\langle \mathfrak{m}^{j+1} \rangle}$, 
i.e. a direct sum of copies of $k[G/I]$, as claimed.

 It remains to check assertion (ii) in the {\em Claim.} For any $y \in \mathfrak{m}^{j-1}$, multiplication by $y$ gives a map
$$F^{j}/F^{j-1} \rightarrow F^1 = F^1/F^0$$
and a suitable sum of such maps is an isomorphism (as we see by checking the corresponding assertion for $k[\Delta]$). 
Since these multiplication maps commute with the action of $k[\widetilde{\mathcal{X}}_*^+]$
we are reduced to computing the action of $I' \chi I'$ on $F^1$;
now $F^1$ is identified with $k[G/I]$ is a natural way and the assertion is clear.
 \fi 
  \qed

      \subsection{Adding level $q$ structure continued:  Galois representations for level  $Y_1(q ,n)$} \label{Lgcompat}
      
      This section and the prior two \S \ref{Morita}, \S \ref{HYln},   \S \ref{Lgcompat} all gather some ``standard'' properties of passing between level $1$ and level $q$.
      We now consider more closely the action of the Iwahori-Hecke algebras at level  $Y_1(q ,n)$
      and formulate local-global compatibility. First let us look at the Galois side.

\begin{lemma} \label{TWdeflem}  Let $q$ be a Taylor-Wiles prime of level $n$
and assume that the unramified representation $\rhobar|G_{\Q_q}$ has image inside $T^{\vee}$. Then any deformation
  of  $\rhobar|_{G_{\Q_{q}}}$ can be conjugated to one
taking values in $T^{\vee}$.    In particular, any such deformation of $\rhobar|_{G_{\Q_{q}}}$ actually factors through $G_{\Q_{q}}^{\ab} \simeq \mbox{profinite completion of }\Q_{q}^{\times}$.  
\end{lemma}

   \proof 
   \iflongproof
     We may present the tame quotient of $G_{\Q_{q}}$ as $\langle F, t: F t F^{-1} = t^{q} \rangle$,
  where $t$ is a generator of tame inertia. 
   
   Suppose that $A$ is an Artin local ring with maximal ideal $\mathfrak{m}$, and $\mathfrak{m}^d=0$.   We are given $F, t$
   in $G^{\vee}(A)$ 
   that satisfy $ F t F^{-1} = t^{q}$,
 where $t$ reduces to the identity in $G^{\vee}(k)$ and $F$ reduces to  (after conjugating) a strongly regular element of $T^{\vee}(k)$. 
 Conjugating, we may suppose
   that $F$ belongs to the maximal torus $T^{\vee}(A)$.  We will prove by induction on $d$ 
   that this  forces $t \in T^{\vee}(A)$ too. 
   By the inductive hypothesis (obvious for $d=1$) 
  the image of $t$ in $G^{\vee}(A/\mathfrak{m}^{d-1})$
   belongs to $T^{\vee}(A/\mathfrak{m}^{d-1})$. 
   Write thus $t =t_0   \delta_t$ where $t_0 \in T^{\vee}(A)$ and $\delta_t \in G^{\vee}(A)$ 
   lies in the kernel of reduction modulo $\mathfrak{m}^{d-1}$.
   
     Now
in fact $t_0 $ lies in the kernel of reduction modulo $\mathfrak{m}$, and so $t_0$ and $\delta_t$ actually commute; indeed, $t_0$
commutes with anything in the kernel of reduction modulo $\mathfrak{m}^{d-1}$.  Also   $\delta_t^{q-1} = e$.   (To check these statements,
just compute in the formal group of $G^{\vee}$ at the identity.) 
 Now,    
$F t F^{-1} = t^{q}$ so that $$ F (t_0 \delta_t)  F^{-1} (t_0 \delta_t)^{-1} =  (t_0 \delta_t)^{q-1}  \in T^{\vee}(A)$$
But the left hand side equals 
 $t_0  \left( \Ad(F)\delta_t  \cdot \delta_t^{-1}\right) t_0^{-1} $, and so 
 $$\Ad(F) \delta_t \cdot \delta_t^{-1}  \in T^{\vee}(A) $$
and  since $F$ is strongly regular this means that $\delta_t \in T^{\vee}(A)$ as desired. \fi
   \qed
   
   \medskip
   
   We now want to connect the Galois deformation ring relevant to $Y_1(q,n)$ with the Iwahori--Hecke algebra.

 Let $q$ be a Taylor--Wiles prime of level $n$. 
   Suppose fixed an element $\Frobl \in T^{\vee}(k)$ conjugate to the Frobenius at $q$. 
Consider a deformation   $\sigma: G_{\Q} \rightarrow G^{\vee}(R)$
    of $\rhobar$, where one allows now ramification at $q$,
    and $R$ is an Artin local ring with residue field $k$. 
     Then we can  uniquely conjugate
   $\sigma$ so its restriction to     to $G_{\Q_{q}}$ factors as 
$$ G_{\Q_{q}}^{\ab} \simeq  \widehat{\Q_{q}^{\times}} \longrightarrow T^{\vee}(R),$$
and the image of a uniformizer in $\Q_q^{\times}$ reduces to $\Frobl$.  This map  factors through (the profinite completion of) 
$\Q_q^{\times}/(1+q \Z_q)$. Restricting to $\mathbf{F}_{q}^{\times}$ we get
$$\mathbf{F}_{q}^{\times} \rightarrow T^{\vee}(R)$$
and pairing with characters of $T^{\vee}$ we get
$ \mathbf{F}_{q}^{\times} \times X^*(T^{\vee}) \rightarrow R^*$;
by the duality of $T^{\vee}$ and $\mathbf{A}$, this is the same thing as
\begin{equation} \label{rgh}  \mathbf{A}(\mathbf{F}_{q}) \rightarrow R^*. \end{equation}
We emphasize that the map \eqref{rgh} depended on the choice of a toral element $\Frobl$ conjugate to Frobenius; changing this element
changes the map through the action of the Weyl group. 

 Let $R^{\univ}_{\ramprimes \cup \{q\}}$ be Mazur's universal deformation ring for $\rhobar$, allowing ramification at $q$.  
By our assumptions  (\S \ref{GaloisAss}) on the existence of Galois representations, there is a map $$R^{\univ}_{\ramprimes \cup \{q\}} \rightarrow \mathbb{T}_{K_1(q,n), \mathfrak{m}},$$
where $K_1(q,n)$ is the level structure for $Y_1(q,n)$.
By means of this map, $R^{\univ}_{\ramprimes \cup \{q\}}$ acts on $H^*(Y_1(q,n), \Z_p)_{\mathfrak{m}}$,
and in particular on the summand $H^*(Y_1(q,n), \Z_p)_{\mathfrak{m}, \Frob_q^T}$ under \eqref{bigsplit}. 
Thus, by \eqref{rgh}, we get an action of $\mathbf{A}(\mathbf{F}_{q})$ this cohomology group. 

Now the assumption of local--global compatibility alluded to in \S \ref{GaloisAss}
is a strengthened version of the following:
\begin{quote}
{\em Local global compatibility: }
The action  $\mathbf{A}(\mathbf{F}_q) \acts H^*(Y_1(q,n), \Z_p)_{\mathfrak{m}, \Frob_q^T}$ 
just defined coincides with the ``geometric'' action, i.e.
wherein $\mathbf{A}(\mathbf{F}_q)$ acts by deck transformations   on $Y_1(q,n)$ (see \eqref{Y1def}).
\end{quote}
By ``strengthened version'', we mean that we require a similar assertion at a derived category level, not just at the level of cohomology,
and we also require the assertion for several auxiliary primes $q$ rather than a single one. 
For details, see \cite[\S 13.5]{GV}.

 We say a deformation  of $\rhobar|G_{\Q_{q}}$ is of ``inertial level $\leq n$'' \index{inertial level}  \label{inertial level}
 if,  when considered as a representation of $\Q_{q}^{\times}$ by  Lemma \ref{TWdeflem},  and restricted to $\mathbf{F}_{q}^{\times}$,  
 it factors through the quotient $\mathbf{F}_{q}^{\times}/p^n$. 
 We denote by $R^{\univ, \leq n}_{\ramprimes \cup \{q\}}$ the quotient of $R^{\univ}_{\ramprimes \cup \{q\}}$ that classifies deformations of $\rhobar$
 such that $\rhobar|G_{\Q_{q}}$ has inertial  level $\leq n$.  
 Explicitly, $R^{\univ,\leq n}_{\ramprimes \cup \{q\}}$ is the quotient of 
\begin{equation} \label{Rndef} R^{\univ}_{\ramprimes \cup \{q\}} / \langle t-1: \mbox{ $t$ is in the image of  $p^n \mathbf{A}(\mathbf{F}_{q})$ under \eqref{rgh}} \rangle.\end{equation}
  Then local-global compatibility implies that the action of $R^{\univ}_{\ramprimes \cup \{q\}}$ on the homology of $Y_1(q, n)$ factors through $R^{\univ, \leq n}_{\ramprimes \cup \{q\}}$, according to our previous discussion.

\section{Patching and the derived Hecke algebra} \label{Patching2}

We continue with the notation  and assumptions of the previous section \S \ref{Patching}; see in particular
\S  \ref{Sec:assumptions} and \S \ref{GaloisAss}. 
That section was primarily setup, and now we get down to proving that the global derived Hecke algebra is ``big enough,''
in the sense discussed around 
\eqref{bigenough}. 
The main result is
Theorem \ref{mindegreeproof}. 

We use the patching of the Taylor--Wiles method; more specifically, we use the 
version of that method  that was discovered \cite{CG} by Calegari and Geraghty,
which applies to situations where the same Hecke eigensystem occurs in multiple degrees.
We also use heavily the presentation of the Calegari--Geraghty method 
given by Khare and Thorne \cite{KT}. 
\subsection{Deformation rings and chain complexes at level $Q_n$.}  \label{Qnchain}

 Fix now a Taylor-Wiles datum $Q_n$ of level $n$. 
  (We will   abusively use $Q_n$ both to denote the Taylor-Wiles datum and simply the set of primes associated to that datum.)
  
  Recall the definition of $Y_1^*(Q_n)$ from \S \ref{TWprimes stuff}: it is the fiber-product of coverings $Y_1(q, n) \rightarrow Y(1)$
  over $q \in Q_n$. 
 We now collect together various results about the homology of $Y_1^*(Q_n)$, which are essentially the same
 results as those already discussed for $Y_1(q,n)$, but using all the primes in $Q_n$ instead of just $\{q\}$.

 \begin{lemma} The homology $H_j(Y_1^*(Q_n), \Z_p)_{\mathfrak{m}}$
  vanishes for $j \notin [\BWq, \BWq+\delta]$. 
 \end{lemma} 
 \proof As in Lemma \ref{vanishing lemma BW}. \qed
 
 Just as in \eqref{bigsplit}, this homology group is split (by ``$U$-operators'') into summands
 indexed by collections $\Frob^T_q \in T^{\vee}(k) \ \ (q \in Q_n)$, where each $\Frob^T_q$
 is conjugate to the Frobenius at $q$. 
  In particular, since the Taylor--Wiles datum is equipped (\S \ref{TWprimes stuff}) with a specific choice of such a $\Frob_q^T$ for each $q \in Q_n$, we
  can consider  the summand 
\begin{equation} \label{Weylsplitting}  H_*(Y_1^*(Q_n), \Z_p)_{\mathfrak{m},\Frob^T_{Q_n}}  \subset H_*(Y_1^*(Q_n), \Z_p)_{\mathfrak{m}}  \end{equation}
indexed by  these prescribed lifts. 

  Recall  from \eqref{idontcare}  that $Y_1^*(Q_n) \rightarrow Y_0(Q_n)$ is Galois, with Galois group $T_n$.  \index{$R_n$} \index{$R_{Q_n}$} 
Let us introduce notation for the deformation rings of interest to us: let 
 $$R_{Q_n} = \mbox{universal deformation ring at level $\ramprimes \coprod Q_n$,}$$
 \begin{equation} \label{RnRn}  R_n = \mbox{ quotient of $R_{Q_n}$ classifying   deformations of inertial level $\leq n$ at primes in $Q_n$}
 \end{equation}
  For example,  in the case when $Q_n = \{q\}$ this was the ring $R^{\univ,\leq n}_{\ramprimes \cup \{q\}}$ discussed around \eqref{Rndef}.  
 
  By the discussion  of  \eqref{rgh}  and after, we get a morphism
   $\mathbf{A}(\mathbf{F}_{q_i})/p^n \rightarrow R_n^{\times}$, and therefore we get (see \eqref{idontcare}): 
 $$T_n \rightarrow R_n^{\times}$$
 What we know (local-global compatibility, assumed in \S \ref{Lgcompat})  is that the natural action of $T_n$  (deck transformations) on homology of $Y_1^*(Q_n)$ 
is compatible with its action via $T_n \rightarrow R_n^{\times}$.  
 To say differently, we get a map \begin{equation} \label{Sndef} S_n := \Z/p^n[T_n] \rightarrow R_n/p^n, \end{equation}  and the natural action of $S_n$ on the $\Z/p^n$ homology of $Y_1^*(Q_n)$ is compatible
 with that via the map to $R_n/p^n$.

  Now consider the complex of singular chains   $$ \widetilde{C}_n = \Chains(Y_1^*(Q_n); \Z/p^n)$$
  with $\Z/p^n$ coefficients.  We think of it as a complex of $S_n$ modules, because of the action of $T_n$ by deck transformations on $Y_1^*(Q_n)$.   It is quasi-isomorphic to a bounded complex of finite free $S_n$-modules
  and we have canonical identifications:
  
 \begin{equation} \label{id1} H_* \widetilde{C}_n \simeq H_*(Y_1^*(Q_n), \Z/p^n)\end{equation}
\begin{equation} \label{id2} H^*\Hom_{S_n}(\widetilde{C}_n, \Z/p^n) \simeq H^*(Y_0(Q_n), \Z/p^n)\end{equation}
\begin{equation} \label{id0} H_*( \widetilde{C_n} \otimes_{S_n} \Z/p^n) \simeq  H_*(Y_0(Q_n), \Z/p^n).\end{equation}

  Note that $\widetilde{C}_n$ is a free $S_n$-module, with basis given by the characteristic functions of an arbitrarily chosen set of representatives
  for $T_n$-orbits on singular simplices. 
  Therefore, the homology  of $\Hom_{S_n}(\widetilde{C}_n, \Z/p^n)$ computes the homomorphisms from $\widetilde{C}_n$ to $\Z/p^n$
  in the derived category of $S_n$-modules. 
  By composition of homomorphisms in this derived category, we get a map
 \begin{equation} \label{above} \underbrace{H^*  \left( \Hom_{S_n}(\widetilde{C}_n, \Z/p^n) \right) }_{\simeq H^*(Y_0(Q_n),\Z/p^n)} \times  \underbrace{ \Ext^*_{S_n}(\Z/p^n, \Z/p^n) }_{\simeq H^*(T_n,\Z/p^n)} \rightarrow  \underbrace{ H^* \left( \Hom_{S_n}(\widetilde{C}_n, \Z/p^n)\right) }_{\simeq H^*(Y_0(Q_n),\Z/p^n)} \end{equation}
   With respect to the identifications noted underneath the respective terms, this is precisely the ``natural'' action of $H^*(T_n, \Z/p^n)$ on $H^*(Y_0(Q_n), \Z/p^n)$.
  This natural action arises thus: 
   the covering $Y_1^*(Q_n) \rightarrow Y_0(Q_n)$ has covering group $T_n$, i.e. can be regarded as a map
\begin{equation} \label{BTmap} Y_0(Q_n) \rightarrow \mathrm{B}T_n\end{equation} 
   from $Y_0(Q_n)$ to the classifying space of $T_n$; this allows one to pull back cohomology classes from $T_n$ and take cup product. 
      The coincidence of \eqref{above} and this ``natural action'' is a general fact; for lack of a reference we sketch a proof in \S \ref{sec:remedialtopology}.

It is possible to ``cut down'' $\widetilde{C}_n$ in a fashion that corresponds to  the summand
  \eqref{Weylsplitting}, as is explained in \cite{KT} (see Lemma 2.12 thereof,  and surrounding discussion).  
  This can be done compatibly for $Y_0(Q_n)$ and $Y_1^*(Q_n)$ and thus one gets 
  a perfect complex $C_n$ of $S_n$-modules, equipped with identifications 
that are analogous to \eqref{id1} and \eqref{id2}: 
 \begin{equation} \label{id3} H_*(C_n; \Z/p^n) \simeq H_*(Y_1^*(Q_n), \Z/p^n)_{\mathfrak{m}, \Frob^T}, \end{equation} 
\begin{equation} \label{id4} H^*(\Hom_{S_n}(C_n, \Z/p^n)) \simeq H^*(Y_0(Q_n), \Z/p^n)_{\mathfrak{m}, \Frob^T}\end{equation}
  and again the action of $\Ext^*_{S_n}(\Z/p^n , \Z/p^n)$ on the latter group corresponds to the natural action    by pulling back cohomology classes via \eqref{BTmap}.\footnote{As a sanity check on this, note that the action of $H^*(T_n, \Z/p^n)$ on $H^*(Y_0(Q_n), \Z/p^n)$ indeed
  does preserve the splitting into  summands of the type \eqref{Weylsplitting}; one can see this directly  
  by seeing that $H^*(T_n, \Z/p^n)$, considered inside the  derived Iwahori-Hecke algebra, commutes   with the ``positive subalgebra'' used to define the splitting \eqref{Weylsplitting}.}

\subsection{Extracting the limit}  \label{extraction}    Now we ``pass to the limit'' as per Taylor--Wiles and  Calegari--Geraghty.
 The idea is roughly speaking to extract, by a compactness argument, a subsequence of $n$ along which
the $C_n, S_n, R_n$ are compatible, and then get limits $\Cinf, \Sinf, \Rinf$ by an inverse limit. Usually in modularity lifting
one is only concerned with the limit of the process; but in our case we also want to remember
some facts about how this relates to the $C_n, S_n, R_n$. 
A discussion of this process which emphasizes exactly what we need is given in \cite[\S 13]{GV}, see in particular Theorem 13.1 therein.

 We choose a sequence of Taylor-Wiles data $Q_n$ with  $n \rightarrow \infty$.  
 After replacing the $Q_n$ by a suitable subsequence and then reindexing  -- that is to say, 
 replacing $Q_i$ by $Q_{n_i}$ for some $n_i > i$, and then regarding $Q_{n_i}$ as a set of Taylor--Wiles primes of level  $i$ -- 
 we can arrange that we  can ``pass to the limit.''
After having done this, we obtain at last the following data:   
  
\begin{itemize}
 
\item[(a)] A sequence of Taylor--Wiles data $Q_n$ of level $n$.

 Recall to this we have associated coverings $Y_1^*(Q_n) \rightarrow  Y_0(Q_n) \rightarrow Y(1)$,
as in \S \ref{TWprimes stuff}, and the Galois group of the former map is called $T_n$; also \eqref{Sndef} we set $S_n = \Z/p^n[T_n]$, the group algebra of $T_n$. 

\item[(b)] With $\Sinf = \Z_p[[x_1, \dots, x_R]]$ as in
\S \ref{sec:diamond}, 
a complex $\Cinf$ of finite free $\Sinf$-modules, 
equipped with 
a quasi-isomorphism \begin{equation} \label{Descent} \Cinf \otimes_{\Sinf} S_n  \simeq C_n,\end{equation}   %
where $C_n$ is as described in \S \ref{Qnchain}: a version of the chain complex of $Y_1^*(Q_n)$ with $\Z/p^n$ coefficients,
but localized at $\mathfrak{m}$ and $\Frob^T$. 
\medskip
\item[(c)]  
 A quotient $\usbR_n$ of $R_n$,   \index{$\usbR_n$} 
defined as follows:

Recall   from \eqref{RnRn} the definition of $R_n$; a quotient of the crystalline deformation ring
at level  $\ramprimes \coprod Q_n$.  We set \begin{equation} \label{goef} \usbR_n = R_n/(p^n, \mathfrak{m}^{K(n)})\end{equation}
for a certain explicit function $K(n)$, chosen so that  e.g.  action of $R_n$ on $H_*(C_n)$  
automatically factors through $\usbR_n$.  We can and will assume $K(n) \geq 2n$.   (The main function of $K(n)$ is to make $\usbR_n$ Artinian, 
while still retaining enough information about all of $R_n$ for our purposes.) 

\item[(d)]

 A ``limit deformation ring'' $\Rinf \simeq \Z_p[[x_1, \dots, x_{R-\delta}]]$  equipped with 
maps  $\Sinf \rightarrow \Rinf$ and maps $\Rinf \rightarrow \usbR_n, \Rinf \twoheadrightarrow \defring$   which are compatible,
in the sense that this diagram commutes: 
    \begin{equation} \label{SinftySn}
 \xymatrix{
\Sinf  \ar[r] \ar[d] &   \Rinf  \ar[d] \ar[r]  & \defring   \ar[d] \\
S_n   \ar[r]   &   \usbR_n   \ar[r] & \defring/  (p^n, \mathfrak{m}^{K(n)}) }
\end{equation}

(Recall here that $\defring$ is the deformation ring of $\rhobar$, with crystalline conditions imposed, without adding any level, 
cf. \S \ref{GaloisAss}). 

Moreover, the composite $\Sinf \rightarrow \Rinf \rightarrow \defring$ factors through 
the augmentation $\Sinf \rightarrow \Z_p$; and also the  left-hand square induces an isomorphism
 \begin{equation} \label{second-iso}  \Rinf \otimes_{\Sinf} S_n \simeq  \usbR_n \end{equation}

\item[(e)] An  action of $\Rinf$ on $H_*(\Cinf)$, compatible with the $\Sinf$ action,
and with the maps $H_*(\Cinf) \rightarrow H_*(C_n)$, where
$\Rinf$ acts on $H_*(C_n)$ via $\Rinf \rightarrow  \usbR_n$.

\item[(f)]   An identification of  \begin{equation} \label{logstructure} H^*(\Hom_{\Sinf}(\Cinf, \Z_p)) \stackrel{\sim}{\longrightarrow} H^*(Y(1), \Z_p)_{\mathfrak{m}},\end{equation}   
compatible under \eqref{Descent} with 
the identification $H^*(\Hom_{S_n}(C_n, \Z/p^n)) \simeq H^*(Y(1), \Z/p^n)_{\mathfrak{m}}$
that is the composition of \eqref{id4} with the pushforward. (This is described in dual form in \cite[Theorem 13.1(d)]{GV} but one gets
similarly this result, and the statement about compatibility is just a matter of looking at the definition of the map \eqref{logstructure}. 

\item[(g)]

(These last results use heavily the formal smoothness, assumption (e) from  \S \ref{GaloisAss}):
 $\Cinf$  has homology only in degree $q$, and its homology there $H_q(\Cinf)$ is free   as $\Rinf$-module. Moreover, one has 
an ``$R=T$ result''
\begin{equation} \label{first-iso} \Rinf \otimes_{\Sinf} \Z_p \simeq \defring \simeq \mbox{ image of $\mathbb{T}_{K_0} $ in  $\End \ H_q(Y(1), \Z_p)_{\mathfrak{m}}$}. \end{equation}
 \end{itemize}

 \subsection{The structure of $\Sinf \rightarrow \Rinf$} \label{SR} 
 
 The limit process has given a map of rings $\Sinf \rightarrow \Rinf$, 
 where $\Sinf$ and $\Rinf$ are formal power series rings that represent, roughly speaking,
 ``limits'' of the rings $S_n, R_n$ as $n \rightarrow \infty$. 
  
As in \eqref{Zplift}, the  representation $\Pi$ gives a lift $\rho$ to $\Z_p$ of the residual representation $\rhobar$;
 this corresponds to an augmentation $\defring \rightarrow \Z_p$.  
Thus we also get an augmentation $$f: \Rinf \rightarrow \Z_p,$$
and the pullback of this to $\Sinf$ is the natural augmentation of $\Sinf \rightarrow \Z_p$ ((d) of \S \ref{extraction}).
In particular, the kernel of $f$ on $\Sinf$ is precisely the ideal
$\Iinf$.

 Our assumptions imply that the map $\Sinf \rightarrow \Rinf$ is surjective.
 Indeed, because $\Sinf$ is complete for the $\Iinf$-adic topology
 it is enough to verify that $\Sinf/\Iinf \rightarrow \Rinf/\Iinf \Rinf$
 is surjective.  But
 $\Rinf/\Iinf$ is a Hecke ring by \eqref{first-iso} and so isomorphic to $\Z_p$
 by \eqref{nocong0}. 
  Note in particular that this also means that $\Iinf \Rinf$ is precisely the kernel of $f$.

 The following easy  lemma is now useful for explicit computations. 
 
 \begin{lemma} \label{straightening} We can choose generators  $x_i, y_j$ for $\Sinf, \Rinf$, i.e.
  $$\Sinf = \Z_p[[x_1, \dots, x_R]],  \ \ \Rinf = \Z_p[[y_1, \dots, y_{R-\delta}]]$$ 
 such that the $x_i, y_j$s lie in the kernel of the compatible augmentations 
 $$\Sinf \rightarrow \Rinf \rightarrow \Z_p,$$  and 
the map $\Sinf \rightarrow \Rinf$  is given by $x_i \mapsto y_i$ for $i \leq R-\delta$ and $x_i \mapsto 0 $ for $i > R-\delta$. 
 \end{lemma}

\proof

 Write $f: \Rinf \rightarrow \Z_p$ for the augmentation. 
  Abstractly, $\Rinf \simeq \Z_p[[u_1, \dots, u_{R-{\delta}}]]$ where all the $u_i$ lie in the maximal ideal.  
Set  $y_i = u_i - f(u_i) \in  \ker(\Rinf \rightarrow \Z_p)$. Then still $\Rinf \simeq \Z_p[[y_1, \dots, y_{R-\delta}]]$. 
We have noted above that $\Jinf := \Iinf \Rinf$ is precisely the kernel of the augmentation $\Rinf \rightarrow \Z_p$; thus, 
  the $y_i$ freely span as $\Z_p$-module the quotient $\Jinf/\Jinf^2$.  
 
 Lift the $y_i$ to $x_1, \dots, x_r \in \ker(\Sinf \rightarrow \Z_p)$. 
Necessarily the $x_i$  span a saturated\footnote{Here we say that a submodule $Q$ of a free $\Z_p$-module $Q'$ is saturated if the quotient
$Q'/Q$ is torsion-free.} $\Z_p$-submodule of rank $s$ inside $\Iinf/\Iinf^2 \simeq \Z_p^s$; 
they are $\Z_p$- independent because any linear relation $\sum a_i x_i \in \Iinf^2$ (with $a_i \in \Z_p$) 
would give rise to a corresponding linear relation in $\Rinf$, a contradiction. Similarly, they are saturated because 
given $x'$  and $(a_1 \dots, a_r)$ with $\gcd(a_1, \dots, a_r) =1$ and  $p x' = \sum a_i x_i + \Iinf^2$, we would get a corresponding relation in $\Rinf$, again a contradiction. 

Now  extend  the $x_i$ to a full $\Z_p$-basis
$x_{r+1}, \dots, x_s$ for $\Iinf/\Iinf^2$.  Each $x_j$ for $j > r$
is sent under $\Sinf \rightarrow \Rinf$ an element of $\Jinf \subset \Rinf$, which means that it can be written
as a formal polynomial  $P_j(y_1, \dots, y_r)$ in  $y_1, \dots, y_r$, with no constant term; so  replacing $x_j$ by $x_j - P_j(x_1, \dots, x_r)$ we may  suppose that $x_j \mapsto 0$ in $\Rinf$. 
\qed

Now we come to the main theorem of the section.  

\begin{theorem}  \label{mindegreeproof} 
Let assumptions be as in \S \ref{Sec:assumptions}
and \S \ref{GaloisAss}. 
The cohomology  $H^*(Y(1), \Z_p)_{\mathfrak{m}}$
  is generated, as a module over the strict global derived Hecke algebra (see \S \ref{limit} and \S \ref{Hecke:resplace} for definition with $\Z_p$ coefficients), by its minimal degree component $H^{\BWq}(Y(1), \Z_p)_{\mathfrak{m}}$.
 \end{theorem}

 \proof 
We use the setup of the Taylor--Wiles limit process (\S \ref{extraction}), beginning   with the fact that the natural map 
\begin{equation} \label{cyclicity2}H^{\BWq}(\Hom_{\Sinf}(\Cinf, \Z_p)) \otimes \Ext^j_{\Sinf}(\Z_p, \Z_p) \twoheadrightarrow   H^{\BWq+j}(\Hom_{\Sinf}(\Cinf, \Z_p))    \end{equation}
is a surjection for all $j$: by   (g) of \S \ref{extraction} and Lemma \ref{straightening}, 
we can choose coordinates so that $\Sinf \simeq \Z_p[[x_1, \dots, x_R]]$,
and the complex $\Cinf$ is quasi-isomorphic to {\color{\changecolor} a sum of copies of } $\Sinf/(x_{R}, \dots, x_{R-\delta+1})$
concentrated in a single degree.   So   the surjectivity of \eqref{cyclicity2} follows from the ``Koszul algebra'' computations in \S \ref{AppendixB} of the Appendix.

   Examine now the diagram, where all the maps are the obvious ones;
{\small     \begin{equation} \label{previousstory}
 \xymatrix{
 H^{\BWq}(\Hom_{\Sinf}(\Cinf, \Z_p))^{} \ar[d]^{=} & \times  &    \Ext^{i}_{\Sinf}(\Z_p, \Z_p)  \ar[d]  \ar[r] &  H^{\BWq+i}(\Hom_{\Sinf}(\Cinf, \Z_p)) \ar[d]^V \\
  H^{\BWq}(\Hom_{\Sinf}(\Cinf, \Z_p))   \ar[d]^{V}  & \times&  \Ext^{i}_{\Sinf}(\Z_p, \Z/p^n)    \ar[r]&   H^{\BWq+i}(\Hom_{\Sinf}(\Cinf,  \Z/p^n))\ar[d]^{=}\\
  H^{\BWq}(\Hom_{\Sinf}(\Cinf, \Z/p^n))^{}   & \times & \left( \Ext^{i}_{\Sinf}(\Z/p^n, \Z/p^n) \right) \ar[r]   \ar[u]^U  &  H^{\BWq+i}(\Hom_{\Sinf}(\Cinf, \Z/p^n))   \\
 H^{\BWq}(\Hom_{S_n}(\Cinf \otimes_{\Sinf} S_n, \Z/p^n))^{} \ar[u]^{f, \sim} & \times  &   \left( \Ext^{i}_{S_n}(\Z/p^n, \Z/p^n) \right) \ar[r]^{\mathfrak{Q} \ \  \ \ \ \ \ } \ar[u]^{U'} &  \ar[u]^{\sim} H^{\BWq+i}(\Hom_{S_n}(\Cinf \otimes_{\Sinf} S_n, \Z/p^n)). \\
  }
  \end{equation}
  }
Here, the middle square ``commutes'' in the sense that the image of $(x, Uy)$ is the same as the image
  of $(Vx, y)$, i.e. $U, V$ are adjoint for the pairing.  The top and bottom squares commute.   
  All this is obvious, except for perhaps the bottom square which involves change of rings, so let us talk through it: 
The map $\Sinf \rightarrow S_n$ induces a forgetful map $T$
from the derived category of
$S_n$-modules to the derived category of $\Sinf$-modules.  
Take $$\alpha \in H^{\BWq}(\Hom_{S_n}(\Cinf \otimes_{\Sinf} S_n, \Z/p^n)), \ \ \beta \in   \Ext^{i}_{S_n}(\Z/p^n, \Z/p^n).$$ 
We can regard $\alpha$ as a map $\Cinf \otimes_{\Sinf} S_n \rightarrow \Z/p^n[q]$
and $\beta$ as a map $\Z/p^n[q] \rightarrow \Z/p^n[q+i]$,  both in the derived category of $S_n$-modules.
Applying the functor $T$, we see that $T \beta \cdot T\alpha  = T(\beta \alpha)$. On the other hand,
$T \alpha$ is  simply the morphism
$\Cinf \otimes_{\Sinf} S_n \rightarrow \Z/p^n[q]$ considered as a map of $\Sinf$-modules.
If we pre-compose with $\gamma: \Cinf \rightarrow \Cinf \otimes_{\Sinf} S_n$, considered
as a map of $\Sinf$-modules, we get $f(\alpha)$. 
 Similarly, $  T(\beta \alpha) \gamma  = f(\beta \alpha) \cdot \gamma$.
So $T \beta . f(\alpha) = f(\beta \alpha)$: that is the  commutativity of the bottom square.

It follows  from the Lemma of \S \ref{sec:diamond}  that the composite $U \circ U'$ is surjective.  
  Also the map $V$ is surjective (because the cohomology of $\Hom_{\Sinf}(\Cinf, \Z_p)$ is torsion-free, by \eqref{logstructure} and
  assumption 7(a) of \S \ref{Sec:assumptions}). 
    Tracing through the above diagram, this is enough to show that the image of  $\mathfrak{Q}$ generates the codomain of $\mathfrak{Q}$.

Now,  recall  from \eqref{Descent} the quasi-isomorphism $\Cinf \otimes_{\Sinf} S_n  \simeq C_n$; 
    we have therefore shown that 
   $   H^*(\Hom_{S_n}(C_n, \Z/p^n))$  is generated  by $ H^{\BWq}(\Hom_{S_n}(C_n, \Z/p^n))$ as a module   over $\Ext^*_{S_n}(\Z/p^n, \Z/p^n) $.
   As in the discussion after \eqref{above}, this  is equivalent to saying that
 $H^*(Y_0(Q_n), \Z/p^n)_{\mathfrak{m},\Frob^T}$ is generated by  its degree $q$ component as a $H^*(T_n, \Z/p^n)$-module, i.e.  
\begin{equation} \label{bencumber}H^{\BWq}(Y_0(Q_n), \Z/p^n)_{\mathfrak{m}, \Frob^T} \otimes H^*(T_n,\Z/p^n) \twoheadrightarrow H^*(Y_0(Q_n),\Z/p^n)_{\mathfrak{m}, \Frob^T}.\end{equation}

\label{randompageref}
In what follows, let us write $\mathscr{H}_I, \mathscr{H}_K$ for the tensor product \index{$\mathscr{H}_I$} \index{$\mathscr{H}_K$} 
of (derived) Iwahori-Hecke algebras  $\mathscr{H}_{I,q}$ and derived Hecke algebra $\mathscr{H}_{q}$  over  $q \in Q_n$;
and write $\mathrm{H}_K, \mathrm{H}_{KI}, \mathrm{H}_{IK}, \mathrm{H}_I$ for the (non-derived) algebras and bimodules
of \S \ref{Morita}, but tensoring over all $q \in Q_n$.  All of these will be taken with $\Z/p^n$ coefficients. 
 
Note that the action of $H^*(T_ n, \Z/p^n)$  on $H^*(Y_0(Q_n), \Z/p^n)$ factors through the action of $\mathscr{H}_I$ (e.g. see the Remark in \S \ref{sec:concrete}). 
  So $H^*(Y_0(Q_n), \Z/p^n)_{\mathfrak{m}, \Frob^T}$ is    generated in degree  $\BWq$ over  $\mathscr{H}_I$. 
  Taking the sum over all possible lifts $\Frob^T$,  as in \eqref{bigsplit}, we see that
  $H^*(Y_0^*(Q_n), \Z/p^n)_{\mathfrak{m}}$
  is also generated in degree $q$ over the derived $\mathscr{H}_I$.

  Now,  each of the following maps are surjective:
\begin{eqnarray*}
H^{\BWq}(Y(1),\Z/p^n)_{\mathfrak{m}} \otimes \mathrm{H}_{KI} \twoheadrightarrow H^{\BWq}(Y_0(Q_n), \Z/p^n)_{\mathfrak{m}}, \\    H^{\BWq}(Y_0(Q_n), \Z/p^n)_{\mathfrak{m}} \otimes \mathscr{H}_I^j \twoheadrightarrow   H^{\BWq+j}(Y_0(Q_n), \Z/p^n)_{\mathfrak{m}},\\
H^{\BWq+j}(Y_0(Q_n) , \Z/p^n)_{\mathfrak{m}}\otimes \mathrm{H}_{IK} \twoheadrightarrow H^{\BWq+j}(Y(1),\Z/p^n)_{\mathfrak{m}}, \end{eqnarray*}
where the second statement is what we just proved, 
whereas the first and third statement come from Lemma \ref{needaref}.
Also there is a map  $\left( \mathrm{H}_{KI}  \otimes \mathscr{H}_I^j  \otimes \mathrm{H}_{IK} \right)  \rightarrow \mathscr{H}_K^j$ compatible with the respective actions, just arising  
from composition of Exts. We get  
$$   H^{\BWq}(Y(1), \Z/p^n)_{\mathfrak{m}}  \otimes  \mathscr{H}_K^j 
\twoheadrightarrow H^{\BWq+j}(Y(1), \Z/p^n)_{\mathfrak{m}}. $$

Passing to the limit (as in  the discussion of \S \ref{limit})    concludes the proof.  \qed

\section{The reciprocity law} \label{reciprocity}

In \S \ref{Patching2} we proved, conditional under assumptions (\S \ref{Sec:assumptions}, \S \ref{GaloisAss}) on the existence of Galois representations attached to modular forms
and other assumptions that simplify the integral situation (\S \ref{Sec:assumptions}), 
that the global derived Hecke algebra is ``big enough,'' in the sense discussed around \eqref{bigenough}.

We now turn to the question mentioned in \S \ref{reindexing}: 
we index elements of this global derived Hecke algebra by means of a certain dual Selmer group.
This is achieved in Theorem \ref{maintheorem}.  
This Theorem is not an end in itself; rather, it just gives the correct language for us to formulate the central conjecture of the paper, 
Conjecture \ref{mainconjecture}.

 \subsection{The coadjoint representations}

We are interested in the co-adjoint representation, i.e. the dual of the 
representation of $G^{\vee}$ on its Lie algebra $\mathfrak{g}^{\vee}$.  
Denoting by $\dualLiedual$ the $\Z$-dual to this Lie algebra, we obtain  \index{$\Ad^*$} 
$$\Ad^*: G^{\vee} \rightarrow GL(\dualLiedual)$$ 
which we regard as a morphism of algebraic groups over $\Z$. 
\footnote{Why the {\em coadjoint} representation rather than the adjoint? They are isomorphic
for $G$ semisimple, at least away from small characteristic. However, canonically what comes up for us is the coadjoint;
for example,  when one works with tori, as in \S \ref{tori}, the difference is important. }

In particular given a representation $\sigma: G_{\Q} \rightarrow G^{\vee}(R)$ we denote by   
$\Ad^* \sigma: G_{\Q} \rightarrow \GL(R \otimes_{\Z} \dualLiedual)$ 
the composition of $\sigma$ with the co-adjoint representation.  
When $\sigma$ is valued in $\Z_p$, we will write
$\Ad ^* \sigma_n$ for the reduction of $\Ad ^* \sigma$ modulo $p^n$.

\subsection{Galois cohomology} 

We will freely use the theory of Fontaine and Laffaille which (in good circumstances) parameterizes
crystalline representations of $\Gal(\overline{\Q_p}/\Q_p)$, even with torsion coefficients. For a summary, see \S 4 of \cite{BK}. 

Fix once and for all the interval $[-\frac{p-3}{2}, \frac{p-3}{2}] \subset \mathbb{N}$ of Hodge weights. 
We will say that a  representation   of $\Gal(\overline{\Q_p}/\Q_p)$  on a finitely generated $\Z_p$-module is  ``crystalline''
if it is isomorphic to a subquotient of a crystalline representation with Hodge weights in $[-\frac{p-3}{2}, \frac{p-3}{2}]$. 
This indexing of Hodge weights is useful for  adjoint representations  which have weights symmetric around $0$.

Recall that for any $p$-torsion  crystalline $M$   \index{$H^1_f$}
we can define
$$  H^1_f(\Q_p, M) \subset H^1(\Q_p, M) $$
which classifies those extensions $M \rightarrow  ? \rightarrow 1$ which are crystalline;
it is in fact a submodule and it is identified (by Fontaine-Laffaille theory) with a corresponding Ext-group in the category
of filtered Dieudonn{\'e} modules. In particular,  this allows one to check that
$H^1_f$ is isomorphic to the cokernel of the map 
\begin{equation} \label{kitty}
\mathrm{F}^0 D(M) \stackrel{1 - \mathrm{Frob}}{\rightarrow} D(M),\end{equation}
where $D(M)$ is the associated filtered Dieudonn{\'e} module. 
Also, the kernel of  $1-\mathrm{Frob}$  in \eqref{kitty} is isomorphic to $H^0(\Q_p, M)$. In particular, for $M$ finite we have 
\begin{equation} \label{goo-goo} \# H^1_f(\Q_p, M) = \frac{ | D(M)|}{|F^0 D(M)|} \cdot \# H^0(\Q_p, M),\end{equation} 
 which can be effectively 
 used to compute the size of $H^1_f$ (note: the size of $D(M)$ and $M$ coincide).
 
We will need to know that the subspaces $H^1_f(\Q_p, M)$ and $H^1_f(\Q_p, M^*)$
(with $M^* := \Hom(M, \mu_{p^{\infty}})$ are each other's annihilators \index{$M^*$, dual of a  Galois module}
under the local duality pairing $H^1(\Q_p, M) \times H^1(\Q_p, M^*) \simeq \Q/\Z$. 
This follows from the fact they annihilate each other  (their product would come from an $\Ext^2$
in the category of Fontaine--Laffaille modules, but the relevant $\Ext^2$ vanishes by \cite[Lemma 4.4]{BK})  and a size computation
using \eqref{goo-goo}.

\subsection{Selmer groups}

Let $\Q_{\ramprimes}$ 
be the largest extension of $\Q$ unramified outside 
$\ramprimes$, and let $M$ be a module for the Galois group of $\Q_{\ramprimes}/\Q$;
thus $M$ defines an {\'e}tale sheaf on $\Z[\frac{1}{\ramprimes}]$. 
We write

$$H^1(\Z[\frac{1}{S}], M) \supset H^1_f(\Z[\frac{1}{S}], M)$$
for (respectively) the {\'e}tale cohomology of $M$ (equivalently the group cohomology
of $\Gal(\Q_{\ramprimes}/\Q)$ with coefficients in $M$),
and the subset of this group consisting of classes that are crystalline at $p$, i.e. classes
whose image in $H^1(\Q_p, M)$ lies in the subgroup $H^1_f(\Q_p, M)$ defined above.

Note that we impose {\em no} local condition on  classes in $H^1_f(\Z[\frac{1}{S}], M)$ except  for the crystalline condition  at $p$. 
 
 We will write $H^1_f(\Q, M)$ for the usual Bloch--Kato Selmer group: this is 
the subgroup of classes in $H^1(\Gal(\overline{\Q}/\Q), M)$ which
are unramified away from $p$, {\em and} crystalline at $p$.
In general, we have an inclusion $H^1_f(\Q, M) \subset H^1_f(\Z[\frac{1}{S}], M)$;
the former is more restrictive, requiring that the cohomology class be unramified
at places of $\ramprimes - \{p\}$. However, in our applications, 
$M$ will be a module such that $H^1(\Q_v, M)$ vanishes for $v \in \ramprimes - \{p\}$,
and so $H^1_f(\Q, M) = H^1_f(\Z[\frac{1}{S}], M)$.

\subsection{}  \label{mMain}

We will follow the notation of the previous section, described in \S \ref{Sec:assumptions};
in particular, we have an arithmetic manifold $Y(1) =Y(K_0)$, an automorphic
representation $\Pi$,  corresponding to a maximal ideal $\mathfrak{m}$ of the Hecke algebra;
and an associated Galois representation $\rho: G_{\Q} \rightarrow G^{\vee}(\Z_p)$. 
We let $\ramprimes$ be the set of ramified primes for $\rho$ or $K_0$, together with $p$.
 
\index{$\Vinf$} 
 Put  \begin{equation} \label{Vdef} \mathsf{V} :=   H^1_f(\Z[\frac{1}{\ramprimes}], \Ad^* \rho(1))^{\vee},\end{equation}
where  we wrote  $-^{\vee}$ for $\Hom(-, \Z_p)$. 
We will prove in Lemma \ref{DeltaDimensional} that both the $H^1_f$ above and  $\mathsf{V}$   are (under our assumptions)   free $\Z_p$-modules of rank $\delta$.

We will produce an action of $\mathsf{V}$ on   $H^*(Y(K), \Z_p)_{\mathfrak{m}}$.
To explain it, fix $\mathbf{A}$   a maximal torus of $\mathbf{G}$ and
let $q$ be a Taylor--Wiles prime of level $n$, 
equipped with an element of $T^{\vee}(k)$ conjugate to  Frobenius at $q$. 
Let $$T_q = \mathbf{A}(\F_q)/p^n.$$  \index{$T_q$} 
From this data we  will construct: \index{$\iota_{q,n}$} 
\begin{itemize}
\item[-] 
A natural embedding (\S \ref{Vnaction})  of  
 \begin{equation} \label{iotaqndef1} \iota_{q,n}: H^1(T_q,\Z/p^n) \hookrightarrow \left( \mathscr{H}^{(1)}_{q,\Z/p^n}\right)_{\mathfrak{m}} \end{equation} into the  degree $1$ component $ \mathscr{H}^{(1)}$ of the  local, full level,  derived Hecke algebra
$\mathscr{H}_{q, \Z/p^n}  $. (More precisely, we will use its completion at the maximal ideal $\mathfrak{m}$).   
\item[-]  A map \begin{equation} \label{fqndef} f_{q,n}: H^1(T_q, \Z/p^n) \rightarrow \mathsf{V}/p^n, \end{equation} \index{$f_{q,n}$}
We have already explained, in a special case, the constrution of $f_{q,n}$ in \eqref{padef}.
We briefly outline the general case: 
Given $\alpha \in H^1(T_q, \Z/p^n)$,
we obtain, by \eqref{brst}, an element  $\alpha'$ in the quotient of $H^1(\Q_q, \Ad \rho_n)$
by unramified classes;  now,  we associate to $\alpha$
the functional  sending $\beta \in H^1_f(\Z[\frac{1}{\ramprimes}], \Ad^* \rho(1))$ 
to the local pairing $\langle \alpha' , \beta_q \rangle_q  \in \Z/p^n$, where $\beta_q$ is the restriction of $\beta$ to $H^1(\Q_q, \Ad^* \rho(1))$. 
\end{itemize}

Finally recall that under our assumptions (7(a) of \S \ref{Sec:assumptions}) , $H^*(Y(1), \Z_p)_{\mathfrak{m}}$ is torsion-free;
its reduction modulo $p^n$ coincides with  $H^*(Y(1), \Z/p^n)_{\mathfrak{m}}$.
\begin{theorem} \label{maintheorem}
Let notation and assumptions be as established in \S \ref{Patching}
(in particular \S \ref{Sec:assumptions},  \S \ref{GaloisAss}). 
Let $\Vinf$ be as in \eqref{Vdef}.

There  exists a function $a: \mathbf{Z}_{\geq 1} \rightarrow \mathbf{Z}_{\geq 1}$ and 
an action of $\Vinf$ on $H^*(Y(1), \Z_p)_{\mathfrak{m}}$ 
by endomorphisms of degree $+1$ 
with the  following property:

(*) For any $n \geq 1$ and any  prime $q \equiv 1$ modulo $p^{a(n)}$,
equipped with a strongly regular  element  of $T^{\vee}(k)$ conjugate to 
$\rhobar(\mathrm{Frob}_q)$, 
 the two actions of $H^1(T_q, \Z/p^{n})$  on $H^*(Y(1), \Z/p^{n})_{\mathfrak{m}}$ coincide:
one via $f_{q,n}$  and one via $\iota_{q,n}$.
 
The property (*) uniquely characterizes the $\Vinf$ action (this is true for any function $a$). 
 
 In particular, the {\color{\changecolor} strict} global derived Hecke algebra
 contains $\mathsf{V}$, and so also the exterior algebra freely generated by $\mathsf{V}$
(the induced map from $\wedge^* \mathsf{V}$ 
 to endomorphisms of    $H^*(Y(1), \Z_p)_{\mathfrak{m}}$ is injective). 
 \end{theorem}

 Note that the
 uniqueness part of the statement is straightforward, because the condition pins down the action of $\mathsf{V}/p^r$ for arbitrarily large $r$:
by Chebotarev, 
  the images $f_{q,n}(H^1(T_q, \Z/p^n))$ generate $\mathsf{V}/p^n$ even when restricted to primes $q \equiv 1$ modulo $p^{a(n)}$;
  this basically  follows from the existence of Taylor--Wiles data (see in particular  \eqref{tsnvsurj}).

{\color{\changecolor}  Finally,   under a further ``multiplicity one'' assumption, this result is sufficient to force the whole
derived Hecke algebra to be graded commutative.  We separate this result
from the main analysis because it is inessential to our main goals and it requires this additional multiplicity one assumption.

Before we state the proposition, we note that, by the argument of \eqref{removing w}
together with \cite[Lemma 6.20]{KT}, the action of global Hecke algebra $\gHecke$ 
on cohomology induces an action of $\gHecke$ on 
  on $H^*(Y(1), \Z_p)_{\mathfrak{m}}$.

\begin{prop} \label{P853}
  If   $H^{\BWq}(Y(1), \Z_p)_{\mathfrak{m}}  = \Z_p$,
then the 
image of the full global derived Hecke algebra $\gHecke$ inside $\End \ H^*(Y(1), \Z_p)_{\mathfrak{m}}$ in fact precisely coincides with the exterior algebra generated by $\mathsf{V}$,
and in particular the global derived Hecke algebra (acting on $\mathfrak{m}$-part of cohomology) is graded commutative.
 
This same conclusion of graded commutativity holds, more generally, when there exists a semisimple  $\Q_p$-algebra $\mathsf{S}$ of   (degree-preserving) endomorphisms of $H^*(Y(1), \Q_p)$
 commuting with $\gHecke$ and such that $H^{\BWq}(Y(1), \Z_p)_{\mathfrak{m}} \otimes \Q_p$ has multiplicity one as a $\mathsf{S}$-module.   \end{prop}}

 \begin{remark}
 \begin{itemize}
\item[-] {\color{\changecolor}  A natural choice for the commuting subalgebra $\mathsf{S}$ in Proposition \ref{P853} is
  the $\Z_p$-algebra generated by the local underived Hecke algebras at all primes that are not good.
  However, the multiplicity one condition is only realistic  when there is a  {\em unique} tempered representation $\pi_{\infty}$ of $\mathbf{G}(\R)$ with nontrivial $(\mathfrak{g}, \mathrm{K})$-cohomology.
  If this uniqueness condition is satisfied, then we would likely expect the multiplicity one condition to be valid in the great majority of cases. 
 
%

Under the assumptions we are currently working (simply connected group, and working over $\Q$) only $\SL_{2n+1}$
has this property. However, this unicity of $\pi_{\infty}$ in fact applies
whenever we work over a CM base field, or for $\PGL_n$ over $\Q$. \footnote{Note that, for $\PGL_n$ with $n$ even, 
one should enlarge $\mathsf{S}$ using the action of the component group of the archimedean maximal compact.} 
The analysis of \S \ref{Patching} -- \S \ref{reciprocity}
would extend to those cases with only some  notational changes.   }

\item[-] That we get an {\em integral} isomorphism of $\gHecke$ and $\wedge^* \Vinf$, in the first statement of the Proposition, is an artifact of our simplifying hypotheses. We don't expect $\gHecke$
to be an integral exterior algebra in general, but corresponding statements should remain valid $\otimes \Q$. 
 
 In general,
 we would expect the $\gHecke$ to be $\mathbb{T} \otimes \wedge^* \mathsf{V}$, where $\mathbb{T}$
 is the usual Hecke algebra,   after tensoring with $\Q$. But here our assumptions mean that $\mathbb{T}$ is just $\Z_p$,
 and moreover that the conclusion is true integrally.
\end{itemize}
\end{remark}

\medskip

\subsection{Formulation of the conjecture}
We are now ready to state the conjecture, the formulation of which is the main point of this paper. 
(Although the trip was fun too.) 

The formulation of the conjecture itself rests on the conjecture of Langlands that associates  
to $\Pi$ a motive, or more precisely a system of motives indexed by representations of the dual group. 
 Unfortunately it is difficult to find a comprehensive account of this conjecture in print;  the reader may consult
the brief remarks in \cite{Langlands} or the appendix of \cite{PV}.

Continue with notation as in Theorem \ref{maintheorem}.   
As in    the discussion of \S \ref{sec:MV} and \S \ref{reindexing},
let $M_{\coad}$ be the motive  with $\Q$ coefficients associated to $\Pi$ and the co-adjoint representation of $G^{\vee}$, if it  exists. {\em A priori},
one may not always be able to descend the coefficients of $M_{\coad}$ to $\Q$, although we expect this is possible in most if not all cases.  (See discussion
in \cite[Appendix, A.3]{PV}).  In what follows we assume that $M_{\coad}$ can indeed be descended to $\Q$ coefficients; if not 
one can simply reformulate the conjecture by replacing $\Q$ by a field extension.

Thus there is an identification of Galois modules
$$ \mbox{ {\'e}tale realization of $M_{\coad}$}  \simeq \Ad^* \rho  \otimes \Q_p$$
 and there is 
 a regulator map   from the motivic cohomology 
 $$H^1_{\mot}(\Q, M_{\coad,\Z}(1)) \longrightarrow H^1_f(\Z[\frac{1}{S}], \Ad^* \rho(1)) \otimes  \Q_p.$$
  As in \S \ref{motivic chi}, the motivic cohomology group on the left-hand side has been restricted
to classes that extend to an integral model.
  We assume that this regulator map is an isomorphism.\footnote{
  In general, the appropriate conjecture is that the regulator, taken to the group of classes that are crystalline at $p$
  {\em and unramified at all primes in $S$}, is an isomorphism. However, by virtue of our assumptions in 
  \S \ref{GaloisAss}, part (e), it is not necessary to explicitly impose ``unramified at primes in $S$.''}
 Let $\mathsf{V}_{\Q_p} = \mathsf{V} \otimes  \Q_p$,
 and let $\mathsf{V}_{\Q}$ be those classes in $\mathsf{V}_{\Q_p}$
 whose pairing with  motivic cohomology
 lies in $\Q$. 
 
 Write $H^*(Y(1), -)_{\Pi}$ for the 
Hecke eigenspace for the character $\mathbb{T}_{K_0} \rightarrow \mathbf{Z}$
 associated with $\Pi$ (see \eqref{oort}).   Our assumptions imply that $H^*(Y(1), \Z_p)_{\Pi}
 = H^*(Y(1), \Z_p)_{\mathfrak{m}}$. 

\begin{conjecture}     \label{mainconjecture} 
 Notation as above. With reference to the action
 $$ \wedge^* \mathsf{V}_{\Q_p} \acts H^*(Y(1), \Q_p)_{\Pi}$$
 furnished by Theorem \ref{maintheorem}, the action of $\mathsf{V}_{\Q}$
 preserves $H^*(Y(1), \Q)_{\Pi}$. 
\end{conjecture}

  Some rather scant evidence is discussed in the next section  (\S \ref{pisspoor}).   As we have mentioned in the introduction, 
  much more compelling is that we have been able
to obtain numerical evidence for a coherent analog of the conjecture, in a joint work with Michael Harris.

Recall  (\S \ref{GaloisAss}) 
we assume that 
\begin{equation} \label{recall} H^0(\Q_p, \Ad \rhobar) = H^2(\Q_p, \Ad \rhobar) = 0\end{equation} 
which implies the same conclusions for $\Ad^* \rhobar(1)$. 
In  particular,  $H^1(\Q_p, \Ad^* \rho(1))$ is torsion-free and surjects
 onto $H^1(\Q_p, \Ad^* \rhobar(1))$. 
  Finally, $H^1_f(\Q_p, \Ad^* \rho(1))$ is a saturated submodule of $H^1(\Q_p, \Ad^* \rho(1))$ and we have an equality of ranks
 $$\mathrm{rank}_{\Z_p} H^1_f(\Q_p, \Ad^* \rho(1)) = \mathrm{rank}_{\mathbf{F}_p} H^1_f(\Q_p, \Ad^* \rhobar(1)),$$
 as follows from explicit computation. In particular, $H^1_f(\Q_p,\Ad^* \rho(1))$ surjects onto $H^1_f(\Q_p, \Ad^* \rhobar(1))$.  
Also   observe that, because of the assumed ``big image'' (\S \ref{GaloisAss} assumption (b)) of $\rhobar$, we have 
\begin{equation} \label{h00} H^0(\Z[\frac{1}{\ramprimes}], \Ad^* \rhobar(1)) = 0,\end{equation}
  and so $H^1(\Z[\frac{1}{\ramprimes}], \Ad^* \rho(1))$ is torsion-free. 

 \begin{lemma} \label{DeltaDimensional}  Both $ H^1_f(\Z[\frac{1}{\ramprimes}], \Ad^* \rho(1))$ and $\mathsf{V}$ are  free  $\Z_p$-modules of rank $\delta$.
 \end{lemma}
 \proof  
 
 First of all, because the Taylor--Wiles method in this case implies an $R=T$ theorem (see \eqref{first-iso})
 and  we are assuming   that the Hecke algebra is isomorphic to $\Z_p$  (see \eqref{nocong0}) we get from a tangent space computation that
 $H^1_f ( \Z[\frac{1}{\ramprimes}], \Ad \rhobar ) = 0.$
 We now apply Tate global duality to this statement. It implies  both the surjectivity of
\begin{equation} \label{sarge} H^1(\Z[\frac{1}{\ramprimes}], \Ad ^*  \rhobar(1)) \twoheadrightarrow \frac{ H^1(\Q_p, \Ad^* \rhobar (1)) }{ H^1_f(\Q_p, \Ad^* \rhobar (1)) }, \end{equation}
and  the injectivity  of
\begin{equation} \label{Vanishing} H^2(\Z[\frac{1}{\ramprimes}],\Ad^* \rhobar(1))  \hookrightarrow \underbrace{ \prod_{v \in \ramprimes} H^2(\Q_v, \Ad^* \rhobar(1)) }_{=0 \mathrm{  \ by  \ \S \ref{GaloisAss}}} \end{equation}
  so in fact $ H^2(\Z[\frac{1}{\ramprimes}],\Ad^* \rhobar(1))=0$.

The surjectivity \eqref{sarge} holds also for $\Ad^* \rho(1)$, not just the mod $p$ reduction.
This follows because $H^1(\Z[\frac{1}{\ramprimes}], \Ad^* \rho(1))$ surjects onto  $H^1(\Z[\frac{1}{\ramprimes}], \Ad^* \rhobar(1))$, by \eqref{Vanishing};
and the induced map
$$\frac{ H^1(\Q_p, \Ad^* \rho (1)) }{ H^1_f(\Q_p, \Ad^* \rho (1)) }/p  \rightarrow \frac{ H^1(\Q_p, \Ad^* \rhobar (1)) }{ H^1_f(\Q_p, \Ad^* \rhobar (1)) } $$
is an isomorphism. 

The Euler characteristic formula, taken together with \eqref{h00}, \eqref{sarge} and \eqref{Vanishing}, allows one to compute  \begin{equation}
\label{EC} \dim H^1_f(\Z[\frac{1}{\ramprimes}], \Ad^* \rhobar \ (1)) = \delta.\end{equation}
Now examine  the short exact sequences  
 
 \begin{equation} \label{feeling_groovy1}
 \xymatrix{
H^1(\Z[\frac{1}{\ramprimes}], \Ad^*  \rho (1)) \ar[r]^{p}\ar[d]^j   & H^1(\Z[\frac{1}{\ramprimes}], \Ad^* \rho  (1))\ar[d]   \ar[r]  &H^1(\Z[\frac{1}{\ramprimes}], \Ad^* \rhobar (1))   \ar[d] \\
\frac{H^1(\Q_p,\Ad^* \rho (1))}{H^1_f (\Q_p, \Ad^* \rho(1))} \ar[r]^{p}   & \frac{H^1(\Q_p,\Ad^* \rho  (1))}{H^1_f(\Q_p,\Ad^* \rho  (1))} \ar[r] &\frac{H^1(\Q_p,\Ad^* \rhobar (1))}{H^1(\Q_p,\Ad^* \rhobar (1))} \\
 }
 \end{equation} 
 Since we have seen that $j$ is onto, 
 it follows that the induced maps of vertical kernels is a short exact sequence;
that and \eqref{EC} imply that
 $$ H^1_f(\Z[\frac{1}{\ramprimes}], \Ad^* \rho(1)) \simeq \Z_p^{\delta},$$
 as claimed. 
   \qed

\subsection{ Cohomological vanishing in the  Taylor--Wiles method} \label{TWrecall} 
In  the Taylor--Wiles method we choose a set of primes $Q$ such that, with $SQ=  \ramprimes \cup Q$, we have the following properties:
\begin{itemize}
\item[(a)]  $Q$ is a Taylor--Wiles datum of some level, and
\item[(b)]  The map $H^1(\Z[\frac{1}{SQ}], \Ad \rhobar) \rightarrow \frac{H^1(\Q_p, \Ad \rhobar)}{H^1_f(\Q_p, \Ad \rhobar)}$ is surjective,
and the map $H^2(\Z[\frac{1}{SQ}], \Ad  \rhobar) \rightarrow \prod_{v \in Q} H^2(\Q_p, \Ad \rhobar )$ is injective.

\end{itemize} 
(Recall that our local assumptions at $S$ mean there is no local cohomology there: \S \ref{GaloisAss}, assumption (e)). 

Observe also that if $Q$ is such a set of primes, and $Q'$ is a further set satisfying (a) and (b), then certainly $Q \cup Q'$ satisfies (a) and (b) too.
Indeed, the cohomological criteria of (b) are equivalent to asking that \begin{equation} \label{dual-inj} H^1_f(\Z[\frac{1}{SQ}], \Ad^* \rhobar \ (1)) \rightarrow \prod_{v \in Q} H^1(\Q_v, \Ad^* \rhobar(1))\end{equation}  is injective,  and this  is stable under enlarging $Q$ (it is equivalent to the same injectivity on $H^1_f(\Z[\frac{1}{S}], -)$,
since anything in the kernel would be unramified at $Q$). 

Now, if we choose a system of such data $Q_n$ of level $n$, we can 
(by passing to a subsequence and reindexing, e.g. as in \cite[\S 13.10]{GV}) achieve 
a new sequence $Q_n$ which satisfy the ``limit properties'' of \S \ref{extraction}.  

\begin{Definition} \label{TWconvgt}  \index{convergent (Taylor--Wiles data)}  
A sequence of Taylor--Wiles data $Q_n$ of level $n$ is called {\em convergent} if:
\begin{itemize} \item[-] $Q_n$ have  the cohomological properties stated in (b) above and, 
\item[-] One can pass to the limit in the sense of \S \ref{extraction}, i.e. there exists data $\Rinf, \Sinf, \Cinf, f_n, g_n$ etc. 
satisfying all the properties enumerated in \S \ref{extraction}.
\end{itemize}
 \end{Definition}
 
 In particular, any sequence of Taylor--Wiles data has a convergent subsequence, after reindexing the subsequence.

\subsection{The tangent spaces to $R$ and $S$} \label{tangentRS}
As in \S \ref{SR} both $\Rinf$ and $\Sinf$ are augmented to $\Z_p$:
$$ \Sinf \rightarrow \Rinf \rightarrow \Z_p,$$
and the composite $\Sinf \rightarrow \Z_p$ is the standard augmentation of $\Sinf$. 
 The kernel of these augmentations are denoted by $\Jinf \subset \Rinf$ and $\Iinf \subset \Sinf$.

  First of all, we set let $\mathfrak{t}_{\Rinf}$ be the ``tangent space'' to $\Rinf$  {\em over $\Z_p$}, that is to say \index{$\mathfrak{t}_{\Rinf}$}
$$\mathfrak{t}_{\Rinf} \simeq \Hom_{*}(\Rinf, \Z_p[\varepsilon]/\varepsilon^2),$$ 
where the subscript $*$ means that the homomorphism lifts the natural augmentation $\Rinf \rightarrow \Z_p$.   Equivalently, $\Rinf$ is the derivations of $\Rinf/\Z_p$ into $\Z_p$, 
or the $\Z_p$-linear dual of $\Jinf/\Jinf^2$. 

We can make exactly the same definition for $\Sinf$.  The surjection $\Sinf \twoheadrightarrow \Rinf$ induces 
a surjection $\Iinf/\Iinf^2 \twoheadrightarrow \Jinf/\Jinf^2$ and thus 
 a natural injection $\mathfrak{t}_{\Rinf} \rightarrow \mathfrak{t}_{\Sinf}$ with saturated image (i.e. split). 
We write  $\Winf$ for the cokernel of the map on tangent spaces, so we have an exact sequence 
\begin{equation} \label{Vdef2}  \mathfrak{t}_{\Rinf} \hookrightarrow \mathfrak{t}_{\Sinf}  \twoheadrightarrow \Winf.  \end{equation} 
 Then $\Winf$ is a free $\Z_p$-module of rank $\delta$.

 \subsection{Tangent spaces to $\usbR_n$ and $S_n$}  \label{tnRnSn}
 
 We suppose now that $Q_n$ are a convergent sequence (\S \ref{TWrecall}) of Taylor--Wiles data.  
  
Recall that $\Sinf, \Rinf$ are defined as ``limits,'' roughly speaking,
of rings $S_n \rightarrow \usbR_n$ that occur at level $Q_n$ in the Taylor--Wiles process.
We recall that $\usbR_n$ is not the full (crystalline at $p$) deformation ring  $R_{Q_n}$ at level $Q_n$, but is a ``very deep'' Artinian quotient of 
it. 

 These rings are also compatibly  augmented over $\Z/p^n$, i.e.
$$S_ n  \longrightarrow   \usbR_n \longrightarrow \Z/p^n$$
(see the bottom row of  \eqref{SinftySn}, and compose with the reduction of the map $\defring \rightarrow \Z_p$
which arises from our fixed automorphic representation $\Pi$).

\begin{lemma}
The map  $R_{Q_n} \rightarrow \usbR_n$ 
induces an isomorphism  upon applying $\Hom_*(-, \Z/p^n[\varepsilon]/\varepsilon^2)$, 
where $\Hom_*$ means that the map lifts the natural augmentations to $\Z/p^n$. 
\end{lemma}
\proof
Write $A :=   \Z/p^n[\varepsilon]/\varepsilon^2$. 
A map $R_{Q_n} \rightarrow  A$ gives rise to a deformation 
 $\widetilde{\rho}_n: G_{\Q} \rightarrow G^{\vee}(\Z/p^n)$
 that lifts the modulo $p^n$ reduction of $\rho$.   %

We want to show that $R_{Q_n} \rightarrow A$ must factor through  $\usbR_n$. 
To do so we must show (see \eqref{goef})  that
the map dies on the $K(n)$th power of the maximal ideal $\mathfrak{m}_{R_{Q_n}}$,
and also that $\rho_n$ automatically has inertial level $\leq n$ at primes in $Q_n$ (see page \pageref{inertial level} for definition). 

Note that the maximal ideal $(p, \varepsilon)$ of $\Z/p^n[\varepsilon]/\varepsilon^2$
satisfies $(p, \varepsilon)^{n+1} = 0$. Since $K(n)  \geq n+1$ by assumption,  we only need verify that $\widetilde{\rho}_n$ has inertial level $\leq n$
at all Taylor--Wiles primes. 
After suitably conjugating, we can suppose $\widetilde{\rho}_n|\Q_q$
to have image in $T^{\vee}$ (see Lemma \ref{TWdeflem}).  Restricted to inertia it takes image inside 
 the kernel of $T^{\vee}(\Z/p^n[\varepsilon]/\varepsilon^2) \rightarrow T^{\vee}(\Z/p^n)$; 
 this group has exponent $p^n$, and  so $\widetilde{\rho}_n$ has inertial conductor $\leq n$. 
\qed

We also define $$\mathfrak{t}_{R_n} = \Hom_*(R_{Q_n}, \Z/p^n[\varepsilon]/\varepsilon^2).$$ 
  \index{$\mathfrak{t}_{R_n}$}    It doesn't matter whether we use $R_{Q_n}$ or $\usbR_n$ in this definition, as we just showed.
  Similarly we define $$\mathfrak{t}_{S_n} = \Hom_*(S_n, \Z/p^n[\varepsilon]/\varepsilon^2).$$  There is a natural map $\mathfrak{t}_{R_n} \rightarrow \mathfrak{t}_{S_n}$
  induced by $S_n \rightarrow \usbR_n$.  \index{$\mathfrak{t}_{S_n}$} 

Finally define  \index{$W_n$} 
\begin{equation} \label{Wndef} W_n = \mathrm{cokernel} \left( \mathfrak{t}_{R_n} \rightarrow \mathfrak{t}_{S_n}\right). 
\end{equation}

The maps $\Rinf \rightarrow \usbR_n$ and $\Sinf \rightarrow S_n$ give rise to an isomorphism 
of short exact sequences as below: 
 
 \begin{equation} \label{feeling_groovy2}
 \xymatrix{
 0 \ar[r] & \mathfrak{t}_{R_n} \ar[r]\ar[d]^{\alpha, \sim}  & \mathfrak{t}_{S_n} \ar[d]^{\beta,\sim}  \ar[r]  & W_n   \ar[d]^{\gamma, \sim} \ar[r] & 0  \\
 0 \ar[r] & \mathfrak{t}_{\Rinf}/p^n \ar[r]    &\mathfrak{t}_{\Sinf}/p^n  \ar[r] & \Winf/p^n  \ar[r] &0 \\
 }
 \end{equation} 
 
This requires some explanation.   

First of all,  we explain the maps.  Note first of all
that there is a natural map  $\mathfrak{t}_S/p^n \simeq \Hom_{*}(\Sinf, \Z/p^n[\varepsilon]/\varepsilon^2)$,
which is an isomorphism. Similarly for $\Rinf$. 
 This means that there are maps 
\begin{equation} \label{gogo}   \alpha: \mathfrak{t}_{S_n} \rightarrow \mathfrak{t}_{\Sinf}/p^n,   \beta:  \mathfrak{t}_{R_n} \rightarrow \mathfrak{t}_{\Rinf}/p^n\end{equation} 
that are induced by the projections $\Sinf \rightarrow S_n$ and $\Rinf \rightarrow \usbR_n$.
This explains $\alpha, \beta$; and the map $\gamma$ is the induced map on cokernels.  

Next we see that $\alpha, \beta, \gamma$ are isomorphisms. 

For $S$ this is the assertion that homomorphism $\Sinf \rightarrow \Z/p^n[\varepsilon]/\varepsilon^2$ 
factors through $S_n$. Indeed,  referring to the coordinate presentation \eqref{explicit} and \eqref{Sinfdef},  each element $x_i$ must go to   $a\varepsilon$ for some $a \in \Z/p^n$, 
and then $(1+x_i)^{p^n}$ is carried to $(1+a \varepsilon)^{p^n} = 1$. 
For $R$, we use the fact \eqref{second-iso}    that   we can identify $\usbR_n$ with $\Rinf\otimes_{\Sinf} S_n$. 
As above, any homomorphism $\Sinf \rightarrow \Z/p^n[\varepsilon]/\varepsilon^2$  (lifting the augmentation) factors through
$S_n$,  and in particular any homomorphism $\Rinf \rightarrow \Z/p^n$, lifting the natural one,
factors  through $\usbR_n$. 

Now, the bottom row is exact, by definition and the freeness of $\Winf$. The top row
is exact at the left because the  vertical maps $\alpha, \beta$ are isomorphisms,
and exact at the right by definition. 
Then it follows that $\gamma$ is an isomorphism too. 

This concludes the explanation of diagram \eqref{feeling_groovy2}.

\subsection{Tangent spaces  to $S_n$ reinterpreted}
Let us reinterpret the tangent space to $S_n$ in a few different  (canonical) ways. 

\label{DDR REF1}

By \eqref{coocoo} we have an isomorphism
\begin{equation}\label{ExtT} \mathfrak{t}_{S_n} \simeq \Ext^1_{S_n}(\Z/p^n, \Z/p^n)\end{equation} 

Now $S_n$ was, by definition (\S \ref{sec:diamond}) , the $\Z/p^n$- group algebra of the group $T_n$; thus, from the above equation,  we get a canonical isomorphism
\begin{equation} \label{tsnTn}\mathfrak{t}_{S_n} \simeq H^1(T_n, \Z/p^n).\end{equation}
 
 Next we connect $\mathfrak{t}_{S_n}$ to Galois cohomology. Recall from  Lemma \ref{TWdeflem}
  that, for any $q \in Q_n$, a deformation of $\rho|G_{q}$
can be conjugated to lie in the torus, and in particular factors through the profinite completion of $\Q_q^{\times}$. 
Now we can identify $H^1(G_q, \Ad \ \rho_n)$
with the set of lifts  of $\rho_n$ to $G^{\vee}(\Z/p^n[\varepsilon]/\varepsilon^2)$, modulo conjugacy. 
This lift sends tame inertia to the kernel of reduction modulo $\varepsilon$. 
In particular, having fixed an element of $T^{\vee}(k)$ conjugate to $\rhobar(\mathrm{Frob}_q)$, we get a canonical isomorphism 
\begin{equation} \label{TqH} \frac{ H^1(\Q_q, \Ad \rho_n)}{H^1_{\ur}(\Q_q, \Ad \rho_n) } \simeq 
\mathrm{Hom}\left( \ \mathbf{F}_q^{\times},  
\mathrm{Lie}(T^{\vee}) \otimes \Z/p^n\right)\end{equation}
Identifying the the Lie algebra
with $X_*(T^{\vee})$, we get 
\begin{equation}   \label{brst}
 \frac{ H^1(\Q_q, \Ad \rho_n)}{H^1_{\ur}(\Q_q, \Ad \rho_n) }  \simeq \Hom(\mathbf{F}_q^{\times}/p^n,  \Z/p^n) \otimes X_*(T^{\vee})  \simeq \Hom(\underbrace{X_*(\mathbf{A}) \otimes \mathbf{F}_q^{\times}/p^n}_{\mathbf{A}(\mathbf{F}_q)/p^n}, \Z/p^n)\end{equation} 
(here the subscript ``ur'' means unramified) and thus, from \eqref{idontcare} and  \eqref{tsnTn} 
\begin{equation} \label{tSn Galois} \mathfrak{t}_{S_n} \simeq  \bigoplus_{q \in Q_n}  \frac{ H^1(\Q_q, \Ad \rho_n)}{H^1_{\ur}(\Q_q, \Ad \rho_n)} \end{equation} 
where we emphasize that the isomorphism depends on the choice of  an element of $T^{\vee}(k)$ conjugate to Frobenius at $q$, for each $q \in Q_n$.

  There is an isomorphism similar to \eqref{ExtT} for
 $\mathfrak{t}_{\Sinf}$; in particular,
 $ \mathfrak{t}_{\Sinf} \simeq \Ext^1_{\Sinf}(\Z_p, \Z_p)$
 and more usefully
$$ \mathfrak{t}_{\Sinf}/p^n \simeq \Ext^1_{\Sinf}(\Z_p, \Z/p^n)$$

The composite
 \begin{equation} \label{compositemap} \Ext^1_{S_n}(\Z/p^n, \Z/p^n)  {\longrightarrow} \Ext^1_{\Sinf}(\Z/p^n, \Z/p^n) \rightarrow \Ext^1_{\Sinf}(\Z_p, \Z/p^n).\end{equation}
gives the natural identification  -- map $\beta$ from \eqref{feeling_groovy2}
 -- of $\Ext^1_{S_n}(\Z/p^n, \Z/p^n) = \mathfrak{t}_{S_n}$
 with $ \Ext^1_{\Sinf}(\Z_p, \Z/p^n) = \mathfrak{t}_{\Sinf}/p^n$.

\subsection{$V_n$ and Galois cohomology} \label{Vnhstar}

\label{DDR REF2}

  We  exhibit now a canonical surjection 
\begin{equation} \label{tsnvsurj}  \mathfrak{t}_{S_n} \twoheadrightarrow \Vinf/p^n\end{equation}
In fact, this surjection uses no more than the fact that $Q_n$ is a Taylor--Wiles datum. If we suppose
that $Q_n$ are a convergent sequence (\S \ref{TWrecall}) of Taylor--Wiles data,   we will see that this actually descends to an isomorphism
\begin{equation} \label{brrr} W_n  \simeq \Vinf/p^n,\end{equation}
    i.e. $W_n$ (defined in \eqref{Wndef})  is isomorphic to $ \Hom( H^1_f(\Z[\frac{1}{S}], \Ad^* \rho \ (1)  ), \Z/p^n)$.

  As before, let $SQ_n$ be the union of the set $\ramprimes$ with the set $Q_n$.   Examine: 
  \begin{equation}
 \xymatrix{
 \Hom_{*}(R_{Q_n}, \Z/p^n[\varepsilon]/\varepsilon^2) \ar[r] \ar[d]^{\sim} & \Hom_{*}(S_n, \Z/p^n[\varepsilon]/\varepsilon^2) \ar[r]\ar[d]^{\sim} &  W_n \ar[d]^{=}\\ 
  H^1_f(\Z[\frac{1}{SQ_n}], \Ad \rho_n)  \ar[r]^{\varphi} &     \prod_{v \in Q_n} \frac{H^1(\Q_v, \Ad \rho_n)}{H^1_{\mathrm{ur}}(\Q_v, \Ad \rho_n)}   \ar[r] &  W_n.%
 }
 \end{equation} 
 The first  vertical map is just the computation of tangent spaces  to deformation rings (working over $\Z/p^n$ rather than a field),  the second vertical map 
 is \eqref{tSn Galois}, and $\varphi$ is restriction in Galois cohomology.
 
 There is now a natural pairing :\begin{equation} \label{pairing}  H^1_f(\Z[\frac{1}{\ramprimes}], \ \Ad^* \rho_n \ (1)) \times W_n \longrightarrow \Z/p^n \end{equation}
which, we emphasize again, depends on the choice of toral elements conjugate to Frobenius at each prime in $Q_n$. 
To  be explicit,  
an element of $W_n$ is represented by a collection
 $$(\beta_v) \in  \prod_{v \in Q_n} \frac{H^1(\Q_v, \Ad \rho_n)}{H^1_{\mathrm{ur}}(\Q_v, \Ad \rho_n)},$$  modulo $\mathrm{image}(\varphi)$;  
to pair $\alpha \in H^1_f(\Z[\frac{1}{\ramprimes}], \Ad ^* \rho_n(1))$ with  $(\beta_v)_{v \in Q_n}$ 
 we take the sum of local pairings
\begin{equation} \label{abovemoo}   (\alpha, (\beta_v)_{v \in Q_n}) \mapsto \sum (\alpha_v, \beta_v)_v\end{equation}
where the local pairing is defined by restricting $\alpha$ to $\Q_v$ and using local reciprocity. 
This pairing \eqref{abovemoo} is well-defined because each $\alpha_v$, i.e. the restriction of $\alpha$ to $\Q_v$,  is actually unramified.   Moreover, if the collection $(\beta_v)$ come from a class in $H^1_f(\Z[\frac{1}{SQ_n}], \Ad \rho_n)$, 
the value of \eqref{abovemoo} is zero by global reciprocity: our local assumptions means that the local pairings
for $v \in \ramprimes$ vanishes.   Thus the pairing \eqref{abovemoo} descends to the quotient of
$\prod_{v \in Q_n} H^1$ by $\mathrm{image}(\varphi)$. 
 This concludes our discussion of \eqref{pairing}. 
 
We have also seen (\S  \ref{tnRnSn},  and similar arguments to Lemma \ref{DeltaDimensional}) that both $W_n$ and $H^1_f(\Z[\frac{1}{\ramprimes}], \Ad^*  \rho_n \ (1))$ are free $\Z/p^n$-modules of dimension $\delta$. 
Let us check that     \eqref{pairing} is a perfect pairing of $\Z/p^n$ modules, i.e.
the map
$$ W_n \rightarrow H^1_f(\Z[\frac{1}{\ramprimes}], \Ad^* \rho_n (1))^{\vee}$$
is an isomorphism, where $\vee$ means homomorphisms to $\Z/p^n$.   Since both sides have the same size, it is enough to check that the map  
is surjective, and thus enough to show that the induced map
\begin{equation} \label{lunytunes} \mathfrak{t}_{S_n} \rightarrow H^1_f(\Z[\frac{1}{\ramprimes}], \Ad^* \rhobar \ (1))^{\vee}\end{equation} 
is surjective, where $\vee$ now means homomorphisms to $\Z/p$. 

  Now the Taylor-Wiles set $Q_n$ is chosen \eqref{dual-inj} so that 
 $ H^1_f(\Z[\frac{1}{\ramprimes}],  \Ad^* \rhobar \ (1)) \hookrightarrow \prod_{v \in Q_n} H^1(\Q_v, \Ad^* \rhobar(1))$.
 The image of this map consists of classes unramified at $Q_n$, so we also have
  $$H^1_f(\Z[\frac{1}{\ramprimes}], \Ad^* \rhobar\ (1)) \hookrightarrow \prod_{v \in Q_n} H^1(\Q_v, \Ad^* \rhobar(1))^{\mathrm{ur}}$$
When we dualize this, and apply local duality at primes in $Q_n$,  we get the surjectivity of \eqref{lunytunes}.

 \subsection{The injection $H^1(T_q) \hookrightarrow \mathscr{H}_q^{(1)}$} \label{Vnaction}

 In this section, we suppose that $Q_n$ is a Taylor--Wiles datum of level $n$, but do not assume  that it is part of a convergent sequence (\S \ref{TWrecall}) of Taylor--Wiles data. 
 Let $q \in Q_n$, thus equipped with $\Frob_q^T \in T^{\vee}(k)$.   We will work exclusively with $S=\Z/p^n$ coefficients. 
   As before, $\mathscr{H}_q$ denotes the local derived Hecke algebra at $q$, and $\mathscr{H}_q^{(1)}$ denotes its degree $1$ component. 
   Let
   $T_q = \mathbf{A}(\F_q)/p^n.$

 We are going to describe the map \eqref{iotaqndef1}, which is necessary for the formulation of the Theorem. 
More precisely, we are going to describe a map \begin{equation} \label{thetaref} \theta: H^1(T_q) \rightarrow  \left( \mathscr{H}_q^{(1)}\right)_{\mathfrak{m}},\end{equation}
  where the subscript means that we complete at the ideal of the Hecke algebra  (i.e., the degree zero component of $\mathscr{H}_q$) 
 corresponding to $\mathfrak{m}$. 
 
 The easiest way to think about $\theta$ is probably through  the following property: 
for each $\alpha \in H^1(T_q)$, the action of $\theta(\alpha)$ on $H^*(Y(1), \Z/p^n)_{\mathfrak{m}}$ is thus: 
   \begin{equation}  \label{explicit-action} \mbox{Pullback to $Y_0(q)$, project
 to $\Frob^T_q$-eigenspace,  cup with  $\alpha$,  pushdown to $Y(1)$. } \end{equation} 
 where the projection is done
 with reference to the splitting 
 of Corollary \ref{MoritaCor}.  
 
  The formal definition of the map $\theta$ is given in
\eqref{embedH1T}, and the validity of \eqref{explicit-action} will follow from the Lemma below.

 \medskip
 
 The Satake isomorphism of \S \ref{Satake} gives
 $$\mathscr{H}_q \stackrel{\sim}{\longrightarrow}  \left(  S[X_*] \otimes  H^*(\mathbf{A}(\F_q) ) \right)^W$$
 and in particular with $\Z/p^n$ coefficients the map $H^1(T_q) \rightarrow H^1(\mathbf{A}(\F_q))$ is an isomorphism, so 
  $$\mathscr{H}_q ^{(1)}\stackrel{\sim}{\longrightarrow}  \left(  S[X_*] \otimes  H^1(T_q ) \right)^W$$

 Now $\Frob_q^T \in T^{\vee}(k)$ gives a map $X_*(\mathbf{A}) = X^*(T^{\vee}) \rightarrow k^{\times}$, 
 i.e. it gives rise to a character $\chi_{\Frob_q^T}: S[X_*] \rightarrow k$.  The pullback of this to $S[X_*]^W$
 defines the maximal ideal $\mathfrak{m}$ (using the Satake isomorphism). 
 Let us denote by $\widetilde{\mathfrak{m}}$ the extension of the ideal $\mathfrak{m}$ back to $S[X_*]$;
 we caution that it is  no longer maximal, and rather it cuts out $\Frob_q^T$ together with all its $W$-conjugates. 
 We have an identification of completions
\begin{equation} \label{wsplit} S[X_*]_{\widetilde{\mathfrak{m}}} \cong \bigoplus_{w \in W} S[X_*]_{w \chi},\end{equation}
where we have denoted by $S[X_*]_{w\chi}$ the completion of $S[X_*]$ at the maximal ideal
that is the kernel of $w \chi$.

Next 
the natural inclusion  
\begin{equation} \label{inc}   \left( S[X_*] \otimes H^*(T_q) \right) ^{W} \hookrightarrow S[X_*] \otimes H^*(T_q) \end{equation}
induces the first map of  
\begin{equation} \label{composite iso}  \left( S[X_*] \otimes H^*(T_q) \right) ^{W}_{\mathfrak{m}} \stackrel{\eqref{inc}}{   \rightarrow}   S[X_*]_{\widetilde{\mathfrak{m}}} \otimes H^*(T_q) \stackrel{\eqref{wsplit}}{\rightarrow}  S[X_*]_{\chi} \otimes H^*(T_q). \end{equation}  
The composite map of \eqref{composite iso} is an isomorphism, and thus by composing with the Satake isomorphism we get an isomorphism
\begin{equation} \label{ciso2} \left( \mathscr{H}_q^{(1)} \right)_{\mathfrak{m}} \stackrel{\sim}{\rightarrow} S[X_*]_{\chi} \otimes H^1(T_q). \end{equation}

 We then define the map $\theta: H^1(T_q) \rightarrow \left( \mathscr{H}_{q}^{(1)}\right)_{\mathfrak{m}} $ by the rule  
\begin{equation} \label{embedH1T} \theta:  h \in H^1(T_q) \mapsto 1 \otimes h \in S[X_*]_{\chi} \otimes H^1(T_q) \stackrel{\eqref{ciso2}^{-1}}{\longrightarrow}  \left( \mathscr{H}_{q}^{(1)}\right)_{\mathfrak{m}} \end{equation} 
where the superscript $1$ refers to cohomological degree, and the final map is the inverse of \eqref{ciso2}.
   Note that this embedding depends on the choice of $\Frob^T_q$; if we 
replace $\Frob^T_q$ by $w \Frob^T_q$ then the embedding is modified by means of the natural action of $w$ on $T_q$.

This concludes the description of the map $H^1(T_q) \rightarrow \mathscr{H}_q^{(1)}$.  We now want    to justify \eqref{explicit-action}.     For this
we will describe an explicit preimage of $\theta(\alpha)$ under the map  
\begin{equation} \label{blastfromthepast} \left( \mathrm{H}_{KI}  \otimes \mathscr{H}_I  \otimes \mathrm{H}_{IK} \right)  \rightarrow   \mathscr{H}_K \end{equation}
in the case of the group $\G(\Q_q)$, where notation is as before (see, e.g. the proof of Theorem \ref{mindegreeproof});
note that $\mathscr{H}_K$ coincides
with what was previously called $\mathscr{H}_q$. 
Observe that everything here is a compatibly a module under the center of $\mathrm{H}_{II}$,   which 
is identified (\S 
\ref{Zdis})
with $S[X_*]^W$. In particular, it makes sense to complete at $\mathfrak{m}$.

Now let $e_{\lambda} \in \mathrm{H}_{II}$ correspond to the characteristic function of $I_q \lambda  I_q$, where $\lambda \in X_*$.
Then $\lambda \mapsto e_{\lambda}$ defines an embedding $S[X_*] \rightarrow \mathrm{H}_{II}$, and, completing,
an embedding
\begin{equation} \label{completed embedding} S[X_*]_{\widetilde{\mathfrak{m}}} \rightarrow (\mathrm{H}_{II})_{\mathfrak{m}} \end{equation}

As before (\S \ref{Iwahori}) we let $e_K$ be   the characteristic function of $K_q$, divided by its measure; 
if we write $e_w$ for the characteristic function of $I_q w I_q$, we have $e_K = \frac{1}{|W|} \sum_{w \in W} e_w$. 
It can be considered, as in \S \ref{Iwahori}, as an element of $\mathrm{H}_{KI}$ or $\mathrm{H}_{IK}$. 
Also $|W|   e_K e_{\lambda}   \in \mathrm{H}_{KI}$ 
{\color{\changecolor} (i.e.  the product of $|W| e_K \in \mathrm{H}_{KI}$ with $e_{\lambda} \in \mathrm{H}_{II}$)} corresponds to the characteristic function of $I_q \lambda K_q$.  

\begin{lemma} \label{Theta explicit}
 Let $\Theta  \in  S[X_*]_{\widetilde{\mathfrak{m}}}$
 be chosen so that it projects under \eqref{wsplit} to the identity in $S[X_*]_{\chi}$
 and to zero in all $S[X_*]_{w \chi}$, for $w \in W$ not the identity.  
We use the same letter for its image $\Theta \in  (\mathrm{H}_{II})_{\mathfrak{m}}$
 under \eqref{completed embedding}.

For $h \in H^1(T_q)$ let $\langle h \rangle$ be the associated element of $\mathscr{H}_{I}^{(1)}$, i.e.
$\langle h \rangle$ is supported on the identity double coset of $I$, and 
the associated cohomology class is obtained from $h$ by means of the restriction isomorphism
$H^1(I) \stackrel{\sim}{\rightarrow} H^1(T_q)$.

Then the $\mathfrak{m}$-completion of
 \eqref{blastfromthepast} sends  $|W|  e_K \Theta \otimes \langle h  \rangle \otimes  e_K$  to $\theta(h)$.
  \end{lemma}

  {\color{\changecolor} As above, the product $e_K \Theta$ is understood as the product of $e_K \in \mathrm{H}_{KI}$ with $\Theta \in \mathrm{H}_{II}$. }
   In words, this amounts precisely to the description \eqref{explicit-action}   for $\theta$,  taking into account that $\Theta$ realizes
 precisely the  projection onto the $\Frob_q^T$ component for the splitting \eqref{bigsplit}.  
  
 \proof

Consider the map (no completions, at the moment)    
\begin{equation} \label{compo}  \left( \mathrm{H}_{KI}  \otimes \mathscr{H}_I^{(1)}  \otimes \mathrm{H}_{IK}  \right)  \rightarrow   \mathscr{H}_K^{(1)} \rightarrow \left(  S[X_*] \otimes  H^1(T_q ) \right)^W\end{equation} 
  We show it sends $A := e_{K} e_{\lambda} \otimes \langle h \rangle  \otimes e_K$, for $\lambda$ dominant,   to the ``$W$-average'' of $\lambda \otimes h$,
i.e. $ |W|^{-1} \sum_{w \in W} w \cdot (\lambda \otimes h)$.%
  
   The claimed result will follow easily from this: Consider an element $\Theta' = \sum c_{\lambda} \lambda \in S[X_*]$
 with the property that $c_{\lambda}$ is nonzero only for $\lambda$ dominant, and also $\Theta \equiv \Theta'$ modulo some high power of $\widetilde{\mathfrak{m}}$. 
Our claim implies that  $|W|   e_K \Theta'  \otimes \langle h \rangle \otimes e_K$ is sent to the  sum of Weyl translates of $\sum c_{\lambda} (\lambda \otimes h)$. 
The image of $\sum c_{\lambda} (\lambda \otimes h)$ in $\bigoplus_{w \in W} S[X_*]_{w \chi} \otimes H^1(T_q)$  is
 very close to  $1 \otimes h$ in the $w=1$ factor, and very close  to zero in the other factors; 
 after summing over $W$, its projection  
 to the $w=1$ factor remains very close  to $1 \otimes h$.
  Here ``very close'' is taken in the topology
 of the complete local rings $S[X_*]_{w \chi}$. 
In other words, $\theta(h)$ and the image of $|W|   e_K \Theta'  \otimes \langle h \rangle \otimes  e_K$ 
 under \eqref{blastfromthepast}  and \eqref{ciso2} are very close;
 in the limit, this  shows the desired result. 
 
 We will now consider  everything in the ``function model'' of \S \ref{desc2}. %
 Let $a_1 \in \mathrm{H}_{KI}, a_2 \in \mathrm{H}_{IK}$ be the images of $e_K e_{\lambda} , e_K$ in the function model. Then 
 $a_1$ corresponds to  the function sending $(xK, yI)$ to $|W|^{-1}$ precisely when $I y^{-1} x K= I \lambda K$, 
 equivalently $K x^{-1} y I = K \lambdabar I$, where we write $\lambdabar \in X_*$ for the negative of $\lambda$.
 Also $a_2$ 
 corresponds to  the function $(xI, yK)$ which is $1$ exactly when $xK = yK$.   
 Moreover, the function $\langle h \rangle$ is supported on the diagonal in $G/I$ and sends $(I, I)$ to $h \in H^1(T_q) \simeq H^1(I)$. 
  
 The second map of \eqref{compo} is  given by restricting arguments to the torus,
 and restricting cohomology classes to $T_q$. 
   We can compute this restriction using the ``localization'' results of \S \ref{dhLocalization};
   these results assert that  restriction to the torus actually preserves multiplication. 

We restrict to the torus. Suppose that $\mu \in X_*$  is  dominant, with negative $\mubar \in X_*$.  We compute: 

 $$A (K,\mubar K) = \sum_{  y \in   \tilde{W}} a_1(K,yI) \langle h \rangle (yI,yI)  a_2(yI, \mubar K)$$
 The first term is nonzero only 
 for those $y$ satisfying $K y I = K \lambdabar I$, i.e. $y \in \tilde{W} \cap K  \lambdabar I$. 
This implies that
  $y = w \lambdabar$, with $w \in W$. 
So this equals 
 $$ = |W|^{-1} \sum_{w \in W} \langle h \rangle (w \lambdabar I, w \lambdabar I) a_2(w \lambdabar I, \mubar K)  $$
 The final term is nonzero (and equals $1$) exactly when $\mubar K  =  w \lambdabar K$.    Recall  that both $\lambda$ and $\mu$ are dominant.
  Thus  this only happens if $\lambda = \mu$ and $w \in W_{\mu}$: 
\begin{equation} \label{Derived} A(K,  \mubar K) =   |W|^{-1} \delta_{\lambda \mu}  \sum_{w \in W_{\mu}} \langle h \rangle (w \lambdabar I, w \lambdabar I) \end{equation} 
   On the right, $\langle h \rangle (w \lambda I, w \lambda I) \in H^1(T_q)$
equals   $w \cdot h \in H^1(T_q)$.  

Now, the image of $A$ under the derived Satake map
is the $W$-invariant element of the derived Hecke algebra of the torus
defined by
  $$\mu \mapsto \mbox{restriction to $T_q$ of }A(\mu K,   K).$$
As we have just computed,   for $\mu$ dominant, $A(\mu K, K) = A(K, \mubar K)$   is nonvanishing only when $\mu=\lambda$, where its value is $|W|^{-1} \sum_{w \in W_{\mu}} wh$. 
Therefore the image of $A$ under \eqref{compo} is the $W$-average of $\lambda \otimes h$ as claimed. 
\qed

\subsection{Producing an action of $\Vinf/p^n$ on the cohomology at level $1$} \label{tSnaction} 
Let $Q_n$ be a convergent sequence of Taylor--Wiles data (\S \ref{TWrecall}). For each integer $n$ we 
will produce an action of $\Vinf/p^n$
on   automorphic cohomology $H^*(Y(1), \Z/p^n)_{\mathfrak{m}}$. 
 {\em A priori} these actions will not be guaranteed to be compatible with one another; later
 we will see at least that they ``converge as $n \rightarrow \infty$'' to give an action of $\Vinf$ on $H^*(Y(1), \Z_p)_{\mathfrak{m}}$. 

More exactly, we
begin by constructing an action of $\mathfrak{t}_{S_n}$, then prove (the Lemma below)
that it factors through $W_n$, and finally we have identified $W_n \simeq \Vinf/p^n$ in \eqref{brrr}. 
This gives the desired action, and we will discuss the ``convergence as $n \rightarrow \infty$'' in the next section.

Thus, let $Q_n$ be a convergent sequence of Taylor--Wiles data.

We have   \begin{equation} \label{embed} \mathfrak{t}_{S_n} \stackrel{\eqref{tsnTn}}{\cong} H^1(T_n, \Z/p^n) \stackrel{\eqref{thetaref}}{\hookrightarrow} \mbox{degree $1$ component of }\otimes_{q \in Q_n}  \left(  \mathscr{H}_q\right)_{\mathfrak{m}} ,\end{equation}
where   the Hecke algebras are taken with $\Z/p^n$ coefficients. The composite embedding will be denoted
\begin{equation} \label{iotaQndef}  \iota_{Q_n}: \mathfrak{t}_{S_n} \rightarrow  \otimes_{q \in Q_n} \left(  \mathscr{H}_q\right)_{\mathfrak{m}}.\end{equation}

This gives rise to an  action of $ \mathfrak{t}_{S_n}$
by degree $+1$ endomorphisms of  automorphic cohomology
$H^*(Y(1), \Z/p^n)_{\mathfrak{m}}$, whose explicit description is essentially that already given in 
\eqref{explicit-action}, just replacing the role of one prime by many.

The embedding \eqref{tsnTn} and so also this action depends on the choice of elements $\Frob_q^T \in T^{\vee}(k)$ for each prime $q \in Q_n$. 
Should we modify $\Frob_q^T$ by an element $w_q \in W$, the Weyl group, the action of $\mathfrak{t}_{S_n}$
is also modified (see  comments after \eqref{embedH1T})  by the action of $w_q$  in the obvious way.  
 
 \begin{lemma} 
   \label{Vnfactor}
   Let $Q_n$ be a convergent sequence of Taylor--Wiles data, as in Definition \ref{TWconvgt}.    
   Then, for each $n$, 
the action of $\mathfrak{t}_{S_n}$ on $H^*(Y(1), \Z/p^n)_{\mathfrak{m}}$  (via $\iota_{Q_n}$)  is trivial on the image of $\mathfrak{t}_{R_n}$, and thus  factors through the map $\mathfrak{t}_{S_n} \rightarrow  W_n$. 
  \end{lemma}
  
  \proof 
 Consider the diagrams \eqref{previousstory} of our previous story, now with identifications with $\mathfrak{t}_S$ and $\mathfrak{t}_{S_n}$ included: 
 \begin{equation} \label{BigDog}
 \xymatrix{
 H^*(\Hom_{\Sinf}(\Cinf, \Z_p))^{} \ar[d]^{=} & \times  &   \protect\overbrace{ \Ext^1_{\Sinf}(\Z_p, \Z_p) }^{\mathfrak{t}_{\Sinf}}\owns \tilde{v} \ar[d]  \ar[r] &  H^*(\Hom_{\Sinf}(\Cinf, \Z_p)) \ar[d] \\ 
   \tilde{x} \in H^*(\Hom_{\Sinf}(\Cinf, \Z_p))   \ar[d]^{V}  &  \times&  \Ext^1_S(\Z_p, \Z/p^n)     \owns  v \ar[r]  &   H^*(\Hom_{\Sinf}(\Cinf,  \Z/p^n))\ar[d]^{=}  \\
x \in  H^*(\Hom_{\Sinf}(\Cinf, \Z/p^n))^{}   & \times & \left( \Ext^1_{\Sinf}(\Z/p^n, \Z/p^n) \right) \ar[r]   \ar[u]^U  &  H^*(\Hom_{\Sinf}(\Cinf, \Z/p^n))   \\
 H^*(\Hom_{S_n}(\Cinf \otimes_{\Sinf} S_n, \Z/p^n))^{} \ar[u]^{\sim} & \times  &  \protect\underbrace{ \left( \Ext^1_{S_n}(\Z/p^n, \Z/p^n) \right) }_{\mathfrak{t}_{S_n}} \ar[r] \ar[u] &  \ar[u]^{\sim} H^*(\Hom_{S_n}(\Cinf \otimes_{\Sinf} S_n, \Z/p^n)). \\
  }
  \end{equation}

   In particular,   let $x \in  H^*(\Hom_{S_n}(\Cinf \otimes_{\Sinf} S_n, \Z/p^n)) $ be liftable to $\tilde{x} \in  H^*(\Hom_{\Sinf}(\Cinf, \Z_p))$;
   let $v \in \mathfrak{t}_{S_n} \simeq \mathfrak{t}_{\Sinf/p^n}$ be lifted to $\tilde{v} \in \mathfrak{t}_S$. Then the image of $(x,v)$
   in the bottom row is obtained from projecting the image of $(\tilde{x}, \tilde{v})$ at the top row. 

   Let us recall from \S \ref{extraction} part (g)  that $\mathsf{C}$ is quasi-isomorphic to $\mathsf{R}$ as an $\mathsf{S}$-module in a single degree.  
   Thus we can explicitly compute what goes on in the top row. This explicit computation
   (see Lemma  \ref{KoszulExtComputation} in Appendix \S \ref{AppendixB}) shows that 
   any element   $\tilde{v} \in \mathfrak{t}_{\Sinf}$
that lies in the image of $\mathfrak{t}_{\Rinf}$  acts {\em trivially}  on $H^*(\Hom_{\Sinf}(\Cinf, \Z_p))$.

 For later use, note that these  explicit computations also show that 
\begin{equation} \label{watney} \mbox{ $H^*(\Hom_{\Sinf}(\Cinf, \Z_p))$ 
is free over
$\wedge^* \mathfrak{t}_{\Sinf}/\mathfrak{t}_{\Rinf}$.} \end{equation}

From 
\eqref{feeling_groovy2}
it then follows that 
 $\mathrm{image}(\mathfrak{t}_{R_n} \rightarrow \mathfrak{t}_{S_n})$  
acts trivially on  $$H^*(\Hom_{S_n}(\Cinf \otimes_{\Sinf} S_n, \Z/p^n)) \simeq H^*(Y_0(Q_n), \Z/p^n)_{\mathfrak{m}, \Frob^T},$$
where the action of $\mathfrak{t}_{S_n} = H^1(T_n, \Z/p^n)$  is by cup product, as in \eqref{bencumber}. 

By \eqref{explicit-action} this means that the action of $\iota_{Q_n}(\mathfrak{t}_{S_n})$
on $H^*(Y(1), \Z/p^n)_{\mathfrak{m}}$   is trivial on  $\mathrm{image}(\mathfrak{t}_{R_n} \rightarrow \mathfrak{t}_{S_n})$. 
 Thus, this action of $\mathfrak{t}_{S_n}$ on $H^*(Y(1), \Z/p^n)_{\mathfrak{m}}$ factors  
   through $W_n$ as claimed.  \qed

\subsection{Summary} \label{Strictsequence}
Let us summarize more carefully what we have said to date:

For any Taylor--Wiles datum $Q_n$ we have an action of $\mathfrak{t}_{S_n} = \Hom(T_n, \Z/p^n)$ on $H^*(Y(1), \Z/p^n)_{\mathfrak{m}}$  
constructed via an embedding $$\iota_{Q_n}:\Hom(T_n, \Z/p^n) \hookrightarrow \mbox{derived Hecke algebra} $$
 (see \eqref{embed}). 
On the other hand, we have a  surjective morphism  (see \S \ref{Vnhstar}) 
$$ f_{Q_n}: \Hom(T_n, \Z/p^n) \twoheadrightarrow \Vinf/p^n.$$ 
 These constructions, for a given $n$, depend only on $Q_n$; they do not involve the Taylor--Wiles limit process.

\begin{Definition} \label{strictdef}  \index{strict Taylor--Wiles datum}
We say that a Taylor--Wiles datum $Q_n$ of level $n$ is {\em strict of level $n$} (or just {\em strict})
if the map   
$$ \mathfrak{t}_{S_n} \stackrel{ \iota_{Q_n}}{\longrightarrow} \End \ H^*(Y(1), \Z/p^n)_{\mathfrak{m}} $$
   factors through $f_{Q_n}$. Thus, a strict Taylor-Wiles datum of level $n$
gives rise to an action of $\Vinf/p^n$ on $H^*(Y(1), \Z/p^n)_{\mathfrak{m}}$.  
\end{Definition}

What we have proved, then, amounts to the following:
\begin{lemma}  \label{Convergent implies strict} 
If the $Q_n$ are a convergent sequence of Taylor--Wiles data (Definition \ref{TWconvgt}), then 
each $Q_n$ is strict, in the sense of Definition \ref{strictdef}. 
\end{lemma}

Note we do not know that the resulting actions  of $\Vinf/p^n$ are compatible for different $n$, in any sense.

\subsection{Dependence of our construction on  choices} \label{obstructions}
We now study dependence on choices.  Using the results of this \S, we will conclude the proof of  Theorem  \ref{maintheorem}
in \S \ref{finalesec}. 

First we discuss a minor point,  the choice of $\Frob_q^T$s: 
Suppose we choose two different such choices for a given set $Q_n$,   differing by the action of $w \in W_s$.  (Recall from 
\S \ref{TWprimes stuff} that $W_s$ is just a product of copies of the Weyl group, one copy for each prime in $Q_n$). 
Then  the actions of $\mathfrak{t}_{S_n} $ on cohomology differ
by the action of $w \in W_s$ (comment after \eqref{embedH1T}). 
Also, $w : \mathfrak{t}_{S_n} \rightarrow \mathfrak{t}_{S_n}$  is compatible with the pairings previously constructed, i.e. 
this diagram commutes
 \begin{equation} \xymatrix{
\mathfrak{t}_{S_n}  \ar[d]^{w} & \times  & H^1_f(\Z[\frac{1}{S}], \ \Ad^* \rho_n \ (1))  \ar[d]^{=}   \ar[r]^{ \ \ \qquad \Frob_q^T}  & \Z/p^n  \\ 
\mathfrak{t}_{S_n} & \times  & H^1_f(\Z[\frac{1}{S}], \ \Ad^* \rho_n \ (1))    \ar[r]^{\ \ \qquad w \Frob_q^T} & \Z/p^n  \\ 
  } \end{equation}
 This shows that the action of $\Vinf/p^n$ on $H^*(Y(1), \Z/p^n)_{\mathfrak{m}}$
 did not depend on the choice of $\Frob_q^T$s for $q \in Q_n$.

We now discuss the more serious issue of choice of Taylor--Wiles data. 
 
 \begin{lemma}  \label{EqualityOfTwo}
Given  two sequences $Q_n, Q_n'$ of {\em strict} Taylor--Wiles data, 
there is a subsequence $\mathcal{J}$ of  the integers with following property:

For each $k \geq 1$,  there is $j_0$ such that, for  each $j \in \mathcal{J}, j \geq j_0$,  the  two actions $\Vinf$ on $H^*(Y(1),  \Z/p^{k})_{\mathfrak{m}}$
-- arising from reducing modulo $p^k$ the ``$Q_{j}$-action" and the ``$Q_{j}'$-action"  -- coincide with one another.\footnote{Recall that we are supposing that $H^*(Y(1), \Z_p)$ is
free over $\Z_p$.}

\end{lemma}

 \proof 
 It will be convenient to relabel the sequences of strict Taylor--Wiles data as $Q_n^{(1)},  Q_n^{(2)}$.  It will be harmless to suppose that  the sets of primes underlying 
  $Q_n^{(1)}$ and $Q_n^{(2)}$ are disjoint (otherwise, we can  e.g. just compare both of them with a third set, disjoint from both of them). 

We will compare them both to $Q_n := Q_n^{(1)} \coprod Q_n^{(2)}$  (with the obvious choice of $\Frob_q^T$ for $q \in Q_n$). 
  Of course $Q_n$ is bigger than either $Q_n^{(1)}$ or $Q_n^{(2)}$. However it still a
  sequence of Taylor--Wiles data.   
  Let $T_n, T_n^{(1)}, T_n^{(2)}$ be the analogues of $T_n$ for $Q_n, Q_n^{(1)}, Q_n^{(2)}$ respectively;
  then $T_n = T_n^{(1)} \times T_n^{(2)}$,
  and correspondingly $\mathfrak{t}_{S_n} = \mathfrak{t}_{S_n}^{(1)} \oplus \mathfrak{t}_{S_n}^{(2)}$,
  where $\mathfrak{t}_{S_n} = \Hom(T_n, \Z/p^n)$, etc.

 We have a diagram 
 \begin{equation} \label{tvlift} 
 \xymatrix{
\mathfrak{t}_{S_n}^{(1)} \ar[r] \ar[rd]^{\alpha^{(1)}}   &     \mathfrak{t}_{S_n}  \ar[d]^{\alpha}     &  \mathfrak{t}_{S_n}^{(2)} \ar[ld]^{\alpha^{(2)}} \ar[l].\\
  & H^1_f(\Z[\frac{1}{S}], \Ad^* \rho_n(1))^{\vee}  & 
 }
 \end{equation} 
 where all the $\alpha$-maps are as in   \S \ref{Vnhstar}.

 The upper maps are compatible for the actions on cohomology previously defined (\S \ref{tSnaction}),
 and everything maps compatibly to the bottom group $H^1_f(\Z[\frac{1}{S}], \Ad^* \rho_n(1))^{\vee} $. Moreover the  action of $\mathfrak{t}_{S_n}^{(1)}$ and $\mathfrak{t}_{S_n}^{(2)}$ on  
mod $p^n$ cohomology factors through the bottom row  by the assumed strictness. However, we do not know that the action of $\mathfrak{t}_{S_n}$ factors through $\alpha$.

What is missing is control of the deformation ring after adding level $Q_n$. To obtain this,  we must run now the Taylor--Wiles limit process for $Q_n$. That involves passing   to a subsequence. 
In other words, all we are guaranteed is that there is a subsequence $n_j$ such 
that $(Q_{n_j}, j)$ form a  convergent sequence  of Taylor--Wiles  data of level $j$. 
It is possible that $n_j$ is very much larger than $j$. Our prior analysis of convergent data (Lemma \ref{Convergent implies strict})  implies
that the action of $\Hom(T_{n_j}, \Z/p^j)$  on mod $p^j$ cohomology of $Y(1)$
factors through $\Vinf/p^j$, or to say it explicitly: 

 \begin{quote} (*)  The action of  $\mathfrak{t}_{S_{n_j}}$   
on $H^*(Y(1), \Z/p^j)_{\mathfrak{m}}$, 
via its embedding $\iota_{Q_{n_j}}$  into the derived Hecke algebra,  followed by reduction to $\Z/p^j$ coefficients, 
factors through the map
$$ f_{Q_{n_j}} :  \mathfrak{t}_{S_{n_j}} \rightarrow H^1_f(\Z[\frac{1}{\ramprimes}], \Ad^* \rho_j(1))^{\vee}.$$
\end{quote}
 
The proof  of the Lemma easily follows. We take $\mathcal{J}$ to be the subsequence of $n_j$s.  Let $k$ be as in the Lemma.
Take $n= n_j$, for any $j \geq k$,  and take $w^{(i)} \in \mathfrak{t}_{S_n}^{(i)}$ that have the same image in  $H^1_f(\Z[\frac{1}{S}], \Ad^* \rho_n(1))^{\vee}$.  
The images of $w^{(i)}$ in $\mathfrak{t}_{S_n}$ have the same image in $H^1_f(\Z[\frac{1}{S}], \Ad^* \rho_j(1))^{\vee}$, 
 (we are using the fact that the map  $H^1(\Z[\frac{1}{S}], \Ad \rho_n(1)) \rightarrow H^1(\Z[\frac{1}{S}], \Ad \rho_j (1))$ is surjective, by discussion  before \eqref{h00}) and therefore they act the same way on mod $p^k$ cohomology by  (*) above. 
 \qed

 \subsection{Conclusion of the proof of  Theorem \ref{maintheorem}. } \label{finalesec}  
 
 \label{DDR REF3} 
 Let us call a sequence of Taylor--Wiles data $Q_n$ of level $n$ (where we do not 
 require $n$ to vary through all the integers, but possibly some subsequence thereof)
 {\em $\Vinf$-convergent} if: \index{$\Vinf$-convergent} 
 \begin{itemize} 
 \item Each $Q_n$ is strict (Definition \ref{strictdef}) thus giving an action of $\Vinf$ on $H^*(Y(1), \Z/p^n)_{\mathfrak{m}}$. 
 \item The actions  converge to  an action of of $\Vinf$ on $H^*(Y(1), \Z_p)_{\mathfrak{m}}$.
In other words, if we fix $k$,  the action of $\Vinf$ on $H^*(Y(1), \Z/p^k)_{\mathfrak{m}}$ 
 arising from reducing the $Q_n$-action   is eventually constant.
 
 \end{itemize}

By Lemma \ref{Convergent implies strict} and passing to a further subsequence,
we see that $\Vinf$-convergent sequences exist. 
By Lemma \ref{EqualityOfTwo}, if $Q, Q'$ are two $\Vinf$-convergent sequences, 
  the resulting actions of $\Vinf$ on cohomology coincide.  
Thus at this point we have defined an action of $\Vinf$ on cohomology that is independent of choices, namely,
the action arising from any $\Vinf$-convergent sequence. 
This 
  action has the following property: 
 \begin{quote}  $(\dagger)$:  For any sequence $Q_n$ of Taylor--Wiles data,  
there is a subsequence $Q_{n_r}$ 
 such that  (for every $r$) the following two actions of $\mathfrak{t}_{S_{n_r}}$ on  $H^*(Y(1), \Z/p^r)_{\mathfrak{m}}$ coincide:
 
 \begin{itemize}
\item  The action via $\iota: \mathfrak{t}_{S_{n_r}} \rightarrow \mbox{ derived Hecke algebra with $\Z/p^r$ coefficients}$ (see \eqref{iotaQndef}). 
 \item The action obtained from the $\Vinf$-action, via   $f: \mathfrak{t}_{S_{n_r}} \twoheadrightarrow \Vinf/p^r$
(see   \eqref{tsnvsurj} ). 
\end{itemize} 
   \end{quote}  
  To see this, we first pass from $Q_n$ to a convergent subsequence 
  $(Q_{m_r}, r)$, where we regard $Q_{m_r}$ as having level $r$;
by Lemma \ref{Convergent implies strict} this means 
that  $Q_{m_r}$ is
  a strict datum of level $r$. We then pass to a further  subsequence $m_r'$
  to extract a $\Vinf$-convergent sequence; this gives the assertions above,
  but with $\Z/p^r$ and $\Vinf/p^r$ replaced by $\Z/p^{k(r)}, \Vinf/p^{k(r)}$
  where $k(r) \rightarrow \infty$ with $r$. Passing to a further subsequence gives the desired result.

 \proof (of Theorem \ref{maintheorem}, using $(\dagger)$): 
We have already constructed an action of $\Vinf$; let us prove, by contradiction, that it has property (*)
from the Theorem. 
Suppose that there is an integer $A$ and an infinite sequence of primes $q_n \equiv 1$ modulo $ p^n$
such that the pullback of the action via $ f_{q_n,A}: H^1(T_q, \Z/p^A) \rightarrow \Vinf/p^A$ 
fails to coincide with the action of $H^1(T_q, \Z/p^A)$ via the  embedding $\iota_{q_n, A}$ into the derived Hecke algebra with $\Z/p^A$ coefficients. 
We can choose a Taylor--Wiles system $Q_n$  containing $q_n$ and then get a contradiction to $(\dagger)$
as soon as $r > A$.  This proves (*). 
 
Now let us show that the image of $\wedge^* \Vinf$  in endomorphisms
of cohomology coincides with the global derived Hecke algebra.

Refer to the diagram \eqref{BigDog}, constructed with a   convergent sequence of Taylor--Wiles data $Q_n$. 
We will only use a subsequence of $n$s which is $\Vinf$-convergent. 
Consider for $n \geq k$ the map
\begin{equation} \label{Banjo} \Winf   = \mathfrak{t}_{S}/\mathfrak{t}_R  \stackrel{\eqref{feeling_groovy2}}{\longrightarrow} \mathfrak{t}_{S_n}/\mathfrak{t}_{R_n} 
\stackrel{\eqref{tsnvsurj}}{ \longrightarrow} \Vinf/p^k. \end{equation}
By \eqref{feeling_groovy2} and the discussion after \eqref{tsnvsurj}, 
the composite actually gives an isomorphism \begin{equation} \label{Banjo2} \Winf/p^k \simeq \Vinf/p^k.\end{equation}  
\label{DDR REF4} For fixed $k$ and large $n$, the map \eqref{Banjo} is independent of $n$: Choose $\tilde{v} \in \mathfrak{t}_S$. Let $v_n, v_m$
 be its image in $\mathfrak{t}_{S_n}, \mathfrak{t}_{S_m}$.
 As we saw in the diagram \eqref{BigDog}, the actions of $v_n, v_m$ on mod $p^k$ cohomology must coincide, 
 because both can be computed by means of the lift $\tilde{v}$. 
For this we are implicitly using  \eqref{logstructure} to see that  the composite map
{\small  \begin{equation}
 H^*(\Hom_{\Sinf}(\Cinf, \Z_p)) \rightarrow 
 H^*(\Hom_{S_n}(\Cinf \otimes_{\Sinf} S_n, \Z/p^n))   \stackrel{\sim}{ \rightarrow}
  H^*(Y(1), \Z/p^n)_{\mathfrak{m}}  
 \rightarrow   H^*(Y(1), \Z/p^k)_{\mathfrak{m}}
  \end{equation}}
is {\em independent of $n$}, for $n \geq k$.

 So the images of $v_n, v_m$ in $\Vinf/p^k$ have the same actions on mod $p^k$ cohomology. 
The $\Vinf/p^k$ action on mod $p^k$ cohomology is faithful (by \eqref{Banjo2} and \eqref{watney}) so this forces  the image of $v_n, v_m$
in $\Vinf/p^k$ to coincide as claimed.  
  
Therefore, passing to the limit over $n$, we get a map $\Winf/p^k \rightarrow \Vinf/p^k$, 
which is easily seen to be compatible as we increase $k$; thus the inverse limit over $k$
defines an isomorphism \begin{equation} \label{bo} \Winf  \stackrel{\sim}{\rightarrow} \Vinf.\end{equation}

Next, referring to \eqref{BigDog} the action of $\Ext_{\Sinf}^*(\Sinf/\Iinf, \Sinf/\Iinf) \simeq \wedge^* \mathfrak{t}_S$
on $H^*(\Hom_{\Sinf}(\Cinf,  \Z_p))$ certainly factors through $\wedge^*  \Winf $;
this action of $\wedge^* \Winf$ on $H^*(\Hom_{\Sinf}(\Cinf, \Z_p))$  is compatible
under the identifications  \eqref{logstructure} and  \eqref{bo} with the action of $\wedge^* \Vinf$ on  $H^*(Y(1), \Z_p)_{\mathfrak{m}}$. 

  Therefore, by \eqref{watney},  $\Vinf$ freely generates an exterior algebra inside $\mathrm{End}(H^*)$, and $H^*$ is freely generated over $\wedge^* \Vinf$
  by $H^{\BWq}(Y(1), \Z_p)_{\mathfrak{m}}$, its minimal degree component. 
 On the other hand, the image of $\wedge^* \Vinf$ in $\mathrm{End}(H^*)$ is contained in the global derived Hecke algebra --
 the action of an element $\Vinf$ is, by definition, a limit of actions of elements in the derived Hecke algebra.  Indeed,
 $\Vinf$ lies inside the strict global derived Hecke algebra. 
  \qed 

\subsection{Proof of Proposition \ref{P853}} \label{P853proof}

 {\color{\changecolor} . In our current situation, we have an inclusion
\begin{equation} \label{xyyzzww} \wedge^* \mathsf{V} \subset \gHecke\end{equation}
and we have seen that $H^*(Y(1), \Z_p)_{\mathfrak{m}}$
 is generated by its lowest nonvanishing degree $H^{\BWq}$ as a  
 $\wedge^* \mathsf{V}$-module.  
 If we know that $H^*(Y(1), \Z_p)_{\mathfrak{m}} = \Z_p$, 
 we argue that \eqref{xyyzzww} is an equality
 just as in the corresponding argument after  \eqref{inclusions BHH}.
 
 If we are given $\mathsf{S}$ as in the statement of the Proposition, 
 let $\mathsf{M} = H^{\BWq}(Y(1), \Z_p)_{\mathfrak{m}} \otimes_{\Z_p} \Q_p$
 as a $\mathsf{S}$-module. It is semisimple by assumption, and the natural map
 $$ \mathsf{M} \otimes \wedge^* \mathsf{V}_{\Q_p}      \rightarrow H^*(Y(1), \Z_p)_{\mathfrak{m}} \otimes \Q_p$$
 is now an isomorphism of $\mathsf{S} \otimes \wedge^* \mathsf{V}_{\Q_p}$ modules. 
 However,
$$\mbox{(graded) commutant of $\mathsf{S} \otimes \wedge^* \mathsf{V}_{\Q_p}$   on $\mathsf{M} \otimes \wedge^* \mathsf{V}_{\Q_p}$} =   \mathsf{S}' \otimes  \wedge^* \mathsf{V}_{\Q_p},$$
as can be verified by computing in steps:  the commutant of $\mathsf{S}$ alone
 equals $  \mathsf{S}' \otimes \mathrm{End}(\wedge^* \mathsf{V}_{\Q_p})$,
 and then inside here the (graded) commutant of $\wedge^* \mathsf{V}_{\Q_p}$   is $\mathsf{S}' \otimes \wedge^* \mathsf{V}_{\Q_p} $.

Since this commutant contains  the image of $\gHecke$
acting on $\Q_p$-cohomology, the latter is also graded commutative, since $\mathsf{S}'$ is commutative by the multiplicity one hypothesis.   Note that the action on $\Q_p$ cohomology
is faithful because of our torsion freeness assumption, see 7(a) of  \S \ref{Sec:assumptions}.
%
%
  }
\qed

 \subsection{The action of Hecke operators}  \label{Galoisindexingproofs}
 To conclude, let us translate what we have proved into a more concrete   assertion about the action of a derived Hecke operator. 
 
 Let $q$ be a unramified prime for $\rho$,  with $q \equiv 1$ modulo $p^n$, and such that
 $\rhobar(\Frob_q)$ is conjugate to a strongly regular element of $T^{\vee}(k)$. Let $\nu \in X_*(\mathbf{A})^+$ be strictly dominant 
and let
  $$\alpha \in H^1(\mathbf{A}(\F_q), \Z/p^n).$$ 
  To this we can associate in a natural way (see below) a derived Hecke operator $T_{q, \nu, \alpha}$
  as well as an element $[q, \nu, \alpha] \in \Vinf/p^n$; we will prove
  that the actions of these are compatible (see Lemma below), 
 justifying the assertions made in \S \ref{Galoisindexing}.

First of all,
a small piece of linear algebra. Let $k$ be a field.  Suppose given a  fixed character $\psi \in X^*(T^{\vee})$. 
Let $g \in G^{\vee}(k)$ be regular semisimple, with centralizer $Z_g$; 
this data allows us to construct a homomorphism of $k$-vector spaces
$$e_{\psi,g}: \mathrm{Lie}(T^{\vee}) \rightarrow \Lie(Z_g),$$ 
$$ e_{\psi,g}:  \sum_{\phi: T^{\vee} \stackrel{\sim}{\rightarrow} Z_g} \langle \psi, \phi^{-1}(g) \rangle \cdot d\phi$$
where the sum is taken over all conjugations of $T^{\vee}$ to $Z_g$ over $\bar{k}$;
the morphism is nonetheless defined over $k$. 

{\em Example:}   if $G^{\vee}=\SL_2$, take $T^{\vee}, B^{\vee}$ in the standard way to be the diagonal subgroup and upper triangular matrices, and take $\psi: \left( \begin{array}{cc} x & 0 \\ 0 & x^{-1} \end{array} \right) \mapsto x$.  Then $e_{\psi, g}$ sends $  \left( \begin{array}{cc} 1 & 0 \\ 0 & -1 \end{array} \right) \in \Lie(T^{\vee})$ 
to the element  $  2g - \mathrm{trace}(g) \in M_2(k)$: it's enough to check this for $g \in T^{\vee}$, where the result is clear.

\medskip

Let $q, \nu, \alpha$ be as described at the start of this subsection. We can then construct a class 
  $$ [q, \nu, \alpha] \in \Vinf/p^n$$
  in the following way:  regarding $\nu$ as a character of
  $T^{\vee}$,  and use the linear algebra construction mentioned with $k = \Q_p$ to make the first map of   \begin{equation} \label{cheeseboard} \mathrm{Lie}(T^{\vee}) \stackrel{e_{\nu, \Frob_q}}{\longrightarrow } \Lie(Z_{\rho(\Frob_q)}) \hookrightarrow \Lie(G^{\vee})\end{equation}
  (at first we get this $\otimes \Q_p$ but then it preserves the integral structures, with reference to the natural $\Z_p$-models of the three groups above). 
  The resulting embedding $\mathrm{Lie}(T^{\vee}) \rightarrow  \Ad \rho$
    is a morphism of $\Gal(\overline{\Q_q}/\Q_q)$-modules,
  where $\mathrm{Lie}(T^{\vee})$ is taken to have the trivial action.

As before, we may identify
  $$H^1(\mathbf{A}(\F_q), \Z/p^n) = H^1(X_*(\mathbf{A}) \otimes \F_q^{\times} , \Z/p^n) = \Hom(\F_q^{\times}, X_*(T^{\vee}) /p^n)$$
and so from $\alpha$ we obtain a class  
$$\alpha' \in   \frac{ H^1(\Q_q, \Lie(T^{\vee})/p^n)}{H^1_{\ur}(\Q_q, \Lie(T^{\vee})/p^n)}.$$
   Here $\Lie(T^{\vee})$ is taken as a trivial Galois module.
We can then form 
$$\mbox{pushforward of } \alpha' \mbox{via \eqref{cheeseboard}}  \in   \frac{ H^1(\Q_q, \Ad \rho_n)}{H^1_{\ur}(\Q_q, \Ad \rho_n)}. $$ 
  and, as usual, this can be paired with  $H^1_f(\Z[\frac{1}{S}], \Ad^* \rho (1))$ by means of local reciprocity.
  In this way we obtain a functional $H^1_f(\Z[\frac{1}{S}], \Ad^* \rho (1)) \rightarrow \Z/p^n$, which we denote as
$$[q, \nu, \alpha] \in \Vinf/p^n.$$ 
  
  \begin{lemma} \label{reciprocity law}
Let $q, \nu, \alpha$ be as above.   Let $[q, \nu, \alpha] \in \Vinf/p^n$ be as defined above.

Let $T_{q, \nu,\alpha}$ be the derived Hecke operator with $\Z/p^n$ coefficients which is 
supported on the $G_q$-orbit of $(\nu K_q, K_q)$ and 
whose value at $(\nu K_q, K_q)$
which corresponds to $\alpha$ under the   cohomology isomorphism $H^*(K_q \cap \nu K_q \nu^{-1}, \Z/p^n) \cong H^*(\mathbf{A}(\F_q), \Z/p^n)$. 

Then $T_{q,\nu,\alpha}$ corresponds to $[q, \nu,\alpha] \in \Vinf/p^n$, in the following asymptotic sense:

 There is $N_0(m)$ such that  for  $q ,\nu, \alpha$ as above with 
 $q \equiv 1$ modulo $p^{N_0(m)}$, 
 the actions of  $T_{q, \nu, \alpha}$ and $[q, \nu,\alpha]$
 on $H^*(Y(1), \Z/p^m)_{\mathfrak{m}}$ coincide. 
\end{lemma}
  \proof
Under the derived Satake isomorphism, $T_{q,\nu, \alpha}$ is sent to  
$$\sum_{w} w \nu\otimes w \alpha
\in \left(S[X_*] \otimes H^1(\mathbf{A}(\F_q))\right)^W.$$   

With notation as in \S \ref{Vnaction}, 
let $\Theta_{\nu} \in S[X_*]_{\mathfrak{m}}^W$ be  defined
so that its image under $S[X_*]_{\mathfrak{m}}^W \hookrightarrow S[X_*]_{\tilde{\mathfrak{m}}} \rightarrow S[X_*]_{\chi}$
is equal to $\nu$. Here, we regard $\nu \in X_* \hookrightarrow S[X_*]$. 

 
 Then, after completing at $\mathfrak{m}$, we have an equality
  $$\mathrm{Satake}(T_{q,\nu,\alpha}) =   \sum_{w \in W}  \Theta_{w \nu} \cdot  \mbox{Satake}(\theta(w \alpha)) \in (S[X_*] \otimes H^1(\mathbf{A}(\F_q)))^W_{\mathfrak{m}}$$
 where $\theta$ is as in \eqref{embedH1T}.  (We can check this using the isomorphism  \eqref{composite iso}:
 it gives an isomorphism of the target group with  $S[X_*]_{\chi} \otimes H^1(T_q)$,
 under which $\theta(w \alpha)$ is, by its very definition,  sent to $1 \otimes w \alpha$; under the same isomorphism $\Theta_{w \nu}$ 
 is sent to $w \nu \otimes 1$, and the result follows.) 
 

As before we have fixed $\Frob^T_q \in T^{\vee}(\F_p)$ an element conjugate to $\rhobar(\Frob_q)$;
fix a lift $t_q \in T^{\vee}(\Z_p)$ that is conjugate to $\rho(\Frob_q)$. 
Then  $ \Theta_{w \nu}$  (more exactly, its preimage under Satake) acts on $H^*(Y(1), \Z_p)_{\mathfrak{m}}$
by  $\langle w \nu, t_q \rangle$ (this makes sense: $w \nu \in X^*(T^{\vee})$ and $t_q \in T^{\vee}(\Z_p)$, so they can be paired to get an element of $\Z_p^{\times}$).   Using Theorem \ref{maintheorem}, we see that the action of $T_{q,\nu, \alpha}$ on $H^*(Y(1), \Z_p)_{\mathfrak{m}}$
corresponds (in the sense  of the lemma statement) to the element
$$ \sum_{w \in W} \langle w \nu, t_q \rangle \cdot f_{q,n}(w \alpha)  =  \sum_{w \in W} \langle \nu, w^{-1} t_q \rangle \ \  f_{q,n}(w \alpha) \in  \Vinf/p^n.$$
Winding through the definitions, this element of $\Vinf/p^n$ is  exactly $[q, \nu, \alpha]$. \qed

 \section{Some very poor evidence for the main conjecture: Tori and the trivial representation} \label{pisspoor}
 
 We verify that the main conjecture  (Conjecture \ref{mainconjecture})
  holds in the case when $\mathbf{G}$ is an anisotropic torus.  This is straightforward,  but still slightly comforting. 
  
  One may also verify that  a certain analogous statement to Conjecture \ref{mainconjecture}
  holds in the situation studied in  \S \ref{Quillen},
but there we do not understand the situation clearly at present --  hopefully it will eventually prove to be a specialization
 of the general conjecture to the nontempered case.
    
  \subsection{Setup} \label{tori}

Let $\mathbf{T}$ be an anisotropic $F$-torus; let $\mathcal{O}$ be the ring of integers of $F$. 
Let us fix a finite set of places $\ramprimes$ such that $\mathbf{T}$ admits a smooth model over $\mathcal{O}[\frac{1}{\ramprimes}]$. 
We assume it contains all places $\wp$ above the rational prime $p$.

  The associated symmetric space 
$$ \mathcal{S} =  \mathbf{T}(F \otimes \R)/\mbox{maximal compact}$$
has  $q,\delta$ invariants (see \eqref{qdef}):
$$ q= 0 , \delta = \dim(\mathcal{S}).$$ 
The arithmetic manifold $Y(K)$ associated to a level structure $K$  is
a disjoint union of copies of $\mathcal{S}/\Delta$, where 
  $\Delta$ is a congruence subgroup of $\mathbf{T}(\mathcal{O})$. 
Moreover  the quotient $ \mathcal{S}/\Delta$ is a union of compact tori, and thus 
the rank of $\Delta$, i.e. $\dim_{\Q} (\Delta \otimes \Q)$, equals $\delta$.  
  We will suppose that $K$ is chosen so small that $\Delta$ is free of $p$-torsion.

\subsection{Galois cohomology}

Let $M$ be the motive associated to the first homology group $H_1(\mathbf{T})$ of $\mathbf{T}$. 

Let $X_*(\mathbf{T})$ be the co-character group of $\mathbf{T}$. It carries an action of $\Gal(\overline{F}/F)$. 
Coming from $X_* \otimes_{\Z} \mathbb{G}_m \stackrel{\sim}{\rightarrow} \mathbf{T}$, we get  an isomorphism of $\Gal(\overline{F}/F)$-modules 
$$M_p := \mbox{$p$-adic realization of $M$} = X_*(\mathbf{T}) \otimes_{\Z} \Z_{p}(1) \simeq \varprojlim \mathbf{T}[p^n],$$
the $p$-adic Tate module of $\mathbf{T}$. 
Computing with the Kummer sequence,
\begin{equation} \label{cococo} H^1(\mathcal{O}[\frac{1}{\ramprimes}] ,  M_p) \simeq \mathbf{T}(\mathcal{O}[\frac{1}{\ramprimes}]) \otimes \Z_p.  \end{equation}
Inside the left-hand side, we can consider those classes  that are crystalline at $p$ and 
 unramified at {\em all primes away from $p$, including primes in $S$.} We will denote this subgroup by $H^1_{f, \mathrm{ur}}(F, M_p)$. \footnote{For a precise
 definition we refer to \cite{BK}: the group we have called $H^1_{f, \mathrm{ur}}$  is the group denoted by $H^1_{f, U}$ in \cite[Definition 5.1]{BK} with $U$ taken to be all finite places. This group is often 
simply written $H^1_f(F, -)$ in the literature.}

 The subgroup on the right-hand side corresponding to $H^1_{f, \mathrm{ur}}(F, M_p)$, call it
$$ \Delta' \subset  \mathbf{T}(\mathcal{O}[\frac{1}{\ramprimes}]) \otimes \Z_p,$$ 
 is  commensurable to the image of  $\Delta \otimes \Z_p$ in the right-hand side of \eqref{cococo}.
{\color{\changecolor} Indeed, the constraint that a point $t \in  \mathbf{T}(\mathcal{O}[\frac{1}{\ramprimes}])$
have cohomology class in $H^1_{f, \mathrm{ur}}$ is equivalent
to forcing $t$ to belong to a suitable open compact subgroup of $\mathbf{T}(F_v)$
for each $v \in S$; see   \cite[Theorem 2.3.1]{HK} for this statement in the trickier case when $v$ is above $p$.} 
We suppose (shrinking $\Delta$ a little if necessary) that $\Delta \otimes \Z_p \subset \Delta'$.
  
Also, the {\em motivic} cohomology $H^1_{\mot}(M, \Q(1))$
is identified with  with the $\mathbf{T}(\Q) \otimes \Q$,  as we may see by first passing to an extension that trivializes the torus $\mathbf{T}$. Presumably
the following is valid, but I did not try to check it:

\begin{quote} {\em Assumption:} The subgroup of ``integral classes'' $H^1_{\mot}(M_{\Z}, \Q(1))$ (see discussion after \eqref{saxophone}) is identified  with the image of $\Delta \otimes \Q$ inside $\mathbf{T}(\Q) \otimes \Q$. 
\end{quote}

Now a cohomological automorphic form $\Pi$ for $\mathbf{T}$ is trivial on the connected component of $\mathbf{T}(F \otimes \R)$, i.e.
they are the {\em finite order}  id{\`e}le class characters of $\mathbf{T}$.   
However, the associated co-adjoint motive (see \S \ref{sec:MV}) doesn't depend on which idele class character: we have simply 
$$ \left(\mbox{ coadjoint motive for $\Pi$}  \right) (1) \simeq M,$$
the motive $M$ described above. \footnote{
 Indeed, the the $\Z_{p}$-linear dual of $\hat{T}$ is identified with 
$$ \Lie(T^{\vee})^{\vee} \simeq \left( \Lie(\mathbf{G}_m) \otimes X_*(T^{\vee}) \right)^{\vee} \simeq X^*(T^{\vee})  \otimes \Z_p \simeq X_*(\mathbf{T}) \otimes \Z_{p}$$
where we have fixed an isomorphism $\Z \simeq \Lie(\mathbf{G}_m)$.  }

Now let us examine Conjecture \ref{mainconjecture} in this case. Put 
$$\mathsf{V} = H^1_{f, \mathrm{ur}}(F,  M_p)^{\vee},$$
where $M_p$ is the $p$-adic etale realization, and $\vee$ denotes $\Z_p$-linear dual; thus $\mathsf{V} = \Hom(\Delta',\Z_p)$. 
The notation $H^1_{f, \mathrm{ur}}$ was defined after \eqref{cococo}.

There is a natural action of $\mathsf{V}$ on $H^*(Y, \Z_p)$,
obtained from the maps
\begin{equation} \label{frow} \mathsf{V} = H^1(\Delta', \Z_p) \rightarrow H^1(\Delta, \Z_p). \end{equation}
Moreover, motivic cohomology gives a lattice in $\mathsf{V} \otimes \Q = H^1(\Delta', \Q_p)$
(the classes which are $\Q$-valued on $\Delta \subset \Delta'$)
and obviously this lattice indeed preserves $H^*(Y, \Q)$, in the
$\Q$-linear extension of the action of $\mathsf{V}$.

The only point to be discussed is that the action \eqref{frow} is indeed that 
resulting from the  same formalism as \S \ref{reciprocity}.  We describe this only briefly. 
Let $v$ be a finite place not belonging to $S$, so that $\mathbf{T}$ extends to a smooth torus over $\mathcal{O}_v$. 
As usual we have an injection, 
 $$  \Hom(\mathbf{T}(\mathcal{O}_v), \Z/p^n) \hookrightarrow \mbox{ derived Hecke algebra at $v$,}$$
 and thus an action of  the left-hand group on the cohomology
 of $H^*(Y,\Z/p^n)$; explicitly, this action is obtained by pulling back cohomology classes via $\Delta \rightarrow \mathbf{T}(\mathcal{O}_v)$, and cup product. 
 By just the same procedure as  that described  in \eqref{fqndef}, we can construct a map  
 $$ \Hom(\mathbf{T}(\mathcal{O}_v), \Z/p^n)  \rightarrow \underbrace{ \mathsf{V}/p^n}_{\simeq \Hom(\Delta',\Z/p^n)} ,$$
and one verifies this is the  map induced by $\Delta' \rightarrow \mathbf{T}(\mathcal{O}_v) \otimes \Z_p$. 
 Then the action of $\mathsf{V}/p^n$  on $H^*(Y, \Z/p^n)$ is compatible with the ``derived Hecke'' action of $\Hom(\mathbf{T}(\mathcal{O}_v), \Z/p^n) $ for all $v$;
 and in fact this compatibility determines the action of $\mathsf{V}/p^n$.

  \appendix

\section{Remedial algebra} \label{remedial}

In this section we fill in 
some ``intuitively obvious'' claims in the text in grotesque detail, in particular
the identifications between  various different models of the derived Hecke algebra. 
(The  word ``remedial'' in the title of this section refers to my own lack of fluency with 
  homological algebra.)
 
 An action of a topological group will be called {\em smooth} if the stabilizer of every point is an open subgroup.      We will fix a 
   finite ring $S$ of cardinality prime to $p$.  A ``smooth representation of $G$'', in this section,
   will be a smooth action of $G$ on an $S$-module.  We write $SG$ (or occasionally $S[G]$ when typography requires) for the group algebra of $G$ with 
   coefficients in $S$. {\color{\changecolor} We will write $\Hom_{SG}$ for homomorphisms of $SG$-modules; in \S \ref{sec:mult} we abbreviate this to $\Hom_{G}$ because
   other notation becomes very dense.}  
      
   Note that the usage of $U$ and $K$ in this section do not precisely match their usage in the main text. 

\subsection{} \label{Udefsec}
Let $K$ be a profinite group, which admits a pro-$p$ open normal compact subgroup $U$. 
Then the category $\mathcal{C}$ of smooth representations of $K$ 
 is   an abelian category with enough projectives.  
 In fact, 
 if $Q$ is a projective $K/U$ module, then considering $Q$
 as a smooth $K$-module $\tilde{Q}$ it remains projective: 
 $\Hom_{SK}(\tilde{Q}, V) = \Hom_{S[K/U]}(Q, V^U)$ and $V \mapsto V^U$ is exact by the hypothesis on $U$. 
 (One can lift $U$-invariants under a surjection $V_1 \twoheadrightarrow V_2$ by averaging.) 

In this situation, restriction to a finite index subgroup $K' \subset K$ preserves projectivity.  
Indeed $\Hom_{SK'}(\mathrm{Res}^{K}_{K'} A, B) = \Hom_{SK}(A, \mathrm{Ind}_{K'}^K B)$
and the induction functor is exact.

\subsection{} \label{enough_proj}
Now let $G=G_v$ be the points of a reductive group over a $p$-adic field, or any open subgroup thereof. Then 
 the category of  smooth representations of $G_v$ is an abelian category and it has enough projectives.

Indeed,
consider $W = S[G/U]$ for an open pro-$p$ compact $U \subset G$.
Then $$\Hom_{SG}(W, V) \simeq V^{U},$$ which is obviously exact in $V$. 
   So $W$ is projective, and now given any other $V$
we choose generators $v_i$ for $V$, open pro-$p$ compact subgroups $U_i$ fixing $v_i$,  corresponding projectives
$W_i$, and then get $\bigoplus W_i \twoheadrightarrow V$.

 Throughout the remainder of this section, we suppose that $G$ is  as above,
that  $K$ is an open compact subgroup of $G$ (in particular, $K$ is profinite),
and that $U \subset K$ is a pro-$p$ open {\color{\changecolor} normal} compact subgroup of $K$.

\subsection{} \label{Qdef}  
 
 Fix 
a resolution of the trivial $K$-representation by projective smooth $K$-modules: 
\begin{equation} \label{Qdef-eqn} \mathbf{Q}: \cdots \rightarrow Q_i \rightarrow \cdots \rightarrow Q_1 \rightarrow S.\end{equation}
{\color{\changecolor} (Here, and in what follows, we will use boldface letters to denote complexes.) }  To be explicit, let us take $\mathbf{Q}$ to be the standard ``bar'' resolution of $S$ by free $S[K/U]$-modules,
considered as a complex of smooth $K$-representations. 

Then $\Hom(\mathbf{Q}, \mathbf{Q})$  computes $H^*(K, S)$,
the continuous cohomology of the profinite group $K$ with $S$ coefficients: indeed, the   the cohomology of
$\Hom(\mathbf{Q}, \mathbf{Q})$  is identified with $H^*(K/U, S)$,  which is identified by pullback with the continuous cohomology $H^*(K, S)$.

The  complex $\Hom(\mathbf{Q}, \mathbf{Q})$ has the structure of differential graded algebra arising from composition, and    the resulting multiplication on $H^*(K, S)$ 
coincides with the cup product (this reduces to a corresponding statement for a finite group; for that see \cite{Yoneda}). 
 
If $K' \subset K$ is a finite index subgroup, then $\Hom_{SK'} (\mathbf{Q}, \mathbf{Q})$
still computes $H^*(K', S)$ (see remarks above).  Moreover, the averaging operator $\sum_{K/K'}$
-- that is to say, the map sending  
$$ f \in  \Hom_{SK'}(\mathbf{Q}, \mathbf{Q}) \mapsto  \sum_{\kappa} \kappa f \kappa^{-1} \in \Hom_{SK}(\dots) $$ 
realizes the corestriction map $H^*(K', S) \rightarrow H^*(K, S)$, where the $\kappa$ sum is taken over a set of coset representatives for $K/K'$ in $K$.

  \subsection{Induction and Frobenius reciprocity } 
We use the word ``induction'' for the functor from $K$-modules to $G$ modules
 \begin{equation} \label{tensor} M \rightsquigarrow  SG \otimes_{ SK} M.\end{equation}

This is  isomorphic to the usual ``compact''  induction, namely
space of functions
\begin{eqnarray} \label{functionmodel} \ind_K^G(M)  &:=&   \{ f: G \rightarrow M: f(g k) = k^{-1} f(g), 
\\  \nonumber  &&  \mbox{ $f$ is supported on finitely many left translates of $K$} \} \end{eqnarray}
 where the action of $h \in G$ is by left translation, i.e. $l_h f(x) = f(h^{-1} x)$. 
 We will drop the word ``compact'' and simply refer to \eqref{functionmodel}
 or \eqref{tensor} as ``induction''; we refer to the 
model  \eqref{functionmodel} for induction as the ``function model.'' 

We can define inverse isomorphisms between  \eqref{functionmodel} and \eqref{tensor} as follows:  
  define $\ind_K^G(M) \rightarrow S[G]\otimes_{S[K]} M$ via 
 $$ f \mapsto \sum_{x \in (G/K)} g_x \otimes f(g_x), $$ 
 where $g_x \in G$ is a representative for $x \in G/K$;
in the other direction, we send $g \otimes m$ to the function supported on $gK$
whose value on $g$ equals $m$

\subsection{Frobenius reciprocity}

We have Frobenius reciprocity  $$\Hom_{SG}(\mathrm{ind}^G_K \ M, N) \simeq \Hom_{SK}(M, N)$$
and therefore induction carries projective $K$-modules to projective $G$-modules. 
Explicitly  an $SK$-homomorphism $f: M \rightarrow N$
is sent to its obvious $G$-linear extension $SG \otimes_{SK} M \rightarrow N$.

{\em If $G \supset K$ has finite index,} we also have the reverse adjointness (since ``compact induction'' and ``induction'' coincide): 
 to give a $K$-map
$f:  M \rightarrow M'$ is the same as giving a map $F_f: M \rightarrow \ind_{K}^G M'$.
Explicitly, in the function model for the induced representation, $F_f$ is characterized by the property 
$F_f(m)(e) = f(m)$, and so
$$F_f(m)(g) = l_{g^{-1}} F_f(m)(e) = F_f(g^{-1} m)(e) = f(g^{-1} m)$$
and thus in the tensor product model 
\begin{equation} \label{tpmodel} F_f(m) = \sum_{x \in (G/K)} g_x \otimes   f(g_x^{-1} m)\end{equation}

\subsection{Restriction of induced representations}

Let $Q$ be a smooth representation of $K$. 
The restriction of $\ind^G_K Q$ to $K $ is isomorphic 
to \begin{equation} \label{ind-res} \bigoplus_{x \in  K \backslash G/K}  S[K g_x K] \otimes_{SK} Q \simeq \bigoplus_{K} S[K] \otimes_{S K_x}  Q_x  \simeq \bigoplus_{x} \ind_{K_x}^K Q_x \end{equation}  
where $x = g_x K$ runs through a set of representatives for $K$-orbits on $G/K$, and we write $K_x = K \cap g_x K g_x^{-1}$; moreover
for  a $K$-module $Q$, 
we denote by $Q_x$ the $K_x$-module whose underlying space is $Q$, but 
  for which the action $*$ of $K_x$ on $Q$ is defined thus:
\begin{equation} \label{twist} \kappa * q =  \left( \Ad(g_x^{-1}) \kappa \right)  q.\end{equation}

The first map of \eqref{ind-res} is given  explicitly by 
\begin{equation} \label{IR1} k_1 g_x k_2 \otimes q  = k_1 g_x \otimes k_2  q \mapsto k_1 \otimes k_2  q\end{equation} 
and
the inverse map sends \begin{equation} \label{inverse-map} k \otimes q \mapsto k g_x \otimes q.\end{equation}

In the function model of the induced representations, the composite map of \eqref{ind-res} 
sends  $F: K g_x K \rightarrow Q$ to the function $F': k \mapsto F(k g_x)$.   
In the reverse direction, given a function $F'$ in the function model of $\ind_{K_x}^K Q_x $,
the inverse of \eqref{ind-res} sends it to 
\begin{equation} \label{yuk-yuk}  \sum_{k \in K/K_x}  k g_x \otimes f(k)  \in S[K g_x K ] \otimes_{SK} Q.\end{equation}

\subsection{Derived Hecke algebra}  \label{Dha-def}
The {\em derived Hecke algebra} for the pair $(G,K)$ with coefficients in $S$  is defined as 
$$ \bigoplus_i \Ext^i_{SG}(S[G/K], S[G/K])$$
where the $\Ext$-groups are taken in the category of smooth $S$-representations. 

We can construct an explicit model as follows. Let $\mathbf{Q}$ be as in \eqref{Qdef-eqn}. Then 
  $\mathbf{P} = \ind \mathbf{Q} $ is a  projective resolution of $S[G/K]$. 
In particular,  $\Hom_{SG}(\mathbf{P}, \mathbf{P})$
has the structure of a differential graded algebra and its cohomology gives the derived Hecke algebra.

\subsection{} 
We will now explicitly describe the isomorphism \eqref{TTspin} between the derived Hecke algebra  and  its ``double coset model.''

Let $\mathbf{P}, \mathbf{Q}$ be as in \S \ref{Dha-def}. 

We have \begin{equation} \label{dha-additive} \Hom_{SG}(\mathbf{P},\mathbf{P}) = \Hom_{SK}(   \mathbf{Q}, \ind^G_K  \mathbf{Q}) \stackrel{
\eqref{ind-res}}{
\longleftarrow} \underbrace{ \bigoplus_{x \in K \backslash G/K}   \Hom_{SK}(\mathbf{Q}, \ind_{K_x}^K \mathbf{Q}_x) }_{=
\Hom_{SK_x}(\mathbf{Q}, \mathbf{Q}_{x})},\end{equation} 
where   $x$ varies  now through $K \backslash G/K$, again $g_x K$ is a representative for $x$
and $K_x = K \cap \Ad(g_x) K$, and   the twist operation $Q_x$ is as described in \eqref{twist}. 
 
 Note that the last map induces a cohomology isomorphism.   We must
 see that $H^*(\Hom_{SK}(\mathbf{Q}, -))$  commutes with the infinite direct sum $\bigoplus_x \ind_{K_x}^K \mathbf{Q}_x$. 
 However,  
$\mathbf{Q}_x$ is  cohomologically concentrated in degree $0$, and so the same is true for $\ind_{K_x}^K \mathbf{Q}_x$;
it is enough to show, then, for any $K$-modules $M_i$,  the obvious map
$$ \ \bigoplus_{i} \Hom_{SK}(\mathbf{Q}, M_i) \longrightarrow \Hom_{SK}(\mathbf{Q}, \bigoplus_i M_i))$$
is  a quasi-isomorphism.  But this follows from the fact that  taking $U$ invariants commutes with infinite direct sum,
as does    the functor
 $H^*(K/U, -)$.

Note that the cohomology of 
$ \Hom_{SK_x} (\mathbf{Q}, \mathbf{Q}_x) $  is identified with $H^*(K_x, S)$, because 
$\mathbf{Q}$ and $\mathbf{Q}_x$ are resolutions of $S$, and moreover
$\mathbf{Q}$ is a complex of projective $K_x$-modules. Thus, \eqref{dha-additive}
gives rise to an isomorphism: 
\begin{equation} \label{serre} H^* \left( \Hom_{SG}(\mathbf{P},\mathbf{P})  \right) \simeq \bigoplus_{x} H^*(K_x, S).\end{equation}

For later use, we explicate the map of \eqref{dha-additive}, going from right to left: 
An element   $f \in \Hom_{SK_x}(\mathbf{Q}, \mathbf{Q}_{x})$ must satisfy
$f( \kappa q) = (g_x^{-1} \kappa g_x) f(q)$
for $\kappa \in K_x$; the  associated element of $\Hom_{SK}(\mathbf{Q}, \ind^G_K \mathbf{Q})$
is given in the tensor product model  of the induced representation by  the formula of \eqref{tpmodel}: 
$$  q \in  \mathbf{Q} \mapsto \sum_{k \in K/K_x}  k g_x \otimes f(k^{-1} q)$$
which is well-defined.

  \subsection{Action of derived Hecke algebra on derived invariants}
  Now suppose that $\mathbf{M}$ is a {\em complex} of smooth $G$-representations.     There is a natural action of $\End_{SG}(\mathbf{P})$  
 on $\Hom_{SG}(\mathbf{P}, \mathbf{M})$. Moreover, the latter complex computes the hypercohomology \index{$\mathbb{H}$ (hypercohomology)} $\mathbb{H}^*(K, \mathbf{M})$
 of $K$ with coefficients in the complex $\mathbf{M}$. 
 
 Thus, because of \eqref{serre}, we get an action of $H^*(K_x, S) $ on $\mathbb{H}^*(K, \mathbf{M})$. Let us 
describe the action of $h_x \in  H^*(K_x,S)$  on $\mathbb{H}^*(K, \mathbf{M})$ as explicitly as possible,  in particular justifying the claims
of \S \ref{sec:concrete}: 

\begin{lemma} 
With notation as above, the  action of $h_x$
is given by the following composite:
 $$ \mathbb{H}^*(K, \mathbf{M}) \stackrel{\Ad(g_x^{-1})^*}{\longrightarrow} \mathbb{H}^*(K_x, \mathbf{M}_x) \stackrel{m \mapsto g_x m}\longrightarrow \mathbb{H}^*(K_x, \mathbf{M}) \stackrel{\cup h_x}{\longrightarrow} \mathbb{H}^*(K_x, \mathbf{M})
\stackrel{\mathrm{Cores}}{\longrightarrow} \mathbb{H}^*(K, \mathbf{M}).$$
 Here the first map is the pull-back induced by $\Ad(g_x^{-1}): K_x \hookrightarrow K$, which pulls back $\mathbf{M}$ to $\mathbf{M}_x$. 
\end{lemma}

 \proof Choose $h_x'  \in \Hom_{SK}(\mathbf{Q},  \ind_{K_x}^K \mathbf{Q}_x) $ representing $h_x$. For $f \in  \Hom_{SG}(\mathbf{P}, \mathbf{M})$
 we denote by $f_x \in \Hom_{SK}(\ind_{K_x}^{K} \mathbf{Q}_x,\mathbf{M}) $ the restriction. 
 We denote by $[h'_x] \in  \Hom_{SK_x}(\mathbf{Q}, \mathbf{Q}_x) $
 and $[f_x]    \in \Hom_{SK_x}(\mathbf{Q}_x, \mathbf{M}) $ the elements obtained from $h'_x, f_x$ using Frobenius reciprocity (but using the two different versions of Frobenius reciprocity).

  We want to compare the composition  $f_x \circ h_x'$ and $[f_x] \circ [h_x']$
i.e.  
{\small \begin{equation}
 \xymatrix{
h_x'\in \Hom_{SK}(\mathbf{Q},  \ind_{K_x}^K \mathbf{Q}_x)    &  \times &  f_x \in \Hom_{SK}( \ind_{K_x}^K \mathbf{Q}_x, \mathbf{M})  \ar[r] & \Hom_{SK}(\mathbf{Q}, \mathbf{M}) \\ 
[h_x'] \in  \Hom_{SK_x}(\mathbf{Q}, \mathbf{Q}_x) \ar[u]^{\sim}  &  \times & [f_x] \in  \Hom_{SK_x}(\mathbf{Q}_x, \mathbf{M}) \ar[r]  \ar[u]^{\sim} & \Hom_{SK_x}(\mathbf{Q}, \mathbf{M})  
 }
\end{equation} }
We compute
 $$f_x \circ h_x'(q) \stackrel{\eqref{tpmodel}}{=} f_x\left( \sum_{k_i \in K/K_x} k_i\otimes   [h_x'](k_i^{-1} q)  \right) = \sum_{k_i} k_i [f_x] \ \circ   [h_x'](k_i^{-1} q) $$
i.e., this is what we get by averaging $ [f_x] \circ [h_x']$ over the action of $K/K_x$.
The  cohomology class of the composition $[f_x] \circ [h_x']$ simply amounts (at the level of cohomology)
 to the cup product of the class $[f_x] \in \mathbb{H}^*(K_x, \mathbf{M})$ (hypercohomology) with the class
 $[h_x'] \in H^*(K_x, S)$.   So to prove the lemma it remains to show:
 
 \medskip
 
 {\em Subclaim:} The class  $[f_x] \in \mathbb{H}^*(K_x, \mathbf{M})$
 is obtained from $[f] \in \mathbb{H}^*(K, \mathbf{M})$ 
via 
  $$ \mathbb{H}^*(K, \mathbf{M}) \stackrel{\Ad(g_x^{-1})^*}{\longrightarrow} \mathbb{H}^*(K_x, \mathbf{M}_x) \stackrel{m \mapsto g_x m}\longrightarrow \mathbb{H}^*(K_x, \mathbf{M}) $$
  
  \medskip

 At the chain level this map is given by the composite
 $$ \Hom_{SK}(\mathbf{Q}, \mathbf{M}) 
 \rightarrow \Hom_{S K_x}(\mathbf{Q}_x, \mathbf{M}_x) \rightarrow  \Hom_{SK_x}(\mathbf{Q}_x, \mathbf{M})$$
 where the first map is the trivial map, just considering maps of $K$ modules
 as $K_x$-modules via $\Ad(g_x^{-1}):K_x \rightarrow K$; and the second map $\mathbf{M}_x \rightarrow \mathbf{M}$ is given by $m \mapsto g_x m$.  
To check the subclaim, note  that  $[ f_x] \in \Hom_{S K_x}(\mathbf{Q}_x, \mathbf{M})$ 
 sends $q \in \mathbf{Q}_x$ to the value
 of $f_x$ on the element
 $1 \otimes q  \in \ind^{K}_{K_x} \mathbf{Q} $, which is carried by the isomorphism  inverse to \eqref{ind-res}
 to $g_x \otimes q  \in \ind_{K}^G \mathbf{Q}$; thus,  
 $$[f_x]: q \mapsto g_x f(q)$$
 where here, on the right hand side, we regard $f$ as a map $\mathbf{Q} \rightarrow \mathbf{M}$ by Frobenius reciprocity (i.e. we pull it back via the obvious embedding $\mathbf{Q} \hookrightarrow \mathbf{P},q \mapsto 1 \otimes q$). 
 This concludes the justification of the subclaim.  \qed

  \subsection{Multiplication in the derived Hecke algebra}  \label{sec:mult}
  We finally analyze composition (i.e. multiplication) in the derived Hecke algebra, explicating
  it with respect to the isomorphism \eqref{serre}, and thus justifying the description given in \S \ref{desc2}.

 Let $\alpha, \beta, \gamma \in G/K$ with representatives $g_{\alpha}, g_{\beta}, g_{\gamma} \in G$. 
 Suppose given
 $h_{\alpha} \in H^*(K_{g_{\alpha}})$ and similarly for $\beta$.   We will compute  the product $h_{\beta} h_{\alpha}$ considered as elements of the derived Hecke algebra -- or more precisely the $H^*(K_{g_{\gamma}})$ component of their product. 
 
As in \eqref{dha-additive} we represent $h_{\alpha}$
 by an element  $h_{\alpha}' \in 
 \Hom_{K}(\mathbf{Q}, \ind_K^G \mathbf{Q})$,
 and denote by $[h_{\alpha}']$ the corresponding element of
 $
 \Hom_{K_{g_{\alpha}}}(\mathbf{Q}, \mathbf{Q}_{ g_{\alpha}})$
 and similarly for $h_{\beta}$.  
By \eqref{tpmodel} we have the explicit formula 
$$h_{\alpha}': q \in Q \mapsto 
 \sum_{k \in K/K_{\alpha}} k g_{\alpha} \otimes [h_{\alpha}'](k^{-1} q) \in S[K g_{\alpha} K] \otimes_{SK} Q,$$
 where we make a modest abuse of notation by identifying $K/K_{\alpha}$ to a set of representatives for it in $K$. 
 Now apply $h_{\beta}'$ to the right-hand side, regarding $h_{\beta}' \in \Hom_{SG}(\mathbf{P}, \mathbf{P})$. The result is:  
 \begin{equation} \label{bigmess} 
 \sum_{k \in K/K_{\alpha}}  \sum_{k' \in K/K_{\beta}} k g_{\alpha} k' g_{\beta} \otimes  [h'_{\beta}] k'^{-1} [h'_{\alpha}] k^{-1}  
 \in \Hom_{K} (\mathbf{Q}, SG \otimes_{SK} \mathbf{Q})\end{equation}
The desired  
 $H^*(K_{g_{\gamma}})$ component of the product $h_{\beta} \cdot h_{\alpha}$  is given by 
 considering all $k,k'$ for which $k g_{\alpha} k' g_{\beta} \in K g_{\gamma} K$, i.e. it is represented by 
 \begin{equation} \label{bigmess2} 
 \sum_{k g_{\alpha} k' g_{\beta} \in K g_{\gamma} K } k g_{\alpha} k' g_{\beta} \otimes  [h'_{\beta}] k'^{-1} [h'_{\alpha}] k^{-1}  
 \in \Hom_{K} (\mathbf{Q}, S[K g_{\gamma} K] \otimes_{SK} \mathbf{Q})\end{equation}
By ``dual'' Frobenius reciprocity (see before
\eqref{tpmodel})
the right-hand side  can be identified with $\Hom_{K}(\mathbf{Q}, \ind^K_{K_{\gamma}} \mathbf{Q}_{\gamma}) \simeq \Hom_{K_{\gamma}}(\mathbf{Q}, \mathbf{Q}_{\gamma})$.
If we write $kg_{\alpha} k' g_{\beta} = k_1 g_{\gamma} k_2$, an explicit formula for the corresponding
   element of  $\Hom_{K}(\mathbf{Q}, \ind^K_{K_{\gamma}} \mathbf{Q}_{\gamma})$ is given \eqref{IR1} by
$$q \mapsto  \sum_{k g_{\alpha} k' g_{\beta} = k_1 g_{\gamma} k_2 }  \underbrace{k_1 \otimes k_2 [h'_{\beta}] k'^{-1} [h'_{\alpha}] k^{-1} (q) }_{\in SK \otimes_{S K_{\gamma}} \mathbf{Q}_{\gamma}}
 $$
 where the right-hand sum is over the same $k, k'$ as before, and we consider only those $k,k'$ such that $kg_{\alpha} k' g_{\beta}  \in K g_{\gamma} K$. 
Then the corresponding  element of $\Hom_{K_{\gamma}}(\mathbf{Q}, \mathbf{Q}_{\gamma})$
 is given (see \eqref{tpmodel})  by picking out those terms for which $k_1 \in K_{\gamma}$; in that case we 
 can rewrite $k_1 g_{\gamma} k_2 = g_{\gamma} (\Ad(g_{\gamma})^{-1} k_1) k_2$ and so 
we may as well suppose that $k_1 = 1$. Thus, 
 the desired result is 
\begin{equation} \label{finale}  \sum_{k, k': k g_{\alpha} k' g_{\beta} =  g_{\gamma} k''}  k''   \underbrace{[h'_{\beta}] k'^{-1}}_{\Hom_{K_{\Ad(k g_{\alpha}) K_{k' g_{\beta}}}}(\mathbf{Q}_{k g_{\alpha}}, \mathbf{Q}_{k g_{\alpha} k' g_{\beta}}) }  \underbrace{ [h'_{\alpha}] k^{-1}}_{\Hom_{K_{kg_{\alpha}}}(\mathbf{Q}, \mathbf{Q}_{k g_{\alpha}}) }\in \Hom_{K_{\gamma}}(\mathbf{Q},\mathbf{Q}_{\gamma}) \end{equation} 
 
Here we observed 
 that  {\small $$[h'_{\alpha}] \circ k^{-1} \in \Hom_{K_{k g_{\alpha}}}(\mathbf{Q}, \mathbf{Q}_{k g_{\alpha}}),
 [h'_{\beta}] \circ k'^{-1} \in \Hom_{K_{k'g_\beta}}(\mathbf{Q}, \mathbf{Q}_{k'g_{\beta}}) 
\simeq  \Hom_{\Ad(g) K_{k' g_{\beta}}} (\mathbf{Q}_{g}, \mathbf{Q}_{g  k'g_{\beta}}) $$}
(the last isomorphism is the obvious one and we apply it to $g =kg_{\alpha}$). \footnote{For example, to check the first,  
note that for  $z \in K_{k g_{\alpha}}$  and $q \in \mathbf{Q}$ we have
 $$ [h'_{\alpha}]  \circ k^{-1} (zq) = [h'_{\alpha}]   \left( (k^{-1} zk) (k^{-1} q) \right)  = (g_{\alpha}^{-1} k^{-1} z k g_{\alpha}) [h'_{\alpha}]\circ k^{-1}( q)$$
 Indeed, 
 $[h'_{\alpha}]$ represents the class in $H^*(K_{k g_{\alpha}})$, obtained by applying $\Ad(k)$ to the class $h_{\alpha} \in H^*(K_{g_{\alpha}})$.}
Returning to \eqref{finale}, set
$$ x = k g_{\alpha} k' g_{\beta} K = g_{\gamma} K, \ \ y = k g_{\alpha} K, \ \ z = eK, \ \ U = \mathrm{stabilizer}(x, y, z)$$
Then $x,y$ are in relative position $\beta$, and  $y, z$ are in relative position $\alpha$, and $x,z$ are in relative position $\gamma$. 
 
Note also that  $$U= K_{k g_{\alpha}} \cap \Ad(k g_{\alpha}) K_{k' g_{\beta}}  =K \cap \Ad(k g_{\alpha}) K \cap \Ad(k g_{\alpha} k' g_{\beta}) K.$$ 
Therefore, the composite  occurring in \eqref{finale} 
$$    F = \overbrace{[h'_{\beta}] k'^{-1}}^ {\Hom_{K_{\Ad(k g_{\alpha}) K_{k' g_{\beta}}}}(\mathbf{Q}_{k g_{\alpha}}, \mathbf{Q}_{k g_{\alpha} k' g_{\beta}}) }  \circ  \underbrace{ [h'_{\alpha}] k^{-1}}_{\Hom_{K_{g_{\alpha}}}(\mathbf{Q}, \mathbf{Q}_{k g_{\alpha}}) }$$
actually belongs to  $ \Hom_U(\mathbf{Q}, \mathbf{Q}_{k g_{\alpha} k' g_{\beta}}) $; as such $F$ defines a cohomology class for $U$. 
This cohomology class is given by  taking the classes $h_{\alpha}, h_{\beta}$, transporting them to classes in $H^*(K_{k g_{\alpha}})$ and
 $H^*(\Ad(k g_{\alpha}) K'_{k' g_{\beta}})$,
by means of $\Ad(k): K_{g_{\alpha}} \simeq K_{k g_{\alpha}}$  and $\Ad(k g_{\alpha} k'): K_{\beta} \simeq \Ad(k g_{\alpha}) K'_{k' g_{\beta}}$, 
restricting to $U$, and taking the cup product. 

Said differently, let us think of $h_{\alpha}$ as a $G$-invariant association $H_{\alpha}$
from pairs $(u,v) \in G/K \times G/K$ to cohomology classes in $H^*(G_{uv})$ --
the one whose value at $(g_{\alpha} K, e)$ is given by the original cohomology class in $H^*(K_{g_{\alpha}})$. 
Similarly for $h_{\beta}$. Then, 
$$ \mbox{ cohomology class of $F$ } = H_{\beta}(x,y)  \cup H_{\alpha}(y, z) \in H^*(G_{xyz})$$

  Now $K_{\gamma} = G_{xz}$ acts   on the set of $(k, k', k'')$ as above,  i.e. satisfying $k g_{\alpha} k' g_{\beta} = g_{\gamma} k''$, via $$ \kappa: (k, k', k'') \mapsto (\kappa k, k', \Ad(g_{\gamma}^{-1}) \kappa k'').$$ 
For fixed $(k, k', k'')$ the stabilizer of this $K_{\gamma}$-action is just $U$. 
The contribution of a single $K_{\gamma}$-orbit  is given by
$$ \sum_{\kappa \in K_{\gamma}/U} \Ad(g_{\gamma}^{-1})  \kappa \circ F \circ  \kappa^{-1} $$
which is to say that it averages $F$, considered as an element
of $\Hom(\mathbf{Q}, \mathbf{Q}_{\gamma})$, 
over the cosets of $K_{\gamma}/U$.  (The $\Ad g_{\gamma}^{-1}$ accounts for the twisted action on $\mathbf{Q}_{\gamma}$).
This precisely realizes
the  corestriction from $U$ to $K_{\gamma}$.

In summary,  we have recovered the description of multiplication in the derived Hecke algebra given in \S \ref{desc2}.

\section{  Koszul algebra; other odds and ends}\label{AppendixB}

Let $B$ be a commutative ring with $1$. 
Let $$B[[x_1, \dots, x_r]] = S \stackrel{\iota}{\longrightarrow} R = B[[y_1, \dots, y_{r-\delta}]]$$
where $\iota$ sends $x_i$ to $y_i$ for $i \leq r-\delta$,
and $x_i$ to zero for $i > r-\delta$. 

Let $\mathfrak{p}_S$ be the kernel of the natural augmentation $S \rightarrow B$, and similarly for $R$. 
 Write 
 $\mathfrak{t}_S$ for the $B$-linear dual of
for $\mathfrak{p}_S/\mathfrak{p}_S^2$ and similarly for $\mathfrak{t}_R$. 
Just as in  \eqref{coocoo} we have a canonical identification $$\Ext_S^1(B, B) \simeq \mathfrak{t}_S.$$

 We will prove:  

\begin{lemma} \label{KoszulExtComputation}
  $\Ext_S^*(B, B)$ is a free exterior algebra 
  over its degree $1$ component; thus we have $\Ext_S^*(B,B) \simeq \wedge^* \mathfrak{t}_S $
  as graded $B$-algebra.
Moreover, there  is an identification of $\Ext_S^*(R, B)$ with $\wedge^* (\mathfrak{t}_S/\mathfrak{t}_R)$
in such a way that the natural action of $\Ext_S^*(B, B) \simeq \wedge^* \mathfrak{t}_S$ is the natural one
obtained from the algebra map $\wedge^* \mathfrak{t}_S \rightarrow \wedge^* (\mathfrak{t}_S/\mathfrak{t}_R)$. 
\end{lemma}
This will follow from the computations of \S \ref{KA1} (more precisely, with  the precisely analogous computations wherein one replaces the role of symmetric algebras by power series algebras).

\subsection{Koszul algebra}  \label{KA1} 
  Let $W$ be a free module of rank $e$ over a base ring   
  $B$    and consider the ring $ R= \mathrm{Sym}_B(W)$, i.e. ``the ring of functions on $W^{\vee}$.''     In what follows
  we will omit the $B$ subscript on $\mathrm{Sym}_B$. 
  
   We have a resolution 
  $$\underbrace{ \dots \rightarrow \Sym(W)  \otimes \wedge^2 W \rightarrow  \Sym(W) \otimes W  \rightarrow \Sym(W)}_{\mathbf{K}} \rightarrow  B,$$
where the differential sends 
  $$ r \otimes w_1 \wedge \dots \wedge w_r  \in \Sym(W) \otimes \wedge^i W  \mapsto \sum_{i} (-1)^{i-1}  r w_i \otimes w_1\wedge \dots \widehat{w_i} \wedge \dots w_r .$$
 There is a corresponding resolution where we replace
 $\Sym(W)$ by its completion with respect to the augmentation $\Sym(W) \rightarrow B$, i.e.
 when we replace a symmetric algebra by a formal power series algebra. 
 
 In particular, we get
 $$\Hom_{R} (\mathbf{K}, B) \simeq  \left( \wedge^i W \right)^{\vee} \mbox{ with zero differentials.}$$
 and thus an identification of  $\Ext^i_R(B, B)$ with $ \left( \wedge^i W \right)^{\vee} $. 
 
 In fact, $\Ext^*_R(B, B)$ is a free exterior algebra, where  the algebra structure on the $\Ext$-groups arising from their identification with the cohomology of the differential graded algebra 
  $$\Hom_{R} (\mathbf{K}, \mathbf{K}).$$
  To see this, one verifies 
    that  each element of $w \in W^{\vee}$, considered
  as acting on $\mathbf{K}$ by contractions, actually defines a degree $-1$  endomorphism of $\mathbf{K}$;
   the resulting inclusion 
  $$\bigwedge \nolimits^{*}(W^{\vee}) \hookrightarrow \Hom_R(\mathbf{K}, \mathbf{K})$$ 
  gives a quasi-isomorphism of differential graded algebras. 

%
%
%
%
%
%
%
%

  Suppose now that $U$ is a free submodule of $W$ such that $W/U$ is also free.  In this situation we have a quotient map $$R = \Sym(W) \rightarrow \Sym(W/U) := \bar{R}$$ 
  We have a resolution of $R$-modules  (where the differential is given by the same formula as before): 
      $$\underbrace{ \dots \rightarrow \Sym(W)  \otimes \wedge^2 U \rightarrow  \Sym(W) \otimes U \rightarrow  \Sym(W)}_{\mathbf{Q}} \stackrel{\sim}{\longrightarrow} \bar{R},$$
and from this we  identify $\Ext_R^*(\bar{R}, B)$ with $(\wedge^* U)^{\vee}$. 

\begin{lemma}
The action of $\Ext^*_R(B, B) \simeq (\wedge^* W)^{\vee}$ on this is the natural one that arises from the map $W^{\vee} \rightarrow U^{\vee}$. 
\end{lemma}

\proof 
   It is enough to check this for the action of $\Ext^1_R(B, B)$.   
 We have identifications:
   $$ \Ext^*_R(\bar{R}, B) \simeq H^*\left(  \Hom_{R}(\mathbf{Q}, B) \right)  =  H^* \left( \Hom_R(\mathbf{Q}, \mathbf{K}) \right). $$
Explicitly, a class in $ \omega_j \in (\wedge^j U)^{\vee} \simeq \Ext^j_R(\bar{R}, B)$  
is represented by a map of complexes $\mathbf{Q} \rightarrow \mathbf{K}$ as follows:    
 \begin{equation}
 \xymatrix{
\cdots \ar[r]   &  	 R \otimes \wedge^{j+1}U \ar[r] \ar[d]^f 	&   R \otimes  \wedge^j U  \ar[r] \ar[d]^{\omega_j \in \wedge^j U^{\vee}} 	&   R \otimes  \wedge^{j-1} U   \\
\cdots  \ar[r]   &		 K_1=R \otimes W \ar[r] 				& K_0=R				\ar[r] 									& 0. 
 }
\end{equation}
(since this lifts 
  the map  $\mathbf{Q} \rightarrow B[j]$ associated to  $\omega_j$).
  
Fix a basis $e_1, \dots, e_r$ for $U$ and extend it to a basis $e_1, \dots, e_{e}$ for $W$.
For $I \subset \{1, \dots, e\}$ with cardinality $r$, we define $e_I \in \wedge^r W$ thus: 
write  $I = \{i_1, \dots, i_r\}$ with $i_1 < \cdots < i_r$
and put $e_I  = e_{i_1} \wedge e_{i_2} \wedge \dots \wedge e_{i_r}$. 
  We may choose $f$ to be given, explicitly, as  
$$ e_J \in \wedge^{j+1} U \mapsto \sum_{k \in J}  (-1)^{[k]-1} \omega_j(e_{J-k})   \otimes e_k \in R \otimes W,$$ 
where $[k]$ means the position of $k$ in $J$ (i.e.  if $J$ is ordered in increasing order, then $1$ for the smallest element, two for the second smallest, etc.)

To compute, now, the action of $\beta \in W^{\vee} =\Ext^1_R(B, B) $
on the class $\omega_j$, we just 
regard $\beta$ as an $R$-module map $K_1=  R \otimes W \rightarrow R$, 
and then compose 
 $$ \beta \circ f \in \Hom_R(\wedge^{j+1} U \otimes R, R) = \left( \wedge^{j+1} U \right)^{\vee}.$$
 Explicitly,  
   $$\beta \circ f:  e_J \mapsto \sum_{k \in J} (-1)^{[k]-1} \overline{ \beta}(e_k) \omega_j(e_{J-k})  = (\overline{\beta} \wedge \omega_j, e_J)$$
  i.e. $ \overline{\beta} \wedge \omega_j$, where $\overline{\beta}$ is the image of $\beta$ in $U^{\vee}$.  
  This concludes the proof.
  \qed

   \subsection{A result in topology} \label{sec:remedialtopology}  
 Suppose that  $\pi: X \rightarrow Y$ is a covering of pointed Hausdorff topological spaces, with Galois group $\Delta$. 
This covering is classified by a map $Y \rightarrow B \Delta$ from $Y$ to the classifying space of $\Delta$. 

There are  two natural actions of $H^*(\Delta, E)$ (with $E$ a coefficient ring)
on $H^*(Y, E)$:

\begin{itemize}
\item[(a)] The first arises from pullback of cohomology classes under $Y \rightarrow B\Delta$
together with cup product. 
\item[(b)] The second arises from the identification of the cochain complex of $Y$, with $E$ coefficients,  with 
$$C^*(Y; E) \simeq \Hom_{E \Delta}(C_*(X, E); E)$$
where $C_*(X; E)$ is the cochain complex of $X$ (or e.g. the complex of a $\Delta$-equivariant cell structure),  thought of as a complex of $E \Delta$-modules.
Then one composes with self maps of $E$ in the derived category of $E \Delta$-modules. 
 \end{itemize}

For lack of a reference, we will prove the coincidence of these actions. For this we will use the following standard Lemma
concerning the coincidence of singular and sheaf cohomology (see \cite{Sella} for a careful discussion; however this reference does not discuss the product structures): 
\begin{lemma} For any locally contractible Hausdorff space $M$,  and any $E$-module $A$, 
let $\underline{A}$
be the constant sheaf on $M$ with constant value $A$, considered
as an object of the category of sheaves $\mathcal{S}$ of $E$-modules on $M$. 

Then the complex of local chains $U \mapsto C^*(U, A)$ defines a presheaf on $M$;
let $\mathcal{C}^*_A$ be its sheafification.  Then $\underline{A} \rightarrow \mathcal{C}^*_A$ is a flasque resolution of $\underline{A}$. Moreover, the natural maps
$$ C^*(M, A) \rightarrow \Gamma(M, \mathcal{C}^*_A) = 
 \Hom_{\mathcal{S}}(\underline{A},  \mathcal{C}^*_A) $$ 
induces, at the level of cohomology, an isomorphism
\begin{equation} \label{oink boink} H^*(M, \underline{A}) \simeq   \Ext_{\mathcal{S}}^*(\underline{A}, \underline{A})\end{equation}
which carries the cup product on the left to the $\Ext$-product to the right. 
\end{lemma}
 
\proof (that (a) and (b) coincide): 
Observe, first of all,  that every $E\Delta$-module $M$ gives a locally constant sheaf $\underline{M}$ on $Y$, namely, the one represented by the covering $(X \times M)/\Delta \rightarrow Y$. 
The cochains $C^*(Y, \underline{M})$ are then given by $\Hom_{E \Delta}(C_*(X, E), M)$.

 Fix   $\alpha_{\Delta} \in H^m(\Delta, E)$. It gives rise to a homomorphism $\alpha: E \rightarrow E[m]$ in the derived category of $E\Delta$-modules,
 which can be represented by  a diagram
$E \stackrel{\sim}{\leftarrow} P \rightarrow E[m]$ 
 where $P$ is a complex of projective $E\Delta$-modules.
 Thus we get 
a diagram of locally constant sheaves on $Y$:
 $$\underline{\alpha}: \underline{E}  \stackrel{\sim}{\leftarrow} \underline{P} \rightarrow \underline{E}[m].$$

 This gives a map in the derived category of sheaves on $Y$, and thus an element of $\Ext^m_{\mathcal{S}}(\underline{E}, \underline{E})$; this  element represents the pullback $\alpha_Y$ of $\alpha_{\Delta}$ to $Y$. %
 or rather its image under  \eqref{oink boink}.

By the final sentence of the Lemma, the cup product with $\alpha_Y$ is given, at the level of cohomology, by the  $\Ext$-product, which is explicitly
the composite:
 $$H^*(Y, \underline{E}) \stackrel{\sim}{ \leftarrow} H^*(Y, \underline{P}) \rightarrow H^*(Y, \underline{E}[m])$$
 By the lemma, these groups are naturally identified with the cohomology of the corresponding 
 cochain groups; so the above composite coincides with $$\Hom_{E \Delta} (C_*(X, E), E) \stackrel{\sim}{\leftarrow} \Hom_{E \Delta}(C_*(X, E), P) \rightarrow \Hom_{E \Delta}(C_*(X, E), E)$$
 where the middle term is now the $\mathrm{Hom}$-complex between two complexes. 
 
But this composite is also given by the $\Ext$-product, in the category of $E \Delta$-modules, with the class of
$\alpha$. This concludes the proof of the coincidence of (a) and (b).  \qed
   
 \printindex

  \bibliography{derivedHecke}
  \bibliographystyle{plain}
   
  \end{document}